\documentclass[12pt]{article}
\usepackage{amssymb,amsfonts,amsmath,curves}
\usepackage{epsfig}
\usepackage{graphicx}
\usepackage{float}
\usepackage{color}
\usepackage{anysize}
\usepackage{url}
\usepackage[mathscr]{euscript}

\usepackage{dsfont}

\usepackage{psfrag}
\usepackage{hyperref}
\begin{document}

\newtheorem{defin}{Definition}[section]
\newtheorem{Prop}{Proposition}
\newtheorem{theorem}{Theorem}[section]
\newtheorem{conjecture}{Conjecture}[section]
\newtheorem{proposition}{Proposition}[section]
\newtheorem{lemma}{Lemma}[section]
\newtheorem{ml}{Main Lemma}
\newtheorem{con}{Conjecture}
\newtheorem{cond}{Condition}
\newtheorem{conj}{Conjecture}
\newtheorem{prop}[theorem]{Proposition}
\newtheorem{lem}{Lemma}[section]
\newtheorem{rmk}[theorem]{Remark}
\newtheorem{corollary}{Corollary}[section]
\renewcommand{\theequation}{\thesection .\arabic{equation}}
\renewcommand{\labelenumi}{(\alph{enumi})}

\newcommand{\beq}{\begin{equation}}
\newcommand{\eeq}{\end{equation}}
\newcommand{\beqn}{\begin{eqnarray}}
\newcommand{\beqnn}{\begin{eqnarray*}}
\newcommand{\eeqn}{\end{eqnarray}}
\newcommand{\eeqnn}{\end{eqnarray*}}
\newcommand{\bprop}{\begin{prop}}
\newcommand{\eprop}{\end{prop}}
\newcommand{\bteo}{\begin{teo}}
\newcommand{\bcor}{\begin{cor}}
\newcommand{\ecor}{\end{cor}}
\newcommand{\bcon}{\begin{con}}
\newcommand{\econ}{\end{con}}
\newcommand{\bcond}{\begin{cond}}
\newcommand{\econd}{\end{cond}}
\newcommand{\bconj}{\begin{conj}}
\newcommand{\econj}{\end{conj}}
\newcommand{\eteo}{\end{teo}}
\newcommand{\brm}{\begin{rmk}}
\newcommand{\erm}{\end{rmk}}
\newcommand{\blem}{\begin{lem}}
\newcommand{\elem}{\end{lem}}
\newcommand{\ben}{\begin{enumerate}}
\newcommand{\een}{\end{enumerate}}
\newcommand{\bei}{\begin{itemize}}
\newcommand{\eei}{\end{itemize}}
\newcommand{\bdf}{\begin{defin}}
\newcommand{\edf}{\end{defin}}
\newcommand{\bpr}{\begin{proof}}
\newcommand{\epr}{\end{proof}}

\newenvironment{proof}{\noindent {\em Proof}.\,\,}{\hspace*{\fill}$\halmos$\medskip}

\newcommand{\halmos}{\rule{1ex}{1.4ex}}
\def \qed {{\hspace*{\fill}$\halmos$\medskip}}

\newcommand{\fr}{\frac}
\newcommand{\Z}{{\mathbb Z}}
\newcommand{\R}{{\mathbb R}}
\newcommand{\E}{{\mathbb E}}
\newcommand{\C}{{\mathbb C}}
\renewcommand{\P}{{\mathbb P}}
\newcommand{\N}{{\mathbb N}}
\newcommand{\var}{{\mathbb V}}
\renewcommand{\S}{{\cal S}}
\newcommand{\T}{{\cal T}}
\newcommand{\W}{{\cal W}}
\newcommand{\X}{{\cal X}}
\newcommand{\Y}{{\cal Y}}
\newcommand{\h}{{\cal H}}
\newcommand{\f}{{\cal F}}

\renewcommand{\a}{\alpha}
\renewcommand{\b}{\beta}
\newcommand{\g}{\gamma}
\newcommand{\G}{\Gamma}
\renewcommand{\L}{\Lambda}
\renewcommand{\l}{\lambda}
\renewcommand{\d}{\delta}
\newcommand{\D}{\Delta}
\newcommand{\e}{\epsilon}
\newcommand{\s}{\sigma}
\newcommand{\B}{{\cal B}}
\renewcommand{\o}{\omega}

\newcommand{\nn}{\nonumber}
\renewcommand{\=}{&=&}
\renewcommand{\>}{&>&}
\newcommand{\<}{&<&}
\renewcommand{\le}{\leq}
\newcommand{\+}{&+&}

%Differentiation
\newcommand{\pa}{\partial}
\newcommand{\ffrac}[2]{{\textstyle\frac{{#1}}{{#2}}}}
\newcommand{\dif}[1]{\ffrac{\partial}{\partial{#1}}}
\newcommand{\diff}[1]{\ffrac{\partial^2}{{\partial{#1}}^2}}
\newcommand{\difif}[2]{\ffrac{\partial^2}{\partial{#1}\partial{#2}}}

\definecolor{Red}{rgb}{1,0,0}
\def\red{\color{Red}}

\title{Ellipticity criteria for ballistic behavior of random walks
in random environment}
\author{David Campos$^1$\thanks{
Partially supported by fellowship
Consejo Nacional de Ciencia y Tecnolog\'\i a number D-$57080025$,
by Mecesup PUC $0802$
and by  Universidad Costa Rica Oficio No OAICE-$06$-CAB-$058$-$2008$.}
\ \
 and Alejandro F.\ Ram\'irez$^2$\thanks{
Partially supported by Fondo Nacional de Desarrollo Cient\'\i fico
y Tecnol\'ogico grant $1100746$.}}

\maketitle
%\tableofcontents

{\footnotesize

\thanks{$^{1}$ Facultad de Matem\'{a}ticas, Pontificia Universidad Cat\'{o}lica de Chile, Vicu\~{n}a Mackenna 4860, Macul, Santiago, Chile;\\
{\it e-mail:}
 \href{mailto:jdcampo1@mat.puc.cl}{jdcampo1@mat.puc.cl}

\thanks{$^{2}$ Facultad de Matem\'{a}ticas, Pontificia Universidad Cat\'{o}lica de Chile, Vicu\~{n}a Mackenna 4860, Macul, Santiago, Chile;\\
{\it e-mail:}
 \href{mailto:aramirez@mat.puc.cl}{aramirez@mat.puc.cl}

}

\begin{abstract}
We introduce  ellipticity criteria
for random walks in i.i.d. random environments
under which we can extend the ballisticity
conditions of Sznitman's and the polynomial effective
criteria of Berger, Drewitz and Ram\'\i rez
originally defined for uniformly elliptic random walks.
We  prove under them
the equivalence of Sznitman's
$(T')$ condition with the polynomial
effective criterion $(P)_M$, for $M$ large enough.
We furthermore give  ellipticity criteria under
which a random walk satisfying the polynomial effective
criterion, is ballistic, satisfies the annealed
central limit theorem or the quenched central limit theorem.
\end{abstract}

{\footnotesize
{\it 2000 Mathematics Subject Classification.} 60K37, 82D30.

{\it Keywords.} Random walk in random environment, renormalization,
uniform ellipticity, ellipticity.
}
%%%%%%%%%%%%%%%%%%%%%%%%%%%%%%%%%%%%%%%%%%%%%%%%%%%%%%%%%%%%%%%%%%%%%%%%%%%%%%%%%%%%%%%%%%%%%%%%%%%%%%%%%%%%%%%%%%%%%%%%%%%%%%%%%%%%%%%%%%
\section {Introduction}

We introduce  ellipticity criteria for random walks
in random environment which enable us to extend to
environments which are not necessarily uniformly elliptic the
ballisticity conditions for the uniform elliptic case of Sznitman \cite{Sz02} and
of Berger, Drewitz and Ram\'\i rez \cite{BDR12}, their equivalences
and some of their consequences \cite{SZ99,Sz01,Sz02,RAS09, BZ08}.

 For $x\in\mathbb R^d$, denote by $|x|_1$
and $|x|_2$
its $l_1$ and $l_2$ norm respectively. Call $U:=\{e\in\mathbb Z^d:|e|_1=1\}=\{e_1,\ldots,e_{2d}\}$ the canonical vectors
with the convention that $e_{d+i}=-e_i$ for $1\le
i\le d$ and
let ${\mathcal P}:=\{p(e): p(e)\ge 0,\sum_{e\in U}p(e)=1\}$.
 An \textit{environment} is an element $\omega$ of the \textit{environment space} $\Omega:={\cal{P}}^{\Z^d}$ so
that $\omega:=\{\omega(x):x \in \Z^d\}$, where $\omega(x) \in \cal{P}$. We
denote the components of $\omega(x)$ by $\omega(x,e)$.
 The \textit{random walk in the environment $\omega$ starting from $x$}
 is  the Markov
 chain $\{X_n: n \geq 0\}$ in $\Z^d$ with law $P_{x,\omega}$
defined by the condition $P_{x,\omega}(X_0=x)=1$ and the
transition probabilities

$$
P_{x,\omega}(X_{n+1}=x+e|X_n=x)=\omega (x,e)
$$
for each $x \in \Z^d$ and $e\in U$.
Let $\P$ be a probability measure defined on the
environment space $\Omega$ endowed with its Borel $\sigma$-algebra.
We will assume  that $\{\omega(x):x\in\mathbb Z^d\}$ are i.i.d. under $\P$.
We will call $P_{x,\omega}$ the {\it quenched law} of the random walk in
random environment (RWRE) starting from $x$, while $P_x:=\int P_{x,\omega}d\P$
the {\it averaged} or {\it annealed law} of the RWRE starting from $x$.

 We say that the law $\P$  of the RWRE
is {\it elliptic} if for every $x\in\mathbb Z^d$ and $e\in U$
one has that $\P(\omega(x,e)>0)=1$. We say that $\P$ is
{\it uniformly elliptic} if there exists a constant $\kappa>0$
such that for every $x\in\mathbb Z^d$ and $e\in U$
it is true that $\P(\omega(x,e)\ge \kappa)=1$.
Given $l\in\mathbb S^{d-1}$ we say that the RWRE is
{\it transient in direction $l$} if

$$
P_0(A_l)=1,
$$
where

\begin{equation}
\nonumber
A_l:=\{\lim_{n\to\infty}X_n\cdot l=\infty\}
\end{equation}
We say that it is {\it ballistic in direction $l$}
if $P_0$-a.s.

$$
\liminf_{n\to\infty}\frac{X_n\cdot l}{n}>0.
$$
The following is conjectured (see for example \cite{Sz04}).

\medskip

\begin{conjecture} Let $l\in\mathbb S^{d-1}$.
 Consider a random walk in a uniformly elliptic i.i.d. environment  in dimension
 $d\ge 2$, which is transient in direction $l$. Then it is ballistic
in direction $l$.
\end{conjecture}

\medskip

Some partial progress towards the resolution of this conjecture
has been made in \cite{Sz01,Sz02,DR11,DR12,BDR12}.
In 2001 and 2002 Sznitman in  \cite{Sz01,Sz02} introduced  a class
of ballisticity conditions under which he could prove the above statement.
For each subset $A\subset\mathbb Z^d$ define the first exit time
from the set $A$ as

\begin{equation}
\label{exit-time}
T_A:=\inf\{n\ge 0: X_n\notin A\}.
\end{equation}
For $L>0$ and $l\in\mathbb S^{d-1}$
 define the slab

\begin{equation}
\label{slab}
U_{l,L}:=\{x\in\mathbb Z^d:
-L\le x\cdot l\le L\}.
\end{equation}
Given $l\in\mathbb S^{d-1}$ and $\gamma\in (0,1)$, we
say that condition $(T)_\gamma$ in direction $l$ (also written
as $(T)_\gamma|l$) is satisfied if there exists a neighborhood
$V\subset\mathbb S^{d-1}$ of $l$ such that for all $l'\in V$

$$
\limsup_{L\to\infty}\frac{1}{L^\gamma}\log
P_0(X_{T_{U_{l',L}}}\cdot l'<0)<0.
$$
Condition $(T')|l$ is defined as the fulfillment of condition
$(T)_\gamma|l$ for all $\gamma\in (0,1)$. Sznitman \cite{Sz02} proved
that if a random walk in an i.i.d. uniformly elliptic environment
satisfies $(T')|l$ then it is ballistic in direction $l$.
He also showed that if $\gamma\in (0.5,1)$, then
$(T)_\gamma$ implies $(T')$. In 2011, Drewitz and Ram\'\i rez \cite{DR11}
showed that there is a $\gamma_d\in (0.37,0.39)$ such that
if $\gamma\in (\gamma_d,1)$, then $(T)_\gamma$ implies $(T')$.
In  2012, in \cite{DR12}, they were able to show that for
dimensions $d\ge 4$, if $\gamma\in (0,1)$, then $(T)_\gamma$
implies $(T')$.
 Recently in \cite{BDR12}, Berger, Drewitz
and Ram\'\i rez introduced a polynomial ballisticity condition,
weakening further the conditions $(T)_\gamma$. The condition
is effective, in the sense that it can a priori be verified
explicitly for
a given environment.
To define it,  for each $L,\tilde L>0$ and $l\in\mathbb S^{d-1}$
consider the box

\begin{equation}
\label{boxbox}
B_{l,L,\tilde L}:=R\left(\left(-L,L\right)\times\left(-\tilde L,\tilde L\right)^{d-1}
\right)\cap \mathbb Z^d,
\end{equation}
where
$R$ is a rotation of $\mathbb R^d$
defined by the  condition

\begin{equation}
\label{defr}
R(e_1)=l.
\end{equation}
Given $M\ge 1$ and $L\ge 2$, we say that the {\it polynomial
condition $(P)_M$ in direction $l$ is satisfied 
on a box of size $L$} (also
written as $(P)_M|l$) if
there exists  an $\tilde L\le 70 L^3$ 
such that
the following upper bound for the probability that the walk
does not exit the box $B_{l,L,\tilde L}$ through its front side
is satisfied

$$
P_0(X_{T_{B_{l,L,\tilde L}}}\cdot l< L)\le\frac{1}{L^M}.
$$
In \cite{BDR12}, Berger, Drewitz and Ram\'\i rez prove that
every  random walk in an i.i.d. uniformly elliptic environment which satisfies
$(P)_M$ for $M\ge 15 d+5$ is necessarily ballistic.

On the other hand, it is known (see for example Sabot-Tournier \cite{ST11}) that
in dimension $d\ge 2$, there exist elliptic random walks
which are transient in a given direction but not ballistic
in that direction.  The purpose of this paper is to investigate to
which extent can the assumption of uniform ellipticity
 be weakened.  To do this we introduce several classes
of ellipticity conditions on the environment.
Define for each $\alpha\ge 0$,

\begin{equation}
\label{eta-alpha}
\eta_\alpha:=\max_{e\in U}\E\left[\frac{1}{\omega(0,e)^\alpha}
\right]
\end{equation}
and
$$
\bar\alpha:=\sup\left\{
\alpha\ge 0: \eta_\alpha<\infty\right\}.
$$
Let $\alpha\ge 0$.
We say that the law of the environment satisfies
the {\it ellipticity condition $(E)_\alpha$} if
for every $e\in U$ we have that

\begin{equation}
\nonumber
\bar\alpha>\alpha.
\end{equation}
Note that  $(E)_0$ is equivalent to
the existence of an $\alpha>0$ such that $\eta_\alpha<\infty$.
To state the first result of this paper, we need
to define the following constant

\begin{equation}
\label{deflo}
c_0:=\frac{2}{3}3^{120 d^4+3000 d(\log\eta_{\bar\alpha/2})^2}.
\end{equation}
Throughout the rest of this paper, whenever we assume the polynomial
condition $(P)_M$ is satisfied, it will be understood that this
happens on a box of size $L\ge c_0$.

\medskip

\begin{theorem}
\label{theorem1}
 Consider a random walk in an i.i.d. environment
in dimensions $d\ge 2$. 
 Let $l\in\mathbb S^{d-1}$ and $M\ge 15 d+5$.
Assume that the environment
satisfies the ellipticity condition $(E)_0$.
Then the polynomial condition $(P)_M|l$
 is equivalent to $(T')|l$.
\end{theorem}

\medskip

\noindent It is important to remark that we have not made any particular
effort to optimize the values of $15d+5$ and $c_0$ in the above theorem.

 In this paper we go further from Theorem \ref{theorem1}, and we obtain
assuming $(T')$,
good enough tail estimates for  the distribution of
the regeneration times of the random walk.
Let us recall that there exists an {\it asymptotic direction}
if the limit

$$
\hat v:=\lim_{n\to\infty}\frac{X_n}{|X_n|_2}
$$
exists $P_0$-a.s. The polynomial condition $(P)_M$ implies the
existence of an asymptotic direction (see for example
Simenhaus \cite{Sim07}). Whenever the asymptotic direction exists,
let us define the half space

$$
H_{\hat v}:=\{l\in\mathbb S^{d-1}:l\cdot\hat v\ge 0\}.
$$
Let $\alpha>0$.
We say that the law of the environment
satisfies  the {\it ellipticity condition $(E')_\alpha$}
if there exists an $\{\alpha(e):e\in U\}\in
(0,\infty)^{2d}$ such that

\medskip

\begin{equation}
\label{kappa}
\kappa=\kappa(\{\alpha(e)\}):=2\sum_{e}\alpha(e)-\sup_{e\in U}\left(
\alpha(e)+\alpha(-e)\right)>\alpha
\end{equation}
and \medskip
\begin{equation}
\label{finiteness}
\E\left[\prod_{e\in U}\frac{1}{\omega(0,e)^{\alpha(e)}}\right]<\infty {\red .}
\end{equation}
Note that $(E')_{\alpha'}$ implies $(E)_\alpha$ for

\begin{equation}
\label{eprime-eo}
\alpha<\min_e\alpha(e).
\end{equation}
Furthermore, we say that the ellipticity condition $(E')_\alpha$
is satisfied {\it  towards
the asymptotic direction} if  there exists an $\{\alpha(e):e\in U\}$
satisfying (\ref{kappa}) and (\ref{finiteness}) and such that

\begin{equation}
\label{direction}
 {\rm there\ exists}\ \alpha_1>0\ {\rm
   such\ that}\  \alpha(e)=\alpha_1\ {\rm for\  all}\ e\in H_{\hat
  v},\ {\rm while}\  \alpha(e)\le\alpha_1\ {\rm for}\ e\notin H_{\hat v}.
\end{equation}
\medskip

\noindent The second main result of this paper is the following theorem.

\medskip
\begin{theorem}
\label{theorem2} {\bf (Law of large numbers)}
Consider a random walk in an i.i.d. environment
in dimensions $d\ge 2$. Let $l\in\mathbb S^{d-1}$ and $M\ge 15 d+5$.
 Assume that the random walk
satisfies the polynomial condition $(P)_M|l$ and the ellipticity condition
 $(E')_1$ towards the asymptotic direction (cf. (\ref{kappa}), (\ref{finiteness}) and (\ref{direction})).
 Then the random walk is ballistic
in direction $l$ and there is a $v\in\mathbb R^d$,
$v\ne 0$ such that

$$
\lim_{n\to\infty}\frac{X_n}{n}=v,\quad P_0-a.s.
$$

\end{theorem}
\medskip

\noindent We have directly the following two corollaries,
the second one following from H\"older's inequality.

\medskip

\begin{corollary}
\label{corollary2-intro}
Consider a random walk in an i.i.d. environment
in dimensions $d\ge 2$. Let $l\in\mathbb S^{d-1}$ and $M\ge 15 d+5$. Assume that the random walk
satisfies condition $(P)_M|l$ and that there is an $\alpha>\frac{1}{4d-2}$ such that

$$
\E\left[\prod_{ e\in U}\frac{1}{\omega(0,e)^\alpha}\right]<\infty.
$$
 Then the random walk is ballistic
in direction $l$.
\end{corollary}
\medskip
\begin{corollary}
\label{corollary3-intro}
Consider a random walk in an i.i.d. environment
in dimensions $d\ge 2$. Let $l\in\mathbb S^{d-1}$ and $M\ge 15 d+5$. Assume that the random walk
satisfies condition $(P)_M|l$ and the ellipticity condition
 $(E)_{1/2}$.
 Then the random walk is ballistic
in direction $l$.

\end{corollary}
\medskip
\noindent
Theorems \ref{theorem1} and Theorem \ref{theorem2} give
ellipticity criteria for ballistic behavior for general random
walks in i.i.d. environments, which as we will see below,
should not be far from optimal criteria.
Indeed,
 the value $1$ of condition $(E')_1$ in Theorem \ref{theorem2}
is optimal: within the context of random walks in Dirichlet
random environments (RWDRE), it is well known that
there are examples of walks which satisfy $(E')_\alpha$ for
$\alpha$ smaller but arbitrarily close to $1$, towards
the asymptotic direction (cf. \ref{direction}), which
are transient but not ballistic in a given direction
(see \cite{Sa11a,Sa11b,ST11}). A RWDRE with parameters $\{\beta_e:e\in U\}$
is a random walk in
an i.i.d. environment whose law at a single site has the density

$$
\frac{\Gamma(\sum_{e\in U}\beta_e)}{\prod_{e\in U}\Gamma(\beta_e)}
\prod_{e\in U} p(e)^{\beta_e-1}
$$
with respect to the Lebesgue measure on $\mathcal P$.
 We can also
explicitly construct  examples as those mentioned above  in analogy to the random conductance
model (see for example Fribergh  \cite{F11}, where he
characterizes the ballistic random walks within the directional transient cases). In fact,
for every $\epsilon>0$, one can
construct an environment such that condition $(E')_{1-\epsilon}$ is satisfied
towards the asymptotic direction,
but the walk is transient in direction $e_1$  but  not ballistic
in direction $e_1$. Let $\phi$ be any random
variable taking values on the interval $(0,1/4)$ and such
that the expected value of $\phi^{-1/2}$ is infinite,
while for every $\epsilon>0$, the expected value of
$\phi^{-(1/2-\epsilon)}$ is finite.
Let $X$ be a Bernoulli random variable of parameter $1/2$.
We now define $\omega_1(0,e_1)=2\phi$, $\omega_1(0,-e_1)=\phi$,
$\omega_1(0,-e_2)=\phi$ and $\omega_1(0,e_2)=1-4\phi$
and
$\omega_2(0,e_1)=2\phi$, $\omega_2(0,-e_1)=\phi$,
$\omega_2(0,e_2)=\phi$ and $\omega_2(0,-e_2)=1-4\phi$.
We then let the environment at site $0$ be given by
the random variable $\omega(0):=1_X(1)\omega_1(0)
+1_X(0)\omega_2(0)$. This environment has the
property that traps can appear, where the random walk
gets caught in an edge, as shown shown in
Figure \ref{Fribergh} and it does satisfy $(E')_{1-\epsilon}$ towards
the asymptotic direction. Furthermore,  it is not
difficult to check that the random walk in this random
environment is transient
in direction $e_1$ but not ballistic.
It will be shown in a future work, that this environment satisfies
 the polynomial condition $(P)_M$ for
$M\ge 15d+5$.

\begin{figure}[H]
\centering
\includegraphics[width=7cm]{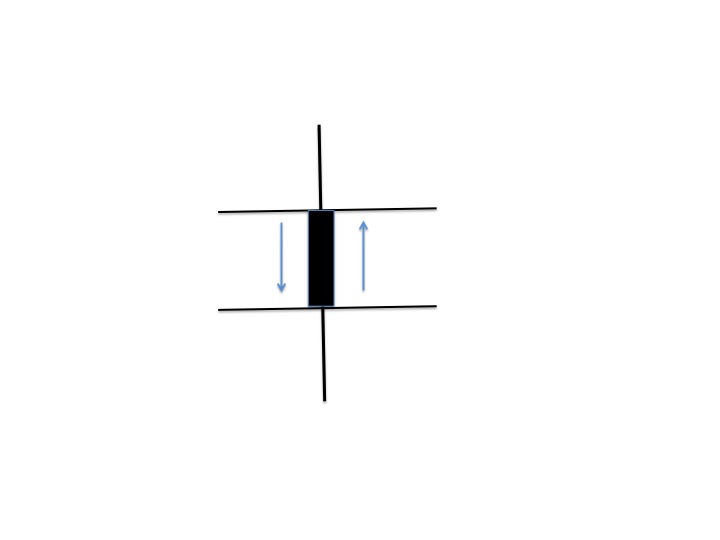}
\caption{A trap produced by an elliptic environment which does not
satisfy $(E')_{1}$.}
\label{Fribergh}
\end{figure}

As mentioned above, similar examples
of random walks in elliptic i.i.d. random environment
which are transient in a given direction but
not ballistic have been exhibited within the
context of the Dirichlet environment.
Here, the environment is chosen i.i.d. with
a Dirichlet distribution at each site
$D(\beta_1,\ldots,\beta_{2d})$ of parameters
$\beta_1,\ldots,\beta_{2d}>0$ (see for example \cite{Sa11a,Sa11b,ST11}),
the parameter  $\beta_i$ being associated with the
direction $e_i$.
For a random walk in Dirichlet random environment,
condition
$(E')_1$
 is equivalent
to

\begin{equation}
\label{lambda-st}
\lambda:=
2\sum_{j=1}^{2d}\beta_j-\sup_{i:1\le i\le 2d}(\beta_i+\beta_{i+d})>1.
\end{equation}
This is the characterization of ballisticity
given by Sabot in \cite{Sa11b} for random walks
in random Dirichlet environments in dimension $d\ge 3$.
Tournier in \cite{T11} proved that if $\lambda\le 1$, then
the RWDRE is not ballistic in any direction. Sabot in \cite{Sa11b},
showed that if $\lambda>1$, and if there is an $i=1,\ldots,d$
such that $\beta_i\ne \beta_{i+d}$, then the random walk
is ballistic. It is thus natural to wonder to what
general condition corresponds (not restricted to random Dirichlet environments),
 the characterization of Sabot and Tournier.
In section \ref{section2}, we will see that there
are several formulations of the necessary and sufficient condition
for ballisticity of Sabot and Tournier
for RWDRE  (cf. (\ref{lambda-st})), but which are not equivalent for general RWRE.
Among these formulations, the following one
is the weakest one in general.
We say that condition $(ES)$ is satisfied if

$$
\max_{i:1\le i\le d} \E\left[\frac{1}{1-\omega(0,e_i)\omega(e_i,-e_i)}
\right]<\infty.
$$

\noindent We have furthermore the
following proposition whose proof will be presented in
section \ref{section2}.
\medskip

\begin{proposition}
\label{proposition11} Consider a random walk in a random environment.
Assume that condition $(ES)$ is not satisfied. Then the random walk
is not ballistic.
\end{proposition}
\medskip

\noindent We will see in the proof of Proposition \ref{proposition11}
how important is the role played by  the edges depicted
in Figure \ref{Fribergh} which play
the role of traps.

\medskip

Another consequence of Theorem \ref{theorem1} and
the machinery that we develop to estimate
the tails of the regeneration times is the
following theorem.

\medskip

\begin{theorem}
\label{theorem3}
Consider a random walk in an i.i.d. environment
in dimensions $d\ge 2$. Let $l\in\mathbb S^{d-1}$ and $M\ge 15 d+5$. Assume that the random walk
satisfies condition $(P)_M|l$.

\begin{itemize}

\item[a)] {\bf (Annealed central limit theorem)} If $(E')_1$
 is satisfied towards the asymptotic direction then

$$
\epsilon^{1/2}(X_{[\epsilon^{-1}n]}-[\epsilon^{-1}n]v)
$$
converges in law under $P_0$ as $\epsilon\to 0$
to a Brownian motion with non-degenerate covariance matrix.

\item[b)] {\bf (Quenched central limit theorem)} If $(E')_{88d}$
 is satisfied towards the asymptotic direction,
 then $\P$-a.s. we have that

$$
\epsilon^{1/2}(X_{[\epsilon^{-1}n]}-[\epsilon^{-1}n]v)
$$
converges in law under $P_{0,\omega}$ as $\epsilon\to 0$
to a Brownian motion with non-degenerate covariance matrix.

\end{itemize}
\end{theorem}
\medskip
Part $(b)$ of the above Theorem is based on
a result of Rassoul-Agha and Sepp\"{a}l\"{a}inen \cite{RAS09},
which gives as a condition so that an elliptic
random walk satisfies the quenched central limit theorem
that the regeneration times have moments of order higher than
$176d$. As they point out in their paper, this particular
lower bound on the moment should not have any meaning
and it is likely that it could be improved. For example,
Berger and Zeitouni in \cite{BZ08}, also prove the
quenched central limit theorem under lower order
moments for the regeneration times but under the
assumption of uniform ellipticity. It should be possible
to extend their methods to elliptic random walks
in order to improve the moment condition of part $(b)$
of Theorem \ref{theorem3}.

It is possible using the methods developed in this paper,
to obtain slowdown large deviation estimates for
the position of the random walk in the spirit
of the estimates obtained by Sznitman in \cite{Sz00}
for the case where the environment is plain nestling
(see also \cite{Ber12}), under hypothesis of the form

$$
\E\left[e^{(\log\omega(0,e))^{\beta}}\right]<\infty,
$$
for some $\beta>1$. Nevertheless, the estimates we would
obtain would not be sharp, in the sense that we would
obtain an upper bound for the
probability that at time $n$ the random walk is slowed down
of the form

$$
e^{-(\log n)^{\beta'(d)}},
$$
where $\beta'(d)>1$, but $\beta'(d)<d$ (as discussed in \cite{Sz04} and
shown in \cite{Ber12}, the exponent $d$ is optimal). Indeed, in the uniformly
elliptic case a sharp bound has been obtained by Berger \cite{Ber12}, only in dimensions
$d\ge 4$ using an approach different from the one presented in this
paper. We have therefore not included
them in this article.

The proof of Theorem \ref{theorem1} requires extending the
methods that have already been developed within the
context of random walks in uniformly elliptic random environments.
Its proof is presented in section \ref{pol-t}.
To do this, we first need to show as in \cite{BDR12},
that the polynomial condition $(P)_M$ for $M\ge 15d+5$,
implies the so called effective criterion, defined
by Sznitman in \cite{Sz02} for random walks
in uniformly elliptic environments, and extended here for
random walks in random environments satisfying condition
$(E)_0$. Two renormalization
methods are employed here, which need to take into account
the fact that the environment is not necessarily uniformly elliptic.
These are developed in subsections \ref{pol-t0} and \ref{t0-aqee}. In subsection \ref{ec-tp}
it is shown, following \cite{Sz02}, that the effective criterion
implies condition $(T')$. The adaptation of the methods
of \cite{BDR12} and \cite{Sz02} from uniformly elliptic
environments to environments satisfying some of the
ellipticity conditions that have been introduced is
far from being straightforward.

The proof of Theorems \ref{theorem2} and \ref{theorem3}, is presented in
sections \ref{section-four} and \ref{tail-regen}. In section \ref{section-four}, an atypical quenched
exit estimate is derived which requires a very careful
choice of the renormalization method, and includes the
definition of an event which we call the {\it constrainment event},
which ensures that the random walk will be able to find a
path to an exit column where it behaves as if  the environment
was uniformly elliptic. In section \ref{tail-regen}, we derive the moments
estimates of the regeneration time of the random walk
using the atypical quenched exit estimate of section \ref{section-four}.
Here,  condition
 $(E')_1$  towards the asymptotic direction is required,
and  appears as the possibility of finding an appropriate
path among $4d-2$ possibilities connecting
two points in the lattice.

\section{Notation and preliminary results}
\label{section2}
\setcounter{equation}{0}
Here we will fix up the notation of the paper and will introduce
the main tools that will be used. In subsection \ref{secprop} we will
prove Proposition \ref{proposition11}. Its proof is straightforward,
but instructive.

\subsection{Setup and background}

Throughout the whole paper we will use letters
without sub-indexes like
$c$, $\rho$ or $\kappa$ to denote any generic constant,
while we will use the notation $c_{3,1},c_{3,2},\ldots,
c_{4,1},c_{4,2},\ldots$ to denote the specific
constants which appear in each section of the paper.
Thus, for example $c_{4,2}$ is the second constant of
section $4$. On the other hand, we will use
$c_1,c_2,c_3,c_4, c_1'$  and $c_2'$ for
specific constants which will appear several
times in several sections.
 Let $c_1\ge 1$ be any constant such that for any pair
of points $x,y\in\mathbb Z^d$, there exists a nearest neighbor path
between $x$ and $y$ with less than

\begin{equation}
\label{c1}
  c_1 |x-y|_2
\end{equation}
sites.
Given $U\subset\mathbb Z^d$, we will
denote its outer boundary by

$$
\partial U:=\{x\notin U: |x-y|_1=1, \ {\rm for}\ {\rm some}\
y\in U\}.
$$
We define $\{\theta_n:n\ge 1\}$ as the canonical time shift
on ${(\mathbb Z^d)}^{\mathbb N}$.
For $l\in\mathbb S^{d-1}$ and $u\ge 0$, we define the times

\begin{equation}
\label{tele}
T^l_u:=\inf\{n\ge 0: X_n\cdot l\ge u\}
\end{equation}
and

$$
\tilde T^l_u:=\inf\{n\ge 0: X_n\cdot l\le u\}.
$$
Throughout, we will denote any nearest neighbor path
with $n$ steps
joining two points $x,y\in\mathbb Z^d$ by $(x_1,x_2,\ldots,x_n)$,
where $x_1=x$ and $x_n=y$. Furthermore,
we will employ the notation

\begin{equation}
\label{delta}
\Delta x_i:=x_{i+1}-x_i,
\end{equation}
for $1\le i\le n-1$,
to denote the directions of the jumps through this path.
Finally, we will call $\{t_x:x\in\mathbb Z^d\}$ the canonical
shift defined on $\Omega$ so that for $\omega=\{\omega(y):y\in
\mathbb Z^d\}$,

\begin{equation}
\label{shift}
t_x(\omega)=\{\omega(x+y):y\in\mathbb Z^d\}.
\end{equation}
Let us now define the concept of {\it regeneration times}
with respect to direction $l$.
Let

\begin{equation}
\label{defa}
a>2\sqrt{d}
\end{equation}
and
$$
D^l:=\min\{n\ge 0: X_n\cdot l<X_0\cdot l\}.
$$
Define $S_0:=0$, $M_0:=X_0\cdot l$,

$$
S_1:=T^l_{M_0+a},\quad R_1:=D^l\circ\theta_{S_1},
$$

$$
M_1:=\sup\{X_n\cdot l:0\le n\le R_1\},
$$
and recursively for $k\ge 1$,

$$
S_{k+1}:=T^l_{M_k+a},\quad R_{k+1}:=D^l\circ\theta_{S_{k+1}}+S_{k+1},
$$

$$
M_{k+1}:=\sup\{X_n\cdot l:0\le n\le R_{k+1}\}.
$$
Define the {\it first regeneration time} as

$$
\tau_1:=\min\{k\ge 1: S_k<\infty, R_k=\infty\}.
$$
The condition (\ref{defa}) on $a$ will be eventually useful
to prove the non-degeneracy of the covariance matrix
of part $(a)$
of Theorem \ref{theorem3}.
Now define recursively in $n$ the $(n + 1)$-st regeneration time $\tau_{n+1}$
 as  $\tau_1(X_\cdot)+\tau_n(X_{\tau_1+\cdot}-X_{\tau_1})$.
Throughout the sequel, we will occasionally write
$\tau_1^l,\tau_2^l,\ldots$ to emphasize the dependence
of the regeneration times with respect to the chosen direction.
It is a standard fact 
to show that the sequence
$( (\tau_1, X_{(\tau_1+\cdot)\land\tau_2}-X_{\tau_1}),
 (\tau_2-\tau_1, X_{(\tau_2+\cdot)\land\tau_3}-X_{\tau_2}),
\ldots )$ is independent and except for its first term also
i.i.d. with the same law as that of $\tau_1$ with respect
to the conditional probability measure $P_0(\cdot|D^l=\infty)$
(see for example Sznitman and Zerner \cite{SZ99}, whose
proof of this fact is also valid without the uniform ellipticity assumption).
This implies the following theorem (see Zerner \cite{Z02} and Sznitman
and Zerner \cite{SZ99} and Sznitman \cite{Sz00}).

\medskip

\begin{theorem} {\bf (Sznitman and Zerner \cite{SZ99}, Zerner
\cite{Z02}, Sznitman \cite{Sz00})}
\label{lwlln}
Consider a RWRE in an elliptic i.i.d. environment.
Let $l\in \mathbb S^{d-1}$ and assume that there is
a neighborhood $V$ of $l$ such that for every $l'\in V$ the
random walk is transient in the direction $l'$. Then there
is a deterministic $v$ such that $P_0$-a.s. one has that

$$
\lim_{n\to\infty}\frac{X_n}{n}=v.
$$
Furthermore, the following are satisfied.

\begin{itemize}

\item[a)] If $E_0[\tau_1]<\infty$, the walk is ballistic and
$v\ne 0$.

\item[b)] If $E_0[\tau_1^2]<\infty$ we have that

$$
\epsilon^{1/2}\left(X_{[\epsilon^{-1}n]}-[\epsilon^{-1}n]v\right)
$$
converges in law under $P_0$ to a Brownian motion with non-degenerate covariance matrix.

\end{itemize}
\end{theorem}

\medskip
In 2009, both Rassoul-Agha and Sepp\"{a}l\"{a}inen in \cite{RAS09} and
Berger and Zeitouni in \cite{BZ08} were able to prove a quenched
central limit theorem under good enough moment conditions
on the regeneration times. The result of Rassoul-Agha and
Sepp\"{a}l\"{a}inen which does not require a uniform ellipticity assumption
is the following one.

\medskip
\begin{theorem}
{\bf (Rassoul-Agha and Sepp\"{a}l\"{a}inen \cite{RAS09})}
\label{lwlln-quenched}
Consider a RWRE in an elliptic i.i.d. environment.
Let $l\in \mathbb S^{d-1}$ and let $\tau_1$ be
the corresponding regeneration time.
Assume that

$$
E_0[\tau_1^p]<\infty,
$$
for some $p>176d$. Then $\P$-a.s. we have that

$$
\epsilon^{1/2}\left(X_{[\epsilon^{-1}n]}-[\epsilon^{-1}n]v\right)
$$
converges in law under $P_{0,\omega}$ to a Brownian motion with non-degenerate covariance matrix.
\end{theorem}
\medskip

\noindent We now define the {\it $n$-th regeneration radius} as

$$
X^{*(n)}:=\max_{\tau_{n-1}\le k\le \tau_n}|X_k-X_{\tau_{n-1}}|_1.
$$

The following theorem was stated and proved
without using uniform ellipticity
by Sznitman as Theorem A.2 of  \cite{Sz02},
and provides a control on the lateral displacement of the
random walk with respect to the asymptotic direction.
We need to define for $z\in\mathbb R^d$

$$
\pi(z):=z-(z\cdot\hat v)\hat v.
$$
\begin{theorem} {\bf (Sznitman \cite{Sz02})}
\label{sznitman01}
Consider a RWRE in an elliptic i.i.d. environment satisfying condition $(T)_\gamma|l$.
Let $l\in\mathbb S^{d-1}$ and $\gamma\in (0,1)$.
Then, for any $c>0$
and $\rho\in(0.5,1)$,

$$
\limsup_{u\to\infty}u^{-(2\rho-1)\land (\gamma \rho)}\log P_0
\left(\sup_{0\le n\le T_u^l}|\pi(X_n)|\ge cu^\rho\right)<0,
$$
where $T_u^l$ is defined in (\ref{tele}).
\end{theorem}
Define the function $\gamma_L: [3,\infty)\to \mathbb R $
as

\begin{equation}
\label{choice-gamma}
\gamma_L:=\frac{\log 2}{\log\log L}.
\end{equation}
Given $l\in\mathbb S^{d-1}$, we
say that condition $(T)_{\gamma_L}$  in direction $l$ (also written
as $(T)_{\gamma_L}|l$) is satisfied if there exists a neighborhood
$V\subset\mathbb S^{d-1}$ of $l$ such that for all $l'\in V$

$$
\limsup_{L\to\infty}\frac{1}{L^{\gamma_L}}\log
P_0(X_{T_{U_{l',L}}}\cdot l'<0)<0,
$$
where the slabs $U_{l',L}$ are defined in (\ref{slab}).

Throughout this paper we will also need
the following generalization of an equivalence proved by Sznitman \cite{Sz02},
for the case $\gamma\in (0,1)$
and which does not require uniform ellipticity.
 It
is easy to extend Sznitman's proof to include the case $\gamma =\gamma_L$.

\medskip

\begin{theorem}
\label{sznitman-equivalence} {\bf (Sznitman \cite{Sz02})}
Consider a RWRE in an elliptic i.i.d. environment. Let $\gamma\in [0, 1)$
 and $l \in \mathbb S^{d-1}$. Then the following are equivalent.

\begin{itemize}
\item[(i)] Condition $(T)_\gamma|l$ is satisfied.
\item[(ii)] $P_0(A_l) = 1$ and if $\gamma>0$ we have that $E_0[ \exp\{c(X^{∗(1)})^\gamma\}] < \infty$
for some $c > 0$, while if $\gamma =\gamma_L$ we have that
$E_0\left[ \exp\left\{c(X^{*(1)})^{\frac{\log 2}{\log\log \left(3\lor X^{*(1)}\right)}}\right\}\right] < \infty$ for some $c>0$.
\item[(iii)] There is an asymptotic direction $\hat v$ such
that $l\cdot\hat v>0$ and
for every $l'$ such that $l'\cdot \hat v>0$ one has that
$(T)_\gamma|l'$ is satisfied.
\end{itemize}
\end{theorem}
\medskip

\noindent The following corollary of Theorem \ref{sznitman-equivalence}
will be important.

\begin{corollary} {\bf (Sznitman \cite{Sz02})}
\label{regeneration time lemma}
Consider a RWRE in an elliptic i.i.d. environment.
Let $\gamma \in (0,1)$  and $l\in\mathbb S^{d-1}$.  Assume that
$(T)_{\gamma}|_l$ holds. Then there exists a constant $c$  such that
for every $L$  and $n\ge 1$ one has that
\begin{equation}
\label{regeneration time estimate}
P_0 \left(X^{*(n)}>L \right) \leq \frac{1}{c} e^{-cL^{\gamma}}.
\end{equation}
\end{corollary}

\subsection{Comments and proof of Proposition \ref{proposition11}}
\label{secprop}

Here we will show that   $(E')_1$ implies $(ES)$.
We will do this passing through another ellipticity condition.
We say that condition $(ES')$ is satisfied if
there
exist non-negative real numbers
 $\alpha_1,\ldots,\alpha_d$ and $\alpha'_1,\ldots,\alpha'_d$
such that

$$
\min_{1\le i\le d}(\alpha_i+\alpha'_i)>1
$$
and

\begin{equation}
\label{2nd condition}
\max_{1\le i\le d}
\E\left[\left(\frac{1}{1-\omega(0,e_i)}\right)^{\alpha_i}\right]<\infty\quad
{\rm and}\quad
\max_{1\le i\le d}
\E\left[\left(\frac{1}{1-\omega(0,e_{i+d})}\right)^{\alpha'_i}\right]<\infty.
\end{equation}
We have the following lemma.

\medskip

\begin{lemma} Consider a random walk in an i.i.d. random
environment. Then  condition $(E')_{1}$
 implies $(ES')$
which in turn implies $(ES)$.  Furthermore, for a random
walk in a random Dirichlet environment, $(E')_1$, $(ES)$ and $(ES')$
are equivalent to $\lambda>1$ (cf. (\ref{lambda-st})).
\end{lemma}

\medskip

\begin{proof}
We  first prove  that $(ES')$ implies $(ES)$. Note first that by
the independence between $\omega(0,e_i)$ and $\omega(e_i,-e_i)$, (\ref{2nd condition}) is equivalent to

 $$
\max_{1\le i\le d} \E\left[\left(\frac{1}{1-\omega(0,e_i)}\right)^{\alpha_i}\left(\frac{1}{1-\omega(e_i,-e_{i})}\right)^{\alpha'_i}\right]<\infty.
 $$
 Then
 it is enough to prove that for each pair of real numbers $u_1, u_2$ in $(0,1)$ one has that

\begin{equation}
\label{suffcond}
\frac{1}{1-u_1u_2} \leq \frac{1}{(1-u_1)^{\alpha}(1-u_2)^{\alpha'}}.
\end{equation}
for any $\alpha$, $\alpha' \geq 0$ such that $\alpha+\alpha'>1$.  Now if we denote by $v_1=1-u_1$ and $v_2=1-u_2$ then (\ref{suffcond}) is equivalent to
\begin{equation}
\label{suffcond1}
v_1^{\alpha}v_2^{\alpha'} \leq v_1+v_2-v_1v_2
\end{equation}
But (\ref{suffcond1}) follows easily by our conditions on $v_1$, $v_2$,
$\alpha$ and $\alpha'$.
To prove that
$(E')_1$ implies $(ES')$, we choose for each $1 \leq i \leq d$
$$
\alpha_i:=\sum_{e\neq e_i}\alpha(e), \quad
\alpha_i':=\sum_{e\neq e_{i+d}}\alpha(e).
$$
Note in particular that

\begin{equation}
\label{1cond}
\alpha_i+\alpha_i'>1, \quad \forall i \in \{1, \ldots, d\}.
\end{equation}

\noindent
Now, it is easy to check that $(E')_1$ implies that

$$
\E\left[e^{\sum_{e\ne e_i}\alpha(e)\log\frac{1}{\omega(0,e)}}\right]<\infty.
$$
Therefore, by the monotonicity of the function $\log x$ we have that
$(E')_1$ implies that

\begin{equation}
\label{2cond}
\E \left(e^{-\sum_{e \neq e_i}\alpha(e)\log \omega(0,e)}\right) <\infty, \quad  \alpha_i \log \sum_{e \neq e_i} \omega(0, e) \geq \sum_{e \neq e_i}\alpha(e)\log \omega(0,e),
\end{equation} 
\noindent and the corresponding inequalities with $e_i$ replaced by $-e_i$, for each $1 \leq i \leq d$. Then $(ES')$ follows by (\ref{1cond})
and (\ref{2cond}). Let us now consider a Dirichlet random environment
of parameters $\{\beta_e:e\in U\}$ and assume that $(ES)$ is satisfied.
It is well known that for each $e\in U$, the random variable
$\omega(0,e)$ is distributed according to a Beta distribution
of parameters $(\beta_e ,\sum_{e'\in U} \beta_{e'}-\beta_e)$: so it
has a density 

\begin{equation}
\label{betalaw}
cx^{\beta_e}(1-x)^{\sum_{e'}\beta_{e'}-\beta_e}
\end{equation}
 with
respect to the Lebesgue measure in $[0,1]$. On the other hand, if we
define $v_1:=1-\omega(0,e)$ and $v_2:=1-\omega(e,-e)$, we see
that 

$$
\E\left[\frac{1}{1-\omega(0,e)\omega(e,-e)}\right]<\infty
$$
 is satisfied if and only if

$$
\frac{1}{v_1+v_2}
$$
is integrable. But by (\ref{betalaw}) this happens whenever

$$
2\sum_{e'}\beta_{e'}-\beta_e-\beta_{-e}>1.
$$
This proves that for a Dirichlet random environment $(ES)$ implies 
that $\lambda>1$. It is obvious that for a Dirichlet random environment
$\lambda>1$ implies $(E')_1$.

\end{proof}

\medskip

Let us now prove Proposition \ref{proposition11}.
If
 the random walk is not transient in any direction,
there is nothing to prove. So assume that
the random walk is transient in a direction $l$
and hence the corresponding regeneration times
are well defined.
Essentially, we will exhibit a trap as the one
depicted in Figure \ref{Fribergh}, in the edge
$\{0,e_i\}$.
 Define the first exit time
of the random walk from the edge $\{0,e_i\}$, so that

$$
F:=\min\left\{n\ge 0: X_n\notin \{0,e_i\}\right\}.
$$
 We then have
for every $k\ge 0$ that

%$$
%P_{0,\omega}(F=2k+2)=\omega_1^{k+1}\omega_2^{k}(1-\omega_2),
%$$
%and

%$$
%P_{0,\omega}(F=2k+1)=\omega_1^k\omega_2^k(1-\omega_1).
%$$
%Hence,

\begin{equation}
\nonumber
P_{0,\omega}(F> 2k)=(\omega_1\omega_2)^k
\end{equation}
and

\begin{equation}
\label{tau1m}
\sum_{k=0}^\infty P_{0,\omega}(F>2k)=\frac{1}{1-\omega_1\omega_2}.
\end{equation}
This proves that under the annealed law,

$$
E_0(F)=\infty.
$$
We can now show using the strong Markov property
under the quenched measure and the i.i.d. nature of
the environment, that for each natural $m>0$, the time
$T_{m}:=\min\{n\ge 0:X_n\cdot l>m\}$ can be bounded from
below by a sequence $F_1,\ldots,F_m$ of random variables which
under the annealed measure are i.i.d. and distributed as $F$.
This proves that $P_0$-a.s. $T_{m}/m\to\infty$ which
implies that the random walk is not ballistic in direction $l$.

\section{Equivalence between the polynomial ballisticity condition
and $(T')$}
\label{pol-t}
\setcounter{equation}{0}
Here we will prove Theorem \ref{theorem1}, establishing the
equivalence between the polynomial condition $(P)_M$ and
condition $(T')$. To do this, we will pass through both
the effective criterion and an version of condition $(T)_\gamma$
which corresponds to the choice of $\gamma=\gamma_L$
according to (\ref{choice-gamma}) (see \cite{BDR12}).
Now, to prove Theorem \ref{theorem1}, we will first
show in subsection \ref{pol-t0} that $(P)_M$ implies $(T)_{\gamma_L}$ for
$M\ge 15d +5$.
 In subsection \ref{t0-aqee}, we will prove that
$(T)_{\gamma_L}$ implies a weak kind of an atypical quenched exit estimate.
In these first two steps, we will generalize
the methods presented in \cite{BDR12} for random walks
satisfying condition $(E)_0$.
In subsection \ref{aqee-ec}, we will see that this estimate implies the
effective criterion. Finally, in subsection \ref{ec-tp}, we will
show that the effective criterion implies $(T')$, generalizing
the method presented by Sznitman \cite{Sz02}, to
random walks satisfying $(E)_0$.

Before we continue, we will need some
additional notation. Let $l\in\mathbb S^{d-1}$.
Let
$L,L'>0$, $\tilde L>0$,

\begin{equation}
\label{box1}
B(R,L,L',\tilde L):=R\left((-L,L')\times (-\tilde L,\tilde L)^{d-1}\right)\cap
\mathbb Z^d
\end{equation}
and

\begin{equation}
\nonumber
\partial_+B(R,L,L',\tilde L):=\partial B\cap\left\{x\in\mathbb Z^d:
x\cdot l\ge L', |R(e_j)\cdot x|<\tilde L,\ {\rm for}\ {\rm each}\
2\le j\le d\right\}.
\end{equation}
Here $R$ is the rotation defined by (\ref{defr}). When there is no risk of confusion, we will
drop the dependence of
$B(R,L,L',\tilde L)$ and
$\partial_+B(R,L,L',\tilde L)$
with respect to $R$, $L$, $L'$ and $\tilde L$ and
write $B$ and $\partial_+B$ respectively.
Let also,

$$
\rho_B:=\frac{P_{0,\omega}(X_{T_B}\notin\partial_+B)}
{P_{0,\omega}(X_{T_B}\in\partial_+B)}=\frac{q_B}{p_B},
$$
where $q_B:=P_{0,\omega}(X_{T_B}\notin\partial_+B)$ and
$p_B:=P_{0,\omega}(X_{T_B}\in\partial_+B)$.

\subsection{Polynomial ballisticity implies $(T)_{\gamma_L}$}
\label{pol-t0}
Here we will prove that the Polynomial ballisticity condition implies
$(T)_{\gamma_L}$. To do this, we will use a multi-scale renormalization scheme
as presented in Section 3 of \cite{BDR12}. Let us note that \cite{BDR12}
assumes that the walk is uniformly elliptic.

\begin{proposition}
\label{PM-T0}
Let $M >15d+5$ and $l \in \mathbb S^{d-1}$. Assume that conditions
$(P)_M|l$ and $(E)_0$ are satisfied. Then $(T)_{\gamma_L}|l$ holds.
\end{proposition}

\noindent Let us now  prove Proposition \ref{PM-T0}.
%\begin{equation}
%\label{c0}
%c_0:=2^{3(d-1)} \vee e^{2 \log 90+\sum_{j=1}^{\infty}\frac{\log j}{2^j}},
%\end{equation}
Let $N_0\ge \frac{3}{2}c_0$, where $c_0$ is defined in (\ref{deflo}). For $k \ge 0$, define recursively the scales

\begin{equation}
\nonumber
N_{k+1}:=3(N_0+k)^2N_k.
\end{equation}
Define also for $k \ge 0$ and $x \in \R^d$ the boxes

\begin{equation}
\nonumber
B(x,k):=\left\{y \in \Z^d: -\frac{N_k}{2} < (y-x)\cdot l <N_k, \,\,
  |(y-x)\cdot R(e_i)|<25 N_k^3\ {\rm{for}}\ 2\le i\le d \right\}
\end{equation}
and their {\it middle frontal part}

\begin{equation}
\nonumber
\tilde{B}(x,k):=\left\{y \in \Z^d: N_k-N_{k-1} \leq (y-x)\cdot l <N_k, \,\,
|(y-x)\cdot R(e_i)|< N_k^3\ {\rm{for}}\ 2\le i\le d \right\}
\end{equation}
with the convention that $N_{-1}:=2N_0/3$. We also define the {\it the front side}

\begin{equation}
\nonumber
\partial_{+}B(x,k):=\{y \in \partial B(x,k): (y-x) \cdot l \geq N_k\},
\end{equation}
the {\it back side}

\begin{equation}
\nonumber
\partial_{-}B(x,k):=\{y \in \partial B(x,k): (y-x) \cdot l \leq -\frac{N_k}{2}\},
\end{equation}
and the {\it lateral sides}

\begin{equation}
\nonumber
\partial_{l}B(x,k):=\{y \in \partial B(x,k): |(y-x)\cdot R(e_i)|\geq 25N_k^3
\ {\rm for}\ 2\le i\le d\}.
\end{equation}

\noindent We need to define for each $n, m \in \N$ the sub-lattices

$$
\mathcal L_{n, m}:=\{x \in \Z^d: [x \cdot l] \in n \Z, \,\, [x \cdot R(e_j)] \in m\Z, \, \, {\rm for} \,\, 2 \leq j \leq d \}
$$
and refer to the elements of

$$
\mathcal B_k:=\left\{B(x,k): x \in \mathcal L_{N_{k-1}-1,N_k^3-1}\right\}
$$
as {\it boxes of scale} $k$. When there is no risk of confusion, we will denote a typical element of this set by $B_k$ or simply $B$ and its middle part as $\tilde B_k$ or $\tilde B$. Furthermore, we have

$$
\cup_{B \in \mathcal B_k}\tilde B=\Z^d,
$$
which will be an important property that will be useful.
In this subsection, it is enough to assume a weaker condition than $(P)_M|l$.
The following lemma gives a practical
version of the polynomial condition which will be used for
the proof of Proposition \ref{PM-T0}.

\medskip
\begin{lemma}
\label{lemma-equiv} Let $M>0$ and $l\in\mathbb S^{d-1}$.
Assume that condition $(P)_M|l$ is satisfied.
 Then, for $N_0\ge \frac{3}{2}c_0$ one has that

\begin{equation}
\nonumber
\sup_{x \in \tilde B_0}P_x \left(X_{T_{B_0}} \not \in \partial_{+}B_0\right) < N_0^{-M}.
\end{equation}
\end{lemma}
\begin{proof} For each $x\in\tilde B_0$, consider the
box 

$$
A_0(x):=\left\{y\in\mathbb Z^d: -\frac{2}{3}N_0\le (y-x)\cdot l\le\frac{2}{3}N_0,
|(y-x)\cdot R(e_i)|\le 24 N_0^3\ {\rm for}\ 2\le i\le d\right\}.
$$
Note that

\begin{equation}
\label{lastlast}
P_0(X_{T_{B_0}}\notin\partial B_0)\le P_0(X_{T_{A_0(x)}}\notin\partial A_0(x)).
\end{equation}
Let $L:=\frac{2}{3}N_0$. Note now that for all $x$,  $A_0(x)$
is a translation of the box $B_{l,L,81L^3}$ (c.f. (\ref{boxbox})) so that the probability
of the right-hand side of (\ref{lastlast}) is bounded from above
by $L^{-M}\le N_0^{-M}$.

\end{proof}
\medskip

\noindent
%Note that similarly to $(P)_M|l$, (\ref{wpolyc}) holds for all $l'$ in a neighbourhood of $l$.

\noindent We now say that box $B \in \mathcal B_0$ is {\it good} if
\begin{equation}
\nonumber
\sup_{x \in \tilde B_0}P_{x,\omega} \left(X_{T_{B_0}} \not \in \partial_{+}B_0\right) < N_0^{-5}.
\end{equation}
Otherwise, we say that the box $B\in\mathcal B_0$ is {\it bad}.
The following lemma appears in \cite{BDR12} as  Lemma 3.7, so
its proof will also be omitted.

\medskip
\begin{lemma}
\label{k=0}
Let $M>0$ and $l \in \mathbb S^{d-1}$.
 Assume that $(P)_M|l$ holds. Let $N_0\ge \frac{3}{2}c_0$.
 Then for
all $B_0 \in \mathcal B_0$ we have that

$$
\P(B_0 \,\, {\rm is} \,\, {\rm good}) \ge 1-2^{d-1}N_0^{3d+3-M}.
$$
\end{lemma}

%\begin{proof}
%Note that

%\begin{equation}
%\label{B0 bad}
%\P(B_0 \,\, {\rm is} \,\, {\rm bad}) \leq \sum_{x \in \tilde B_0}\P\left(P_{x,\omega} \left(X_{T_{B_0}} \not \in \partial_{+}B_0\right)\ge N_0^{-5}\right).
%\end{equation}
%Now by Markov's inequality we have for $x \in \tilde B_0$ that

%\begin{equation}
%\label{B0 bad1}
%\P\left(P_{x,\omega} \left(X_{T_{B_0}} \not \in \partial_{+}B_0\right)\ge N_0^{-5}\right) \leq N_0^5 \sup_{x \in \tilde B_0}P_x \left(X_{T_{B_0}} \not \in \partial_{+}B_0\right).
%\end{equation}
%Now, with the help of Lemma \ref{lemma-equiv}, (\ref{B0 bad}), (\ref{B0 bad1}) and from a routine counting argument we obtain

%$$
%\P(B_0 \,\, {\rm is} \,\, {\rm bad}) \leq 2^{d-1}N_0^{3d+3-M}.
%$$

%\end{proof}

\medskip

\noindent Now, we want to extend the concept of good and bad boxes of scale
$0$
to boxes of any scale $k\ge 1$.
To do this, due to the lack of uniform ellipticity, we need
to modify the notion of good and bad boxes for scales $k\ge 1$
presented in Berger, Drewitz and Ram\'\i rez \cite{BDR12}.
Consider a box $Q_{k-1}$ of scale $k-1\ge 1$.
For each $x \in \tilde Q_{k-1}$ we associate a
natural number $n_x$ and a
self-avoiding path
$\pi^{(x)}:=(\pi^{(x)}_1,\ldots, \pi^{(x)}_{n_x})$
starting from $x$ so that $\pi^{(x)}_1=x$,
such that $(\pi^{(x)}_{n_x}-x)\cdot l\ge N_{k-2}$
and so that

$$
c_{3,1} N_{k-2}\le n_x\le c_{3,2} N_{k-2},
$$
for some 

$$c_{3,1}:=\frac{1}{2\sqrt{d}},$$
 and 

$$c_{3,2}:=\frac{2}{\sqrt{d}}.$$
%We will also choose the path $\pi^{(x)}$ so that
%it either is entirely within $\tilde Q_{k-1}$ or else if it exits
%the box $Q_{k-1}$ it does it through $\partial_+ Q_{k-1}$.
Now, let

\begin{equation}
\label{defxi}
\Xi:=\frac{1}{2}e^{-\frac{c_{3,2}\log \eta_{\alpha}+9d}{c_{3,1}}}.
\end{equation}
We say that the box $Q_{k-1} \in\mathcal{B}_{k-1}$ is {\it elliptically good} if for each $x \in \tilde Q_{k-1}$ one has that

\begin{equation}
\nonumber
 \sum_{i=1}^{n_x}
\log
 \frac{1}{\omega(\pi^{(x)}_i,
 \Delta \pi^{(x)}_i)}\leq n_x \log \left(\frac{1}{\Xi}\right).
 \end{equation}
 Otherwise the box is called {\it elliptically bad}.
We can now recursively define the concept of good and
bad boxes. For $k \geq 1$ we say that a box $B_k
\in \mathcal B_k$ is {\it good}, if the
following are satisfied:

\begin{enumerate}

\item There is a box $Q_{k-1} \in \mathcal
B_{k-1}$
which is elliptically good.

\item Each box $C_{k-1} \in \mathcal B_{k-1}$ of scale $k-1$ satisfying
  $C_{k-1} \cap Q_{k-1}  \neq \emptyset$ and $C_{k-1} \cap B_k \neq \emptyset$
  is  elliptically good.

\item Each box $B_{k-1} \in \mathcal B_{k-1}$ of scale $k-1$ satisfying
$B_{k-1} \cap Q_{k-1} = \emptyset$ and $B_{k-1} \cap B_k \neq \emptyset$, is
good.

\end{enumerate}

\noindent Otherwise, we say that the box $B_k$ is {\it bad}.
 Now we will obtain an important estimate on the probability that a box of
 scale $k\ge 1$ is good, corresponding to Proposition 3.8 of
 \cite{BDR12}. Nevertheless, note that
here we have to deal with our different definition of good and bad boxes
due to the lack of uniform ellipticity. Let

$$c_{3,3}:=c_{3,1}\log \frac{1}{\Xi}-c_{3,2} \log \eta_{\alpha}-9d=c_{3,1} \log 2 >0.$$
We first need the following estimate.

\medskip

\begin{lemma}
\label{nueb3.1}
For each $k \geq 1$   we have that
\begin{equation}
\label{nuebe.3.1}
\P(B_k \,\, {\rm is} \,\, {\rm not}  \,\, {\rm
elliptically}  \,\, {\rm good}) \leq e^{-c_{3,3}N_{k-1}}.
\end{equation}
\end{lemma}

\begin{proof}
By translation invariance and using Chebyshev's inequality as well as independence,
we have that for any $\alpha>0$
\begin{eqnarray*}
 & \P(B_k \,\, {\rm is} \,\, {\rm not}  \,\, {\rm elliptically }  \,\, {\rm good}) \le \sum_{x\in \tilde B_k} \P \left(   \sum_{i=1}^{n_x}
\log
 \frac{1}{\omega(\pi^{(x)}_i,
 \Delta \pi^{(x)}_i)}> n_x \log \left(\frac{1}{\Xi}\right)  \right) \\
  & \leq N_{k-1}N_k^{3(d-1)}
e^{-N_{k-1}\left(c_{3,1}\alpha \log\left(\frac{1}{\Xi}\right)
-c_{3,2}\log\eta_\alpha\right)}\\
& \leq e^{-N_{k-1}\left(c_{3,1}\log\left(\frac{1}{\Xi}\right)
-c_{3,2}\log\eta_\alpha-9d\right)}
\end{eqnarray*}
where $ N_{k-1}N_k^{3(d-1)}$ is an upper bound for $|\tilde B_k|$
and we have used the inequality $N_k\le 12 N_{k-1}^3$
and we have without loss of generality assumed that $\alpha<1$.
 But this expression can
be bounded by $e^{9dN_{k-1}}$ due to our choice of $N_0$. Then, using the definition of $\Xi$ in  (\ref{defxi}), we have that

$$
\P(B_k \,\, {\rm is} \,\, {\rm not}  \,\, {\rm elliptically }  \,\, {\rm good}) \leq e^{-c_{3,3}N_{k-1}}.
$$
\end{proof}

\medskip
\noindent We can now state the following lemma giving an estimate for the probability
that a box of scale $k\ge 0$ is bad.
We will use Lemma \ref{lemma-equiv}.

\medskip

\begin{lemma}
\label{k>0}
Let $l \in \mathbb S^{d-1}$ and $M\ge 15d+5$. Assume that $(P)_M|l$ is satisfied. Then for $N_0\ge \frac{3}{2}c_0$
one has that for all $k \geq 0$ and all $B_k \in \mathcal B_k$,

\begin{equation}
\nonumber
\P(B_k \,\, {\rm is} \,\, {\rm good}) \ge 1-e^{-2^k}.
\end{equation}

\end{lemma}
\begin{proof} By Lemma \ref{k=0} we see that

$$
\P(B_0\ {\rm is}\ {\rm bad})\le
e^{-c'_{3,0}},
$$
where

$$
c'_{3,0}:=\log\frac{N_0^{M-3d-3}}{2^{d-1}}.
$$
We will show that this implies for all $k\ge 1$ that

\begin{equation}
\label{beka}
\P(B_k\ {\rm is}\ {\rm bad})\le e^{-c'_{3,k}2^k},
\end{equation}
for a sequence of constants $\{c'_{3,k}: k\ge 0\}$
defined recursively by

\begin{equation}
\label{ceka}
c'_{3,k+1}:=c'_{3,k}-\frac{\log\left(3^{16d}(N_0+k)^{12d}\right)}{2^{k+1}}.
\end{equation}
 We will now prove (\ref{beka})
using induction on $k$. To simplify notation, we will
denote by $q_k$ for $k\ge 0$, the probability that
the box $B_k$ is bad. Assume that (\ref{beka})  is true for some 
$k\ge 0$.
Let $A$ be the event that all boxes of scale $k$ that intersect $B_{k+1}$
are elliptically good, and $B$ the event that each pair of bad boxes of scale
$k$ that intersect $B_{k+1}$, have a non-empty intersection. Note that the event
$A\cap B$ implies that the box $B_{k+1}$ is good. Therefore,
the probability $q_{k+1}$ that the box $B_{k+1}$ is bad
is bounded by the probability that there are at least two bad
boxes $B_{k}$ which intersect $B_{k+1}$ plus the probability that
there is at least one elliptically bad box of scale $k$, so that by
Lemma \ref{nueb3.1}, for
each
$k \geq 0$ one has that

\begin{equation}
\label{quk}
q_{k+1}\le m_{k}^2
q_k^2+m_k
e^{-c_{3,3}N_k},
\end{equation}
where $m_k$ is the total number of bad boxes of scale $k$ that intersect
$B_{k+1}$. Now note that

\begin{equation}
\label{memk}
\sqrt{2}m_k\le 3^{8d}(N_0+k)^{6d}.
\end{equation}
But by the the fact that $c_{3,3} N_k \ge c'_{3,k}2^{k+1}$ for $k\ge
0$
we have that

$$
e^{-c_{3,3}N_k}\le e^{-c'_{3,k}2^{k+1}}.
$$
Hence, substituting this estimate and estimate (\ref{memk}) back into (\ref{quk}) and using the
induction hypothesis, we conclude that

$$
q_{k+1}\le  3^{16d}(N_0+k)^{12d}e^{-c'_{3,k}2^{k+1}}=e^{-c'_{3,k+1}2^{k+1}}.
$$
 Now note that
the recursive definition (\ref{ceka}) implies that

$$
c'_{3,k}\ge  \log\frac{N_0^{M-3d-3}}{2^{d-1}}
-\sum_{k=0}^\infty \frac{\log\left(3^{16d}(N_0+k)^{12d}\right)}{2^{k+1}}.
$$
Using the inequality $\log(a+b)\le \log a+\log b$ valid for $a,b\ge 2$,
we see that

$$
\sum_{k=0}^\infty \frac{\log\left(3^{16d}(N_0+k)^{12d}\right)}{2^{k+1}}
\le 16d\log 3+12d\log N_0+12d.
$$
From these estimates we see that whenever $M\ge 15d+5$ and

\begin{equation}
\label{ene0}
\log N_0 -\log 2^{d-1}3^{16d}e^{12d+1}\ge 0,
\end{equation}
then for every $k\ge 0$, one has that $c'_{3,k}\ge 1$.
But (\ref{ene0}) is clearly satisfied for $N_0\ge 3^{29d}$.

\end{proof}

\medskip

\noindent The next lemma  establishes that the probability that a random
walk exits a box $B_k$ through its lateral or back side is small if this box is good.

\begin{lemma}
\label{pgb3.1} Assume that $N_0\ge \frac{3}{2}c_0$. Then, there is a constant $c_{3,4}>0$ such that for each $k\ge 0$ and $B_k \in \mathcal B_k$ which is good  one has

\begin{equation}
\nonumber
\sup_{x \in \tilde B_k}P_{x,\omega} \left(X_{T_{B_k}} \not \in \partial_{+}B_k\right)\leq e^{-c_{3,4}N_{k}}.
\end{equation}

\end{lemma}

\begin{proof} Since the proof follows closely that of
Proposition 3.9 of \cite{BDR12}, we will only give a sketch 
indicating the steps where the lack of uniform ellipticity requires
a modification. We first note that
 that for each $k\ge 0$,

$$
P_{x,\omega} \left(X_{T_{B_k}} \not \in \partial_{+}B_k\right)
\le
P_{x,\omega} \left(X_{T_{B_k}} \in \partial_{-}B_k\right)
+P_{x,\omega} \left(X_{T_{B_k}} \in \partial_{l}B_k\right).
$$
\noindent We denote by $p_k:=\sup_{x \in \tilde B_k}P_{x,\omega}
\left(X_{T_{B_k}}\in \partial_{l} B_k \right)$ and
$r_k:=\sup_{x \in \tilde B_k}P_{x,\omega}
\left(X_{T_{B_k}}\in \partial_{-} B_k \right)$.
We   show by induction on $k$ that

\begin{eqnarray}
\label{it3.1l}
&p_k \leq e^{-c''_{3,k}N_k}\quad{\rm and}\\
\label{it3.12}
&r_k \leq e^{-c''_{3,k}N_k},
\end{eqnarray}
where

%\begin{eqnarray*}
%&c''_{3,k}:= \frac{5 \log N_0}{N_0}- \sum_{j=1}^k \frac{\log
 %   27(N_0+j)^4}{N_{j-1}}\\
%&
\begin{equation}
\nonumber
c''_{3,k}:= \frac{5 \log N_0}{N_0}- \sum_{j=1}^k \frac{\log 27(N_0+j)^4}{N_{j-1}}-\sum_{j=1}^k \frac{5N_{j-1}+\log 24+6d(\log \Xi)^2N_{j-1}}{N_j},
\end{equation}
and
$\Xi$ is defined  in (\ref{defxi}).
The case $k=0$ follows easily by the definition of good box at scale $0$ with

$$
c''_{3,0}:= \frac{5 \log N_0}{N_0}.
$$

\noindent We then assume that (\ref{it3.1l})
and (\ref{it3.12}) hold for some $k \geq 0$ and  show that this
 implies that (\ref{it3.1l}) is satisfied for $k+1$, following
the same argument presented in the proof of Proposition 3.9
of \cite{BDR12}.
Next, assuming (\ref{it3.1l})
and (\ref{it3.12})  for some $k \geq 0$, we show that
 (\ref{it3.12}) is satisfied for $k+1$. This is
done essentially as in Proposition 3.9 of \cite{BDR12}, coupling
the random walk to a one-dimensional random walk. To
perform this coupling, we cannot use uniform ellipticity
as in \cite{BDR12}, and should instead use 
the property of elliptical goodness satisfied by the corresponding
suboxes. Note that we only need to use the elliptical goodness of
the bad subox and of those which intersect it. We finally choose

$$
c_{3,4}:=\inf_{k\ge 0}c_{3,k}''.
$$
Note that for $N_0\ge\frac{3}{2}c_0$ we have that $\log 27(N_0+j)^4\le 3(N_0+j-1)^2$
for $j\ge 1$ and hence the second term in the definition of $c_{3,k}''$
is bounded from above by  $4/N_0$. Also, the last term
is bounded from above by $(10+6d(\log\Xi )^2)/N_0$. Hence,
we conclude that $c_{3,4}>0$.
\end{proof}

\medskip

\noindent We can now repeat the last argument of Proposition 2.1  of \cite{BDR12}, which does
not require uniform ellipticity, to finish the proof of Proposition \ref{PM-T0}.

\subsection{Condition $(T)_{\gamma_L}$ implies a weak
atypical quenched exit estimate}
\label{t0-aqee}
In this subsection we will prove that the condition $(T)_{\gamma_L}$
implies a weak atypical quenched exit estimate.
Throughout, we will denote by $B$ the box

\begin{equation}
\label{lb}
B:=B(R,L,L,L),
\end{equation}
as defined in (\ref{box1}), with $R$ the rotation which
maps $e_1$ to $l$.
Let

$$
\epsilon_L:=\frac{1}{(\log\log L)^2}.
$$

\medskip

\begin{proposition}
\label{pqexe}
Let  $l\in\mathbb S^{d-1}$. Assume that the
 ellipticity condition
$(E)_0$ and that $(T)_{\gamma_L}|l$ are fulfilled. Then,  for each function $\beta_L:(0, \infty) \to (0, \infty)$ and each $c>0$ there exists $c_{3,11}>0$ such that

\begin{equation}
\label{qexe}
\P \left(P_{0,\omega}(X_{T_B}\in\partial_+B)
 \leq e^{-cL^{\beta_L+\epsilon_L}}\right) \leq \frac{1}{c_{3,11}}e^{-c_{3,11}L^{\beta_L}}
\end{equation}
where $B$ is the box defined in (\ref{lb}).
\end{proposition}

\medskip

\noindent Let us now prove Proposition \ref{pqexe}.
Let $\varsigma>0$.
 We will perform a one scale
renormalization analysis involving boxes
of side $\varsigma L^{\frac{\epsilon_L}{d+1}}$ which intersect the box $B$.
 Without loss of generality, we  assume that $e_1$ belongs to the
intersection of the half-spaces so that

\begin{equation}
\label{half1}
e_1\in\{x\in\mathbb Z^d: x\cdot l\ge 0\}
\end{equation}
and
\begin{equation}
\label{half2}
e_1\in\{x\in\mathbb Z^d: x\cdot \hat v\ge 0\}.
\end{equation}
Define the hyperplane perpendicular to direction $e_1$ as

\begin{equation}
\label{hiperplane}
H:=\{x\in\mathbb R^d: x\cdot e_1=0\}.
\end{equation}
We will need to work with the projection on the direction $l$ along the
hyperplane $H$ defined for $z\in\mathbb Z^d$ as
\begin{equation}
\label{projl}
P_l z:=\left(\frac{z\cdot e_1}{l\cdot e_1}\right)
\,l,
\end{equation}
and the projection of $z$ on $H$ along $l$ defined by

\begin{equation}
\label{projql}
Q_lz:=z-P_lz.
\end{equation}
Let $r>0$ be a fixed number which will eventually
be chosen large enough. For each $x \in\mathbb Z^d$ and $n$ define the
{\it mesoscopic box}

$$
D_n(x):=\{y\in \mathbb Z^d: -n< (y-x) \cdot e_1<n,\,
-rn\le |Q_l(y-x)|_\infty\le rn \},
$$
and their  front boundary

\begin{equation}
\nonumber
\partial_{+} D_{n}(x):=\{y \in \partial D_{n}(x): (y-x) \cdot e_1\ge n\}.
\end{equation}
Define the set of
mesoscopic boxes intersecting $B$ as

$$
\mathcal D:=\{D_n(x)\ {\rm with}\ x\in\mathbb Z^d: D_n(x)\cap B\ne\emptyset\}.
$$
 From now on, when there is no risk of
confusion, we will write $D$ instead of $D_n$ for a typical box in $\mathcal D$.
 Also, let us set $n:=\varsigma L^{\frac{\epsilon_L}{d+1}}$. We now say that a
box $D(x)\in\mathcal D$ is {\it good} if

\begin{equation}
\label{bbe1}
P_{x,\omega}(X_{T_{D(x)}}\in\partial_+ D(x))
\ge 1-\frac{1}{L}.
\end{equation}
Otherwise we will say that $D(x)$ is {\it bad}.

\medskip

\begin{lemma}
\label{bblemma}
Let $l\in\mathbb S^{d-1}$ and $M > 15d+5$.
Consider a  RWRE satisfying condition $(P)_M|l$
and the ellipticity condition $(E)_0$.
 Then, there is a $c_{3,5}$ such that for $r\ge c_{3,5}$ one has that

\begin{equation}
\label{bbestimate}
\limsup_{L\to\infty}L^{-\frac{\epsilon_L\gamma_L}{d+1}}
\log\mathbb P(D(0)\ {\rm is}\ {\rm bad})<0.
\end{equation}
\end{lemma}

\begin{proof}
By (\ref{bbe1}) and Markov inequality we have that
\begin{equation}
\label{ppp}
\P \left(D(0) \,\, {\rm is} \,\, {\rm bad}\right)
\leq \P \left(P_{0,\omega}(X_{T_{D(0)}}\not \in\partial_+ D(0))>
\frac{1}{L}\right)
\le
LP_0 \left( X_{T_{D(0)}}\not \in\partial_+ D(0)\right).
\end{equation}
Now, by Proposition \ref{PM-T0} of Section \ref{pol-t0}, we
know that the polynomial condition $(P)_M|l$
and the ellipticity condition $(E)_\alpha$
imply $(T)_{\gamma_L}|l$. But by Theorem \ref{sznitman-equivalence},
and the fact that $e_1$ is in the half spaces
determined by $l$ and $\hat v$ (see (\ref{half1}) and
(\ref{half2})), we can conclude that
$(T)_{\gamma_L}|l$ implies $(T)_{\gamma_L}|_{e_1}$.
On the other hand, it is straightforward to check that
there are constants $c_{3,5}$, $c_{3,6}>0$ such that for $r\ge c_{3,5}$, 
$(T)_{\gamma_L}|_{e_1}$ implies
that

$$
P_0 \left( X_{T_{D(0)}}\not \in\partial_+ D(0)\right)
\le \frac{1}{c_{3,6}}e^{-c_{3,6}L^{\frac{\epsilon_L\gamma_L}{d+1}}}.
$$
Substituting this back into inequality (\ref{ppp}) we see that
 (\ref{bbestimate})  follows.

\end{proof}

\noindent For each $m$ such that
$0\le m\le\left\lceil\frac{2L(l\cdot e_1)}{n}\right\rceil$ define the {\it block} $ R_m$  as
the collection of mesoscopic boxes (see Figure \ref{ecg1})

\begin{equation}
 R_m:=\{D(x) \in \mathcal{D}: \ {\rm for}\ {\rm some}\ x\ {\rm such}\ {\rm
  that}\ x\cdot e_1=nm\}.
\end{equation}

\begin{figure}[H]
\centering
\includegraphics[width=10cm]{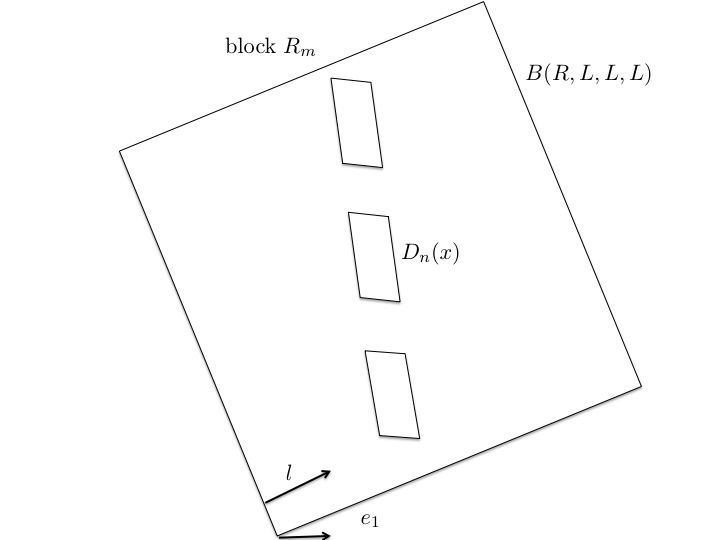}
\caption{A box $B$  with a set of inner boxes $D_n(x)$, which belong to a block $R_m$.}
\label{ecg1}
\end{figure}

\noindent The collection of these blocks is denoted by $\mathcal R$.
 We will say that a block $R_m$ is {\it good} if every box $D \in R_m$ is
 good. Otherwise, we will say that the block $R_m$ is {\it bad}.
 Now, for each $x \in R_m$ we associate a self-avoiding path $\pi^{(x)}$ such that
\begin{enumerate}
\item The path $\pi^{(x)}=(\pi^{(x)}_1,\ldots, \pi^{(x)}_{2n+1})$ has $2n$ steps.
\item $\pi^{(x)}_1=x$ and the end-point  $\pi^{(x)}_{2n+1} \in  D$ for some $D \in R_{m+1}$.
\item Whenever $D(x)$ does not intersect $\partial_+B$,
 the path $\pi^{(x)}$ is contained in $B$. Otherwise, the end-point $\pi_{2n+1}^{(x)}\in \partial_{+}B$.
\end{enumerate}

\noindent Define next $J$ as the total number of bad boxes of the collection
$\mathcal D$ and define

\begin{equation}
\nonumber
G_1:=\{\omega \in \Omega: J \leq L^{\beta_L+\frac{d}{d+1}\epsilon_L}\}.
\end{equation}
We will now denote by $\{m_1, \ldots m_N\}$ a generic subset of $\{0, \ldots, |\mathcal R|-1\}$ having $N$ elements. Let $\xi\in(0,1)$. Define

\begin{equation}
\nonumber
G_2:=\left\{\omega \in \Omega:
\sup_{N,\{m_1, \ldots, m_N\}}
\sum_{j=1}^N
\sup_{x_j\in R_{m_j}}
    \sum_{i=1}^{2n}
\log
 \frac{1}{\omega(\pi^{(x_j)}_i,
 \Delta \pi^{(x_j)}_i)}\leq 2n \log \left(\frac{1}{\xi}\right)L^{\beta_L+\frac{d}{d+1}\epsilon_L}
 \right\},
\end{equation}
where
the first supremum runs over $N\le L^{\beta_L+\frac{d}{d+1}\epsilon_L}$
and all subsets $\{m_1,\ldots,m_N\}$ of the set of blocks.
 Now, we can say that

\begin{equation}
\label{qe3.2}
\P \left(p_B
 \leq e^{-cL^{\beta_L+\epsilon_L}}\right) \leq \P \left(p_B
 \leq e^{-cL^{\beta_L+\epsilon_L}}, G_1 \cap G_2\right)+\P(G_1^c) + \P(G_2^c).
\end{equation}

\noindent Let us now show that the first term
 on the right-hand side of (\ref{qe3.2}) vanishes.
Indeed, on the event $G_1\cap G_2$, the probability
$p_B$ is bounded from below by the probability that the
random walk exits every mesoscopic box from its front side.
 Since $\omega \in G_1$, the random walk will have to do this
for at most $L^{\beta_L+\frac{d}{d+1}\epsilon_L}$ bad boxes. On each bad box $D(x)$ it will
follow the path $\pi^{(x)}$ defined above. But then on the event $G_2$,
we have a control on the product of the probability of traversing
all these paths through the bad boxes.
Hence, applying the strong Markov property and using the definition of good
box,  we conclude that  there is a $c_{3,7}>0$ such that for $0<\varsigma \leq c_{3,7}$ and
on the event $G_1\cap G_2$,

\begin{equation}
\nonumber
p_B \geq e^{-2L^{\beta_L+\epsilon_L}\varsigma \log \left(\frac{1}{\xi}\right)} \left(1-\frac{1}{L}\right)^L>e^{-cL^{\beta_L+\epsilon_L}}.
\end{equation}
Let us now estimate the term $\P(G_1^c)$ of
(\ref{qe3.2}).  Note first that the set $\mathcal D$ of mesoscopic boxes
can be divided into less than $2^dr^{d-1}\varsigma^d L^{\frac{d\epsilon_L}{d+1}}$ collections
of boxes, whose union is $\mathcal D$ and each collection has only
disjoint boxes. Let us call $M$ the number of such collections.
We also denote by
 $\mathcal D_i$ and $J_i$, where $1\le i\le M$, the $i$-th collection
and the number of bad boxes in such a collection respectively.   We then have that

\begin{equation}
\label{firstbound}
\P(G_1^c)\le \sum_{i=1}^M\P\left(J_i\ge \frac{1}{M}
L^{\beta_L+\frac{d}{d+1}\epsilon_L}\right).
\end{equation}
Now, by Chebyshev inequality

\begin{eqnarray}
\nonumber
&\P\left(J_i\ge \frac{1}{M}
L^{\beta_L+\frac{d}{d+1}\epsilon_L}\right)
\le e^{-\frac{L^{\beta_L+\frac{d}{d+1}\epsilon_L
}}{M}}\E[e^{J_i}]\\
\label{ge1}
&= e^{-\frac{L^{\beta_L+\frac{d}{d+1}\epsilon_L
}}{M}}\sum_{n=0}^{|\mathcal D_i|}
{|\mathcal D_i|\choose n} (ep_L)^n(1-ep_L)^{|\mathcal D_i|-n}
\left(\frac{1-p_L}{1-ep_L}\right)^{|\mathcal D_i|-n},
\end{eqnarray}
where $p_L$ is the probability that a box is bad. Now the
last factor of each term after the summation of the right-hand
side of (\ref{ge1}) is bounded by

$$
\left(\frac{1-p_L}{1-ep_L}\right)^{|\mathcal D_i|},
$$
which clearly tends to $1$ as $L\to\infty$ by the fact that
$|\mathcal D_i|\le c_{3,8}L^d$, the definition of $\epsilon_L$ and by Lemma \ref{bblemma} for some $c_{3,8}>0$. Thus,
there is a constant $c_{3,9}>0$ such that
$$
\P\left(J_i\ge \frac{1}{M}
L^{\beta_L+\frac{d}{d+1}\epsilon_L
}\right)
\le c_{3,9}e^{-\frac{L^{\beta_L+\frac{d}{d+1}\epsilon_L
}}{M}}.
$$
Substituting this back into (\ref{firstbound}) we hence see that

\begin{equation}
\label{first1}
\P\left(G_1^c\right)
\le c_{3,9}(2\varsigma)^d r^{d-1}L^{\frac{d}{d+1}\epsilon_L}
e^{-\frac{L^{\beta_L}}{(2\varsigma)^d r^{d-1}}}.
\end{equation}
Let us now bound the term $\P(G_2^c)$ of (\ref{qe3.2}).
Define $\beta'_L:=\beta_L+\frac{d}{d+1}\epsilon_L$.
Note that for each $0<\alpha<\bar\alpha$ one has that

\begin{eqnarray*}
&\P\left(G_2^c\right) \leq \sum_{N=1}^{L^{\beta_L'}} \P \left(\exists \, \{m_1, \ldots, m_N\} \,\, \mbox{and}\,\, x_j \in R_{m_j} \, \mbox{such} \,\,\mbox{that} \right. \nonumber \\
&\left.\sum_{j=1}^
    {N}
    \sum_{i=1}^{2n}\log
 \frac{1}{\omega \left(\pi_i^{(x_j)},
 \Delta
      \pi_i^{(x_j)}\right)}  > 2n \log\left(\frac{1}{\xi}
\right) L^{\beta_L'} \right) \nonumber \\
\nonumber
      &\leq
\sum_{N=1}^{L^{\beta_L'}}{|\mathcal R|\choose N}
r^{d-1}(2\varsigma )^{L^{\beta_L'}} e^{(\log L)\frac{\epsilon_L}{d+1} L^{\beta_L'}}
e^{(\log\eta_\alpha)2nL^{\beta_L'}-2\alpha n \log\left(\frac{1}{\xi}\right) L^{\beta_L'}} \nonumber \\
      &\leq
L^{\beta_L'}\lceil
\frac{2L(l\cdot e_1)}{n}\rceil^{L^{\beta_L'}}
r^{d-1}(2\varsigma )^{L^{\beta_L'}} e^{(\log L)\frac{\epsilon_L}{d+1} L^{\beta_L'}}
e^{(\log\eta_\alpha)2nL^{\beta_L'}-2\alpha n \log\left(\frac{1}{\xi}\right) L^{\beta_L'}},
\end{eqnarray*}
 where in the last line, we have used that $|{\cal R}|$ can be bounded by $\lceil
\frac{2L(l\cdot e_1)}{n}\rceil$. 
It now follows that for $\xi$ such that
$\log\left(\frac{1}{\xi^{2\alpha}\eta_\alpha^3}\right)>0$ one can find
a constant $c_{3,10}$ such that

\begin{equation}
\label{second2}
\P\left(G_2^c\right) \leq \frac{1}{c_{3,10}}e^{-c_{3,10}L^{\beta_L+\epsilon_L}}.
\end{equation}
Substituting back (\ref{first1}) and (\ref{second2}) into (\ref{qe3.2})
we end up the proof of Proposition \ref{pqexe}.

\subsection{Condition $(T)_{\gamma_L}$ implies the effective criterion}
\label{aqee-ec}

Here we will introduce a generalization of the effective criterion
introduced by Sznitman in \cite{Sz02} for RWRE, dropping
the assumption of uniform
ellipticity and replacing it by
the ellipticity condition $(E)_0$.
Let $l \in \mathbb S^{d-1}$ and $d \geq 2$.
 We will say that the {\it effective criterion in direction $l$} holds if

\begin{equation}
\label{ec}
c_2(d)
\inf_{L\ge c_3, 3\sqrt{d}\le\tilde L<L^3}\inf_{\alpha>0}\inf_{0<a\le \alpha}\left\{
\Upsilon^{3(d-1)}
 \tilde L^{d-1}L^{3(d-1)+1}\E[\rho^a_B]
\right\}<1,
\end{equation}
where

\begin{equation}
\label{bspec}
B:=B(R,L-2,L+2,\tilde L) \,\, {\rm and} \,\, \Upsilon:=\max\left\{\frac{\alpha}{24}, \left(\frac{2c_1}{c_1-1}\right)\log \eta_{\alpha}^2\right\},
\end{equation}
while $c_2(d)$ and $c_3(d)$ are dimension dependent constants
that will be introduced in
subsection \ref{ec-tp}.
Note that in particular, the effective criterion in direction $l$
implies that
condition $(E)_0$ is satisfied.
Here we will prove the following proposition.

\medskip

\begin{proposition}
\label{prop-waqec}
Let $l\in\mathbb S^{d-1}$. Assume that the ellipticity condition
$(E)_0$ and that $(T)_{\gamma_L}|l$  are fulfilled.
Then, the effective criterion in direction $l$ is satisfied.
\end{proposition}

\medskip

\noindent To prove Proposition \ref{prop-waqec}, we begin defining the following quantities

\begin{equation}
\nonumber
\beta_1(L):=\frac{\gamma_L}{2}=\frac{\log 2}{2 \log \log L}
\end{equation}

\begin{equation}
\nonumber
\sigma(L):=\frac{\gamma_L}{3}=\frac{\log 2}{3 \log \log L}
\end{equation}

\begin{equation}
\nonumber
a:=L^{-\sigma(L)}.
\end{equation}

\noindent We will write $\rho$ instead of $\rho_{B}$, where $B$
is the box defined in (\ref{bspec}) (see \ref{box1}) with $\tilde L=L^3$.
 Following  \cite{BDR12}, it is convenient to split $\E \rho^a$ according to

\begin{equation}
\label{split}
\E \rho^a=\mathcal{E}_0+\sum_{j=1}^{n-1}
\mathcal{E}_j+\mathcal{E}_{n}
\end{equation}
where

\begin{equation}
\nonumber
n:=n(L):=\left \lceil \frac{4(1-\gamma_L/2)}{\gamma_L}\right \rceil +1,
\end{equation}

\begin{equation}
\nonumber
\mathcal{E}_0:=\E \left(\rho^a, p_B >e^{-cL^{\beta_1}} \right),
\end{equation}

\begin{equation}
\nonumber
\mathcal{E}_j:=\E \left(\rho^a,e^{-cL^{\beta_{j+1}}}< p_B \leq e^{-cL^{\beta_j}} \right)
\end{equation}
for $j \in \{1, \ldots, n-1\}$, and

\begin{equation}
\nonumber
\mathcal{E}_{n}:=\E \left(\rho^a, p_B \leq e^{-cL^{\beta_{n}}} \right)
\end{equation}
with parameters

\begin{equation}
\nonumber
\beta_j(L):=\beta_1(L)+(j-1)\frac{\gamma_L}{4},
\end{equation}
for $2 \leq j \leq n(L)$.
We will now estimate each of the $n$ terms
appearing in (\ref{split}).
For the first $n-1$ terms, we now state two lemmas proved by
 Berger, Drewitz and Ram\'\i rez in \cite{BDR12}, whose proofs we omit.
 The
following lemma is a consequence of Jensen's inequality.

\medskip

\begin{lemma}
\label{split1}
Assume that $(T)_{\gamma_L}$ is satisfied. Then

\begin{equation}
\nonumber
\mathcal{E}_0 \leq e^{cL^{\frac{\gamma_L}{6}}-L^{\frac{2 } {3}\gamma_L(1+o(1))}}
\end{equation}
as $L \to \infty$.
\end{lemma}

\medskip

\noindent The second lemma  follows from Proposition \ref{pqexe}.

\medskip

\begin{lemma}
\label{split2}
Assume that the weak atypical quenched exit estimate
(\ref{qexe}) is satisfied.
 Then there exists a constant $c_{3,12}>0$ such that for all $L$ large enough and all
 $j \in \{1, \ldots, n-1\}$ one has that

\begin{equation}
\nonumber
\mathcal{E}_j \leq \frac{1}{c_{3,12}}e^{cL^{\left(\frac{1}{6}+\frac{j}{4}\right)\gamma_L}
-c_{3,12}L^{\left(\frac{1}{4}+\frac{j}{4}\right)\gamma_L-\epsilon_L}}.
\end{equation}
\end{lemma}

\medskip

\noindent In \cite{BDR12}, where
it is assumed that the environment is uniformly elliptic,
one has that $\mathcal{E}_n=0$. Nevertheless, since here we are not
assuming uniform ellipticity this is not the case.

\begin{lemma}
\label{split3}
Assume that $(E)_0$ and $(T)_{\gamma_L}$ are satisfied. Then there exists a constant $c_{3,16}>0$ such that for all $L$ large enough we have
\begin{equation}
\nonumber
\mathcal{E}_n \leq \frac{1}{c_{3,16}}e^{-c_{3,16}L^{1-\epsilon_L}}.
\end{equation}

\end{lemma}

\begin{proof}
Choose $0<\alpha<\bar\alpha$.
Consider
a nearest neighbor 
path $(x_1,\ldots,x_m)$  from $0$ to $\partial_+ B$,
so that $x_1=0$ and $x_m\in\partial_+ B$,
$x_1,\ldots,x_{m-1} \in B$ and
which has the minimal number of steps $m$ (note that
it is then self-avoiding).
Then,

\begin{eqnarray}
\nonumber
&\E \left[\rho^a, p_B \leq e^{-cL^{\beta_n}}\right] \leq \E \left[e^{
\frac{\alpha}{2}
    \sum_{1}^{m}\log\frac{1}{\omega(x_i,\Delta
      x_i)}},\sum_{1}^{m}\log\frac{1}{\omega(x_i,\Delta x_i)}>
\frac{3m}{\alpha}\log
 \eta_\alpha\right]\\
\label{rhoa}
      & +\E \left[e^{a
    \sum_{1}^{m}\log\frac{1}{\omega(x_i,\Delta
      x_i)}},\sum_{1}^{m}\log\frac{1}{\omega(x_i,\Delta x_i)} \leq
\frac{3m}{\alpha}
\log \eta_\alpha,p_B \leq e^{-cL^{\beta_n}}\right],
         \end{eqnarray}
where in the first line, we have used that for any $\alpha>0$,
 $a \leq \frac{\alpha}{2}$ for $L$ large. Now, using Cauchy-Schwartz
inequality, Chebyshev inequality, (\ref{qexe}) and the
fact that the probability of exiting a square box through its front
side is smaller than the probability of exiting a box of width $L^3$ and
length $L$ through its front side, we can see that the
right-hand side of (\ref{rhoa}) is smaller than 

\begin{eqnarray}
\nonumber
& \E \left[e^{\alpha
    \sum_{1}^{m}\log\frac{1}{\omega(x_i,\Delta
      x_i)}}\right]^{1/2} \P \left(\sum_{1}^{m}\log\frac{1}{\omega(x_i,\Delta
    x_i)}
>\frac{3m}{\alpha}\log  \eta_\alpha\right)^{1/2}
      + e^{\frac{3a m}{\alpha} \log \eta_\alpha}\P \left(p_B \leq e^{-cL^{\beta_n}}\right)
      \\
\label{second-estimate}
      & \leq e^{-m\log \eta_{\alpha}}+\frac{1}{c_{3,13}}
e^{\frac{3am}{\alpha} \log\eta_\alpha-c_{3,13}L^{\beta_n(L)-
      \epsilon_L}},
\end{eqnarray}
for some constant $c_{3,13}>0$.
Now, using the fact that there are constants
$c_{3,14}$ and $c_{3,15}$ such that

$$
c_{3,14}L\le m\le c_{3,15}L,
$$
we can substitute (\ref{second-estimate}) into
(\ref{rhoa}) to conclude that there is
a constant $c_{3,16}$ such that

$$
\E \left[\rho^a, p_B \leq e^{-cL^{\beta_n}}\right] \leq
\frac{1}{c_{3,16}}e^{-c_{3,16}L^{1-\epsilon_L}}.
$$

\end{proof}

\medskip

It is now straightforward to conclude the proof
of Proposition \ref{prop-waqec}
using the estimates of Lemmas \ref{split1}, \ref{split2}
and \ref{split3}.

\subsection{The effective criterion implies $(T')$}
\label{ec-tp}

%Here we will introduce a generalization of the effective criterion
%introduced by Szntiman in \cite{Sz02} for RWRE in uniformly
%elliptic environments, to RWRE in elliptic environments.
We will prove that the generalized effective criterion
and the ellipticity condition $(E)_0$
imply $(T')$. To do this, it is enough to prove the following.

\medskip

\begin{proposition}
\label{pec}
Throughout choose $0<\alpha<\bar\alpha$. Let $l\in\mathbb S^{d-1}$ and $d\ge 2$. If the effective criterion in direction $l$ holds then
there exists a constant $c_{3,28}>0$ and a
neighborhood $V_l$ of direction $l$ such that
for all $l'\in V_l$ one has that
\begin{equation}
\nonumber
\overline{\lim_{L \to \infty}}L^{-1}e^{c_{3,28}(\log L)^{1/2}}\log P_0
\left[\tilde{T}^{l'}_{-\tilde{b}L}<T^{l'}_{bL}\right]<0, \,\,\, \mbox{for}\,\, \mbox{all}\,\,\, b, \tilde{b}>0.
\end{equation}
In particular, if (\ref{ec}) is satisfied, condition $(T')|l$ is satisfied.
\end{proposition}

%The condition defined in inequality (\ref{ec})
%of Proposition \ref{pec}, is called the {\it effective criterion
%in direction $l$}.

\noindent To prove this proposition,
we will follow the same strategy used by Sznitman in \cite{Sz02}
to prove Proposition 2.3 of that paper under
 the assumption of uniform ellipticity.
Firstly we need to define some constants. Let

$$
c_1'(d,\alpha):= 13+\frac{24d}{\alpha}+\frac{24d+12\log\eta_\alpha}{2 \log\eta_\alpha},
$$

$$
c_2'(d,\alpha):= c_1 c_1',
$$
and

$$
c_4(d,\alpha):=\frac{48 c_2'}{\alpha},
$$
where $c_1$ is defined in (\ref{c1}).
Define for $k\ge 0$ the sequence $\{N_k:k\ge 0\}$ by

\begin{equation}
\label{specifications 0}
N_k:=\frac{c_4}{u_0}8^k,
\end{equation}
where $u_0 \in (0,1)$.
Let $L_0,\tilde L_0, L_1$ and $\tilde L_1$ be constants such that

\begin{equation}
\label{l0l1}
3\sqrt{d} \leq \tilde{L}_0\leq L_0^3,\quad L_1=N_0L_0\quad {\rm and}\
\tilde L_1=N_0^3\tilde L_0.
\end{equation}
Now, for $k\ge 0$ define recursively the sequences
$\{L_k:k\ge 0\}$ and $\{\tilde L_k:k\ge 0\}$ by

\begin{equation}
\label{specifications 1}
 L_{k+1}:=N_kL_k, \quad {\rm and}\quad
\tilde L_{k+1}:=N_k^3\tilde L_k.
\end{equation}
 It is straightforward to see that
for each $k\ge 1$
\begin{equation}
\label{specifications 2}
L_k=\left(\frac{c_4}{u_0}\right)^k \,8^{\frac{k(k-1)}{2}}L_0, \quad \tilde{L}_k=\left(\frac{L_k}{L_0}\right)^3 \tilde{L}_0.
\end{equation}
Furthermore, we also consider for $k\ge 0$ the box
\begin{equation}
\nonumber
B_k:=B(R, L_k-1,L_k+1,\tilde L_k),
\end{equation}
and the positive part of its boundary $\partial_+B_k$,
and will use the notations $\rho_k=\rho_{B_k}, p_k=p_{B_k}, q_k=q_{B_k}$ and
$n_k=[N_k]$. Following  Sznitman \cite{Sz02}, we
introduce for each $i\in\mathbb Z$
\begin{equation}
\label{function H}
\h_i:=\{x \in \Z^d,\, \exists \, x'\in \Z^d, \,|x-x'|=1,\, (x \cdot l-iL_0)\,(x' \cdot l-iL_0) \leq 0\}.
\end{equation}
We also define the function $I:\mathbb Z^d\to\mathbb Z$ by

$$I(x):=i, \ \mbox{for} \ x \ \mbox{such}\ \mbox{that}\  x\cdot l \in
\left[iL_0-\frac{L_0}{2}, iL_0+\frac{L_0}{2}\right).
$$
 Consider now  the successive times of visits of the random walk
to the sets $\{\h_i:i\in\mathbb Z\}$,  defined recursively as

$$
V_0:=0, \quad V_1:=\inf\{n \geq 0: \, X_n \in \h_{I(X_0)+1} \cup
\h_{I(X_0)-1}\}
$$
and
$$
V_{k+1}:=V_k+V_1 \circ \theta_{V_k}, \,\, k \geq 0.
$$
For $\omega \in \Omega$, $x \in \Z^d$, $i \in \Z$, let
\begin{equation}
\label{exit prob}
\widehat{q}(x,\omega):=P_{x, \omega}[X_{V_1} \in \h_{I(x)-1}]
\end{equation}
while $\widehat{p}(x, \omega):=1-\widehat{q}(x, \omega)$,
and
\begin{equation}
\label{rho function}
\widehat{\rho}(i,\omega):=\sup \left\{\frac{\widehat{q}(x,
  \omega)}{\widehat{p}(x,\omega)}: x \in \h_i, \; \sup_{2\le j\le d}|R(e_j)\cdot x| <
\tilde{L}_1 \right\}.
\end{equation}
We consider also the stopping time
$$
\tilde{T}:=\inf\left\{n \geq 0: \sup_{2\le j\le d}|X_n \cdot R(e_j)|\geq
\tilde{L}_1\right\},
$$
and the function $f:\{n_0+2, n_0+1, \ldots \} \times \Omega\to\mathbb R$
defined by

\begin{equation}
\nonumber
f(n_0+2,\omega):=0,\quad
f(i,\omega):=\sum_{ m=i}^{n_0+1}\prod_{j=m+1}^{n_0+1} \widehat{\rho}(j, \omega)^{-1}, \quad \mbox{for} \quad i \leq n_0+1.
\end{equation}
We will frequently write $f(n)$ instead $f(n,\omega)$.
 Let us now proceed to prove Proposition \ref{pec}. The following proposition corresponds to the first step
in an induction argument which will be used to prove
Proposition \ref{pec}.

\medskip

\begin{proposition}
\label{prop 1}
 Let $\alpha>0$.
Let $L_0, L_1,\tilde L_0 $ and $\tilde L_1$ be constants
satisfying (\ref{l0l1}), with $N_0\geq 7$.
Then, there exist  $c_{3,17}, c_{3,18}(d), c_{3,19}(d)>0$ such that
 for  $L_0 \geq c_{3,17}$,
 $a \in (0,\min\{\alpha,1\}]$,
$u_0 \in [\xi^{L_0/d},1]$,
 $0<\xi<\frac{1}{\eta_\alpha^{2/\alpha}}$
and
\begin{equation}
\label{induction-in}
N_0\le \frac{1}{L_0}\left(\frac{e}{\xi}\right)^{L_0},
\end{equation}
the following is satisfied
\begin{eqnarray}
\nonumber
&\E[\rho_1^{a/2}]
\leq c_{3, 18} \left\{ \xi^{-c_2'L_1}\left(c_{3,19}
\tilde{L}_1^{(d-2)}\frac{L_1^3}{L_0^2}\tilde{L_0}\E[q_0]
\right)^{\frac{\tilde{L}_1}{12N_0 \tilde{L}_0}}\right.\\
\label{estimate 1}
&\left. + \sum_{m=0}^{ N_0+1} \left(c_{3,19} \tilde{L}_1^{(d-1)}\E[\rho_0^a]\right)^{\frac{N_0+m-1}{2}}+e^{-c_{1}L_1\log\frac1{\xi^\alpha\eta_\alpha^2}}\right\}.
\end{eqnarray}
\end{proposition}

\noindent \begin{proof}
The following inequality   is stated
and proved in \cite{Sz02} by Sznitman without using
any kind of uniform ellipticity assumption (inequality (2.18) in
\cite{Sz02}).
For  every $\omega \in \Omega$
\begin{equation}
\label{chung-formula}
P_{0,\omega}\left(\tilde{T}^l_{1-L_1}<\tilde{T} \wedge T^l_{L_1+1}\right)\leq \frac{f(0)}{f(1-n_0)}.
\end{equation}
\noindent Consider now the event
\begin{equation}
\nonumber
G:=\{\omega: P_{0,\omega}\left(\tilde{T}\leq \tilde{T}^l_{1-L_1} \wedge T^l_{L_1+1}\right) \leq \xi^{(c_1'-1)c_1L_1}\},
\end{equation}
and write

\begin{equation}
\label{dec-g}
\E[\rho_1^{a/2}]=\E[\rho_1^{a/2},G]+\E[\rho_1^{a/2},G^c].
\end{equation}
The first term $\E[\rho_1^{a/2},G]$ of (\ref{dec-g}),
can in turn be decomposed as

\begin{equation}
\label{dec-ga}
\E[\rho_1^{a/2},G]=\E[\rho_1^{a/2},G,A_1]+\E[\rho_1^{a/2},G,A_1^c],
\end{equation}
where we have defined

$$
A_1:= \{\omega \in \Omega:f(2-n_0)-f(0) \geq f(1-n_0) \xi^{(c_1'-1)c_1L_1},
f(0) \geq f(1-n_0)\xi^{(c_1'-1)c_1L_1}\}.
$$
Furthermore, note that

$$
A_1^c\subset A_2\cup A_3,
$$
where

\begin{eqnarray*}
&A_2:=\{\omega\in\Omega: f(2-n_0)-f(0)<f(1-n_0)\xi^{(c_1'-1)c_1 L_1}\},\ {\rm while} \\
&A_3:=\{\omega\in\Omega: f(0)<f(1-n_0)\xi^{(c_1'-1)c_1L_1}\}.
\end{eqnarray*}
Therefore,

\begin{equation}
\label{dec-a2a3}
\E[\rho_1^{a/2}]\le \E[\rho_1^{a/2},G,A_1]+\E[\rho_1^{a/2},A_2]
+\E[\rho_1^{a/2},G, A_3]+\E[\rho_1^{a/2},G^c].
\end{equation}
We now subdivide the rest of the proof in several
steps corresponding to an estimation for each one of the terms in
inequality (\ref{dec-a2a3}).

\smallskip

\noindent {\it Step 1: estimate of $\E[\rho_1^{a/2},G,A_1]$}.
Here we estimate the first term of display (\ref{dec-ga}).
To do this, we can follow the argument presented by Sznitman
in Section 2 of \cite{Sz02}, to prove that inequality
(\ref{chung-formula}) implies that there exist constant
$c_{3,20}(d)$ such that
    \begin{equation}
    \label{first estimate}
    \E \left[\rho_1^{a/2},G, A_1\right] \leq 2\sum_{m=0}^{n_0+1} \left(c_{3,20}(d)\tilde{L}_1^{(d-1)}\E[\rho^a_0]\right)^{\frac{n_0+m-1}{2}}.
    \end{equation}
Indeed on $G \cap A_1$ and with the help of (\ref{chung-formula}) one gets that
\begin{eqnarray}
\nonumber
&\rho_1 = \frac{P_{0,\omega}[\tilde T^l_{-L_1+1}<\tilde T \wedge T^l_{L_1+1}]+P_{0,\omega}[\tilde T \leq \tilde T^l_{-L_1+1}\wedge T^l_{L_1+1}]}{1-P_{0,\omega}[\tilde T^l_{-L_1+1}<\tilde T \wedge T^l_{L_1+1}]-P_{0,\omega}[\tilde T \leq \tilde T^l_{-L_1+1}\wedge T^l_{L_1+1}]}\\
\label{rho1}
&\leq
\frac{f(0)+f(1-n_0)\xi^{(c_1'-1)c_1L_1}}{\left(f(1-n_0)-f(0)-f(1-n_0)\xi^{(c_1'-1)c_1L_1}\right)_+}
 \leq  \frac{2f(0)}{\left(f(1-n_0)-f(0)-f(1-n_0)\xi^{(c_1'-1)c_1L_1}\right)_+},
\end{eqnarray}
where in the first inequality we have used the fact that $\omega\in G$,
while in the second that $\omega\in A_1$.
Regarding the term in the denominator in the last expression, we can use the definition of the function $f$ and obtain
\begin{eqnarray*}
& f(1-n_0)-f(0)-f(1-n_0)\xi^{(c_1'-1)c_1L_1}\\
&= \prod_{j=2-n_0}^{n_0+1}\widehat \rho(j,\omega)^{-1}+
f(2-n_0)-f(0)-f(1-n_0)\xi^{(c_1'-1)c_1L_1}\\
& \geq  \prod_{j=2-n_0}^{n_0+1}\widehat \rho(j,\omega)^{-1},
\end{eqnarray*}
where we have used that $\omega \in A_1$ in the last inequality. Substituting
this estimate in (\ref{rho1}), we conclude that for
  $\omega \in G \cap A_1$ one has that
\begin{equation}
\label{rho2}
\rho_1 \leq 2 \prod_{j=2-n_0}^{ n_0+1}\widehat \rho(j,\omega)f(0)=2
\sum_{m=0}^{n_0+1}
 \prod_{j=2-n_0}^{m}\widehat \rho(j,\omega).
\end{equation}
At this point, using (\ref{rho2}), the fact that since $a\in [0,1]$
one has that $(u+v)^{a/2} \leq
u^{a/2}+v^{a/2}$ for $u,v \geq 0$,
the fact that $\{\widehat
\rho(j,\omega), j \,\mbox{even} \}$ and $\{\widehat \rho(j,\omega), j
\,\mbox{odd} \}$ are two collections of independent random variables
 and the Cauchy-Schwartz's inequality, we can assert that
\begin{eqnarray*}
& &\E[\rho_1(\omega)^{a/2},G, A_1]\\
& \leq & 2 \sum_{0 \leq m \leq n_0+1}\E\left[\prod_{1-n_0 <j\leq m} \widehat \rho (j,\omega)^{a/2}\right]\\
& \leq & 2 \sum_{0 \leq m \leq n_0+1}\E\left[\prod_{\substack{1-n_0 <j\leq m\\ j \, is \, even}} \widehat \rho (j,\omega)^{a}\right]^{1/2}\E\left[\prod_{\substack{1-n_0 <j\leq m\\ j \, is \, odd}} \widehat \rho (j,\omega)^{a}\right]^{1/2}\\
&=& 2 \sum_{0 \leq m \leq n_0+1}\quad \prod_{1-n_0 <j\leq m} \E\left[\widehat \rho (j,\omega)^{a}\right]^{1/2}.
  \end{eqnarray*}

\medskip

\noindent In view of (\ref{exit prob}) one gets easily that for
$i\in\mathbb Z$ and $x \in \mathcal{H}_i$,

$$
\widehat p(x,\omega) \geq p_0 \circ t_x(\omega),
$$
where the canonical shift $\{t_x:x\in\mathbb Z^d\}$ has been
defined in (\ref{shift}).
Hence, for $i \in \Z$ and $x \in \mathcal{H}_i$,

$$\frac{\widehat q(x)}{\widehat p(x)}\leq \rho_0 \circ t_x.$$
 Following Sznitman \cite{Sz02} with the help of (\ref{function H})
the estimate (\ref{first estimate}) follows.

\smallskip

\noindent {\it Step 2: estimate of $\E[\rho_1^{a/2},A_2]$}.
Here we will prove the following estimate for the second term of
inequality (\ref{dec-a2a3}),

$$
\E[\rho^{a/2},A_2]\le 4 e^{-c_1L_1\log\frac{1}{\xi^{\alpha} \eta_\alpha^2}}.
$$
By the definition of $c_1$ (see (\ref{c1})), we know that
necessarily there exists a path with less than
$c_1(L_1+1+\sqrt{d})$ steps between the origin and $\partial_+B_1$.
 Therefore, for $L_0\ge 1+\sqrt{d}$,
there is a nearest neighbor self-avoiding path $(x_1,\ldots,x_n)$ with $n$ steps
from the origin to $\partial_+B_1$, such that
$2c_1L_1\le n\le 2c_1L_1+1$, $x_1, \ldots, x_n \in B_1$ and
$x_n \cdot l\ge L_1+1$. Thus,
for every $r\ge 0$ we have that

\begin{equation}
\label{rest}
\rho_1^r \leq \frac{1}{p_1^r} \leq e^{r\sum_{i=1}^{n}\log \frac{1}{\omega(x_i,
    \Delta x_i)}},
\end{equation}
where $\Delta x_i:=x_{i+1}-x_i$ for $1\le i\le n-1$ as defined in (\ref{delta}).
 We then have applying inequality (\ref{rest}) with $r=a/2$ that
    \begin{eqnarray}
\nonumber
 &   \E\left[\rho_1^{a/2}, A_2\right]  \leq  \E
\left[e^{\alpha/2 \sum_{1}^{n}\log \frac{1}{\omega(x_i, \Delta x_i)}}, A_2,
\sum_{1}^{n}\log \frac{1}{\omega(x_i, \Delta x_i)} \leq n
\log \left(\frac{1}{\xi}\right)
 \right]
 \\
    \label{estimate 2a}
  &+ \E \left[e^{\alpha/2 \sum_{1}^{n}\log \frac{1}{\omega(x_i, \Delta x_i)}},
\sum_{1}^{n}\log \frac{1}{\omega(x_i, \Delta x_i)} > n
\log \left(\frac{1}{\xi}\right)\right].
    \end{eqnarray}

\noindent Regarding the second term of the right side of (\ref{estimate 2a}), we can apply
the Cauchy-Schwarz inequality,  the exponential Chebyshev inequality and use the fact that the jump probabilities
$\{\omega(x_i,\Delta_i):1\le i\le n-1\}$ are independent
to conclude that

\begin{eqnarray}
\nonumber
& \E \left[e^{\alpha \sum_{1}^{n}\log \frac{1}{\omega(x_i, \Delta x_i)}}\right]^{1/2}
\P\left(\sum_{1}^{n}\log \frac{1}{\omega(x_i, \Delta x_i)} > n
\log \left(\frac{1}{\xi}\right)\right)^{1/2}
\\
   \label{estimate 2b}
&\leq e^{\left(2 \log \eta_{\alpha}-\alpha \log\left(\frac{1}{\xi}\right)\right)n/2}.
\end{eqnarray}

\noindent Meanwhile, note that  the first term on the right side of (\ref{estimate 2a}) can be bounded by
\begin{equation}
\label{estimate 2c}
e^{\frac{\alpha \log\left(\frac{1}{\xi}\right)n}{2}}\P(A_2).
\end{equation}
Hence, we need an adequate estimate for $\P(A_2)$. Now,
\begin{eqnarray}
\nonumber
&\P(A_2) = \P \left(f(2-n_0)-f(0)<
  e^{\frac{-(c_1'-1)c_1L_1\log\left(\frac{1}{\xi}\right)}{2}},A_2\right)\\
\nonumber
&+\P
  \left(f(2-n_0)-f(0)\geq e^{\frac{-(c_1'-1)c_1L_1\log(\frac{1}{\xi})}{2}},A_2
  \right)\\
\label{similar}
& \leq \P \left(f(2-n_0)-f(0)< e^{\frac{-(c_1'-1)c_1L_1\log\left(\frac{1}{\xi}\right)}{2}}\right)+\P \left(f(1-n_0) >  e^{\frac{(c_1'-1)c_1L_1\log(\frac{1}{\xi})}{2}}\right).
\end{eqnarray}
The two terms in the rightmost side of display (\ref{similar})
will be estimated by similar methods: in both cases, we will
 use the fact that $\{\widehat
\rho(j,\omega), j \,\mbox{even} \}$ and $\{\widehat \rho(j,\omega), j
\,\mbox{odd} \}$ are two collections of independent random variables, the
Cauchy-Schwartz's inequality and the Chebyshev inequality. Specifically for
the first
term of the rightmost side of (\ref{similar}) we have that

\begin{eqnarray}
\nonumber
&  \P \left(f(2-n_0)-f(0)<
 e^{\frac{-(c_1'-1)c_1L_1\log\left(\frac{1}{\xi}\right)}{2}}\right)
  \leq \P\left(\prod_{j=0}^{n_0+1}\widehat{\rho}(j,\omega)^{-1}<
e^{\frac{-(c_1'-1)c_1L_1\log\left(\frac{1}{\xi}\right)}{2}} \right)\\
\nonumber
&  =\P \left(\prod_{j=0}^{n_0+1}\widehat{\rho}(j,\omega)^{\alpha/2}>
e^{\frac{(c_1'-1)c_1L_1 \alpha \log\left(\frac{1}{\xi}\right)}{4}} \right)\\
\nonumber
&  \leq
 e^{\frac{-(c_1'-1)c_1L_1 \alpha \log\left(\frac{1}{\xi}\right)}{4}}
  \E\left[\prod_{\substack{j=1,\\ j\,
        odd}}^{n_0+1}\widehat{\rho}(j,\omega)^\alpha \right]^{1/2}
  \E\left[\prod_{\substack{j=0,\\ j\,
        even}}^{n_0+1}\widehat{\rho}(j,\omega)^\alpha \right]^{1/2}\\
\label{f2n}
& = e^{\frac{-(c_1'-1)c_1L_1 \alpha \log\left(\frac{1}{\xi}\right)}{4}} \prod_{j=0}^{n_0+1}\E\left[\widehat{\rho}(j,\omega)^\alpha \right]^{1/2}.
\end{eqnarray}
By an estimate analogous to (\ref{rest}), we know that for $L_0\ge 1+\sqrt{d}$,
 for each
$j \in \{0, \ldots, n_0+1\}$ and each $x\in\mathcal H_j$,
 there exists a nearest neighbor self-avoiding path
$(y_1,\ldots,y_m)$ with
$m$ steps, such that $2c_1L_0\le m\le 2c_1L_0+1$, between $x$ and ${\cal{H}}_{j+1}$.
Also, $y_1 \cdot l, \ldots, y_{m-1}\cdot l \in (1-L_0, L_0+1)$ and
$y_m\cdot l \geq L_0+1$.
Then, in view of (\ref{specifications 1}), (\ref{specifications 2}), (\ref{exit prob}) and (\ref{rho function}), we have that for each $j \in \{0,\ldots, n_0+1\}$

\begin{eqnarray}
\nonumber
&\E\left[\widehat{\rho}(j,\omega)^{\alpha}\right]^{1/2} \leq
\sum
\E\left[\widehat{p}(x,\omega)^{-\alpha}\right]^{1/2} \\
\label{estimate 2d aux}
&\leq  2L_1^{3(d-1)}\E \left[e^{\alpha \sum_{1}^{m}\log \frac{1}{\omega(y_i,
\Delta y_i)}}\right]^{1/2}
 \leq 2L_1^{3(d-1)}e^{\frac{m\log \eta_\alpha}{2}},
\end{eqnarray}
where the summation goes over all $x \in {\cal{H}}_j$
such that $ \sup_{2 \leq i \leq d}|R(e_i)\cdot x|<\widetilde{L}_1$.
Substituting the estimate (\ref{estimate 2d aux}) back into
(\ref{f2n}) we see that

\begin{eqnarray}
\nonumber
&\P \left(f(2-n_0)-f(0)<e^{\frac{-(c_1'-1)c_1L_1\log\left(\frac{1}{\xi}\right)}{2}}\right) \leq e^{\frac{-(c_1'-1)c_1L_1 \alpha \log\left(\frac{1}{\xi}\right)}{4}} 2^{(n_0+2)} L_1^{3(d-1)(n_0+2)}e^{\frac{(\log \eta_\alpha)m(n_0+2)}{2}} \\
\nonumber
&\leq  e^{\frac{-(c_1'-1) c_1L_1\alpha \log\left(\frac{1}{\xi}\right)}{4}+\log 2(n_0+2)+3(d-1)\frac{\log L_0N_0}{L_0}L_0(n_0+2)+(\log \eta_{\alpha})\frac{(2c_1L_0+1)(n_0+2)}{2}}\\
\label{estimate 2d}
& \leq e^{-L_1 \left( \frac{(c_1'-1)\alpha \log\left(\frac{1}{\xi}\right)c_1}{4} -1 -6(d-1) \left(1+\log \left(\frac{1}{\xi}\right) \right)-\log \eta_{\alpha}(3c_1+1)\right)},
\end{eqnarray}
where we have used the fact that
for $L_0\ge 2\log c_4$ it is true that
  $\frac{\log N_0L_0}{L_0} \leq 1+\log \left(\frac{1}{\xi}\right)$ for
all  $u_0 \in [\xi^{L_0/d},1]$.
 Meanwhile, for the second term of the
rightmost side of (\ref{similar}), we have that

\begin{eqnarray*}
&  \P\left(f(1-n_0)>e^{\frac{(c_1'-1) c_1L_1\log \left(\frac{1}{\xi}\right)}{2}}\right)=\P\left(\sum_{k=1-n_0}^{1+n_0}\prod_{j=k+1}^{n_0+1}\widehat{\rho}(j,\omega)^{-1}>e^{\frac{(c_1'-1)c_1L_1 \log \left(\frac{1}{\xi}\right)}{2}}\right)\\
&  \leq \sum_{k=1-n_0}^{1+n_0}\P\left(\prod_{j=k+1}^{n_0+1}
\widehat{\rho}(j,\omega)^{-\alpha/2}> \frac{e^{(c_1'-1) c_1L_1 \alpha \log
\left(\frac{1}{\xi}\right)/4}}{(2n_0+1)^{\alpha/2}} \right)
\\
&  = e^{\frac{-(c_1'-1)c_1L_1 \alpha \log \left(\frac{1}{\xi}\right)}{4}}(2n_0+1)^{\alpha/2}\sum_{k=1-n_0}^{1+n_0}\prod_{j=k+1 }^{n_0+1}\E\left[\widehat{\rho}(j,\omega)^{-\alpha}\right]^{1/2}
\end{eqnarray*}
In analogy to  (\ref{estimate 2d aux}), we can conclude that
$\E\left[\widehat{\rho}(j,\omega)^{-\alpha}\right]^{1/2}
\le 2L_1^{3(d-1)}e^{\frac{(\log \eta_\alpha)m}{2}}$.
Therefore,
 for $L_0\ge 2\log c_4$ we see that
\begin{eqnarray}
\nonumber
&\P\left(f(1-n_0)>e^{\frac{(c_1'-1) c_1L_1 \log \left(\frac{1}{\xi}\right)}{2}}\right)\\
\nonumber
&\leq e^{\frac{-(c_1'-1)c_1L_1 \alpha \log \left(\frac{1}{\xi}\right)}{4}} (2n_0+1)^{\alpha/2}\sum_{k=1-n_0}^{1+n_0}2^{n_0+1-k}L_1^{3(d-1)(n_0+1-k)}e^{(\log \eta_{\alpha})(2c_1L_0+1)(n_0+1-k)}\\
\nonumber
&\leq e^{\frac{-(c_1'-1) \alpha \log \left(\frac{1}{\xi}\right)c_1L_1}{4}+\left(\frac{\alpha+2}{2}\right)\log(2n_0+1)+2n_0 \log 2+6(d-1)n_0\left(\frac{\log L_0N_0}{L_0}\right)L_0+(4c_1L_0n_0+2n_0)\log \eta_{\alpha}}\\
\label{estimate 2e}
&\leq e^{-L_1 \left( \frac{(c_1'-1) \alpha \log \left(\frac{1}{\xi}\right)c_1}{4}-1-6(d-1)\left(1+\log \left(\frac{1}{\xi}\right) \right)-(4c_1+1)(\log \eta_{\alpha}) \right)}.
\end{eqnarray}
 Now, in view of (\ref{estimate 2c}), (\ref{similar}), (\ref{estimate 2d}) and (\ref{estimate 2e})
the first term on the right side of (\ref{estimate 2a}) is bounded by
\begin{equation}
\label{estimate 2f}
2e^{-L_1 \left( \frac{(c_1'-1) \alpha \log \left(\frac{1}{\xi}\right)c_1}{4}-2 \alpha c_1 \log \left(\frac{1}{\xi}\right)-1-6(d-1)\left(1+\log \left(\frac{1}{\xi}\right) \right)-(4c_1+1)(\log \eta_{\alpha}) \right)}.
\end{equation}
Now, since $c_1'\ge 13+\frac{24d}{\alpha}+\frac{24d+12\log\eta_\alpha}{\alpha \log \frac{1}{\xi}}$, we conclude that

\begin{eqnarray*}
&\frac{(c_1'-1) \alpha \log \left(\frac{1}{\xi}\right)c_1}{4}-2 \alpha c_1
 \log \left(\frac{1}{\xi}\right)-1-6(d-1)\left(1+\log
 \left(\frac{1}{\xi}\right) \right)-(4c_1+1)(\log \eta_{\alpha})\\
& \geq \alpha c_1 \log \left(\frac{1}{\xi}\right)-2c_1 \log \eta_{\alpha},
\end{eqnarray*}
and therefore, by (\ref{estimate 2b}) and (\ref{estimate 2f}) we have that

\begin{equation}
\label{second estimate}
E\left[\rho_1^{a/2},A_2\right] \leq 4
e^{-c_1L_1 \log \frac{1}{\xi^ \alpha \eta_{\alpha}^2}}.
\end{equation}

\smallskip

\noindent {\it Step 3: estimate of $\E[\rho^{a/2}_1,G, A_3]$}.
Here we will estimate the third term of the inequality
(\ref{dec-a2a3}). Specifically we will show that
\begin{equation}
\label{third estimate}
\E[\rho^{a/2}_1,G, A_3]\leq 2 e^{-c_1L_1 \log \frac{1}{\xi^ \alpha \eta_{\alpha}^2}}.
\end{equation}
This upper bound will be almost obtained as the previous case, where we
achieved (\ref{second estimate}).  Indeed, in analogy to the
development of (\ref{rho1}) in Step 3, one has that for $\omega \in G$,
$$
\rho_1 \leq
\frac{f(0)+f(1-n_0)\xi^{(c_1'-1)c_1L_1}}
{(f(1-n_0)-f(0)-\xi^{(c_1'-1)c_1L_1}f(1-n_0))_{+}}.
$$
 But, if $\omega \in A_3$ also, one easily gets that $0<\rho_1 \leq 1$ if $L_0\geq \frac{\alpha \log 4}{2 \log \eta_{\alpha}}$. Thus,

\begin{equation}
\nonumber
\E\left[\rho_1^{a/2}, G, A_3 \right] \leq \P \left(A_3\right).
\end{equation}
Therefore, since $c_1'\ge 13+\frac{24d}{\alpha}+\frac{24d+12\log\eta_\alpha}{\alpha \log \frac{1}{\xi}}$,
it is enough to prove that

\begin{equation}
\label{a3}
\P(A_3)\le
  2 e^{-L_1 \left( \frac{(c_1'-1) \alpha \log \left(\frac{1}{\xi}\right)c_1}{4}-1-6(d-1)\left(1+\log \left(\frac{1}{\xi}\right) \right)-(4c_1+1)(\log \eta_{\alpha}) \right)}.
\end{equation}
 To justify this inequality,  note that
\begin{eqnarray*}
&\P(A_3) \leq \P \left(f(0)<e^{\frac{-(c_1'-1)c_1L_1 \log \frac{1}{\xi}}{2}}\right)+\P \left(f(1-n_0)>e^{\frac{(c_1'-1)c_1L_1 \log \frac{1}{\xi}}{2}}\right)\\
&\le\P\left(\prod_{j=1}^{n_0+1}\widehat{\rho}(j,\omega)^{-1}<
e^{\frac{-(c_1'-1)c_1L_1\log\left(\frac{1}{\xi}\right)}{2}} \right)+\P \left(f(1-n_0)>e^{\frac{(c_1'-1)c_1L_1 \log \frac{1}{\xi}}{2}}\right)
,
\end{eqnarray*} and hence we are in a very similar situation as in
(\ref{similar}) and development in (\ref{f2n}) and (\ref{estimate 2e}),
from where we derive (\ref{a3}).

\smallskip

\noindent {\it Step 4: estimate of $\E[\rho_1^{a/2}, G^c]$}.
Here we will prove that there exist constants $c_{3,21}(d)$ and $c_{3,22}(d)$
such that
\begin{equation}
\label{fourth estimate}
\E[\rho_1^{a/2}, G^c]\le
c_{3,21}\xi^{-c_1'c_1L_1}\left(c_{3,22}\widetilde{L}_1^{(d-2)}\frac{L_1^3}{L_0^2}\widetilde{L}_0
  \E[q(0)]\right)^{\frac{\widetilde{L}_1}{12N_0\widetilde{L}_0}}+
e^{-c_1L_1\log\frac{1}{\xi^\alpha\eta_\alpha^2}}.
\end{equation}
Firstly, we need to consider the event

 $$
A_4:=\left\{\omega \in \Omega: P_{0,\omega}\left(T^l_{L_1+1} \leq
\widetilde{T} \wedge \widetilde{T}^{l}_{1-L_1}\right) \geq
\xi^{2c_1L_1}\right\}.
$$
In the case that $\omega \in G^c \cap A_4$, the walk behaves
as if effectively it satisfies a uniformly ellipticity condition with
constant $\kappa=\xi$, so that we can follow exactly the
same reasoning presented by Sznitman  in \cite{Sz02} leading to
inequality (2.32) of that paper,
showing that there exist constants $c_{3,21}(d)$, $c_{3,22}(d)$ such that
whenever   $\widetilde{L}_1 \geq 48 N_0 \widetilde{L}_0$
one has that
\begin{equation}
\label{fourth estimate a}
\E\left[\rho_1^{a/2}, G^c, A_4\right] \leq
\xi^{-c_1L_1}\P(G^c)\le
 c_{3,21}\xi^{-c_1'c_1L_1}\left(c_{3,22}\widetilde{L}_1^{(d-2)}\frac{L_1^3}{L_0^2}\widetilde{L}_0 \E[q(0)]\right)^{\frac{\widetilde{L}_1}{12N_0\widetilde{L}_0}}.
\end{equation}
The second inequality of (\ref{fourth estimate a}) does not
use any uniformly ellipticity assumption.
It would be enough now to prove that

\begin{equation}
\label{a4c}
\E(\rho_1^{a/2},A_4^c)\le e^{-c_1L_1\log\frac{1}{\xi^\alpha\eta_\alpha^2}}.
\end{equation}
To do this we will follow the reasoning presented in Step 2.
Namely, for $L_0\ge 1+\sqrt{d}$, there is a nearest neighbor self-avoiding
path $(x_1,\ldots,x_n)$ with $n$ steps from $0$ to $\partial_+B_1$
such that $2c_1L_1\le n\le 2c_1L_1+1$, $x_1,\ldots,x_n\in B_1$ and
$x_n\cdot l\ge L_1+1$. Therefore

$$
A_4^c \subset \left\{\omega \in \Omega: \prod_{1}^{n}\omega(x_i,\Delta x_i)<
\xi^{n}\right\}=\left\{\omega \in \Omega: \sum_{1}^{n}\log\frac{1}{
\omega(x_i,\Delta x_i)}>n\log\frac{1}{\xi}\right\},$$
 so that
\begin{eqnarray*}
&  \E\left[\rho_1^{a/2}, A_4^c\right] \leq \E \left[e^{\alpha/2
    \sum_{1}^{n}\log\frac{1}{\omega(x_i,\Delta
      x_i)}},\sum_{1}^{n}\log\frac{1}{\omega(x_i,\Delta x_i)}>n\log\frac{1}{\xi}\right] \nonumber \\
&  \leq \E \left[e^{\alpha \sum_{1}^{n}\log\frac{1}{\omega(x_i,\Delta
      x_i)}}\right]^{1/2} \P\left(\sum_{1}^{n}\log\frac{1}{\omega(x_i,\Delta
  x_i)}>n\log\frac{1}{\xi}\right)^{1/2} \nonumber \\
&  \leq e^{(2\log\eta_\alpha-\alpha\log\frac{1}{\xi}) c_1L_1},
\end{eqnarray*}
which proves (\ref{a4c}) and finishes Step 4.

\smallskip

\noindent {\it Step 5: conclusion}. Combining the estimates (\ref{first estimate}) of step 1,
(\ref{second estimate}) of step 2, (\ref{third estimate}) of step 3 and (\ref{fourth estimate}) of step 4, we have (\ref{estimate 1}).
\end{proof}

\medskip

\noindent We will now prove a corollary of Proposition \ref{prop 1}, which will
imply Proposition \ref{pec}. For this, it will be important to note that
the statement of Proposition \ref{prop 1} is still valid if given $k\ge 1$ we change
$L_0$ by $L_k$, $L_1$ by $L_{k+1}$, $\tilde L_0$ by $\tilde L_k$
and $\tilde L_1$ by $\tilde L_{k+1}$. In effect, to see this, it is
enough to note that inequality (\ref{induction-in}) is satisfied with these
replacements. 
Define

$$
c_{3,23}:=e^{-\frac{4c_1\log \eta_{\alpha}}{\left(c_1-1\right)\alpha}}.
$$

%\begin{equation}
%\label{induction condition}
%c_9c_{10}\widetilde{L}_{k+2}^{d-1}L_{k+1}e^{-c_1L_{k+1}(K-2\log \eta)} \leq \xi^{u_{k+1}L_{k+1}}
%\end{equation}

\begin{corollary}
\label{inductive cor}
Let $0 <\xi <\min \{c_{3,23},  e^{-1/24} \} $ and
$\alpha >0$. Let $\{L_k:k\ge 0\}$ and  $\{\widetilde{L}_k:k\ge 0\}$
 be sequences satisfying (\ref{specifications 0}), (\ref{l0l1}) and (\ref{specifications 1}). Then
there exists $c_{3,25}(d,\alpha)>0$,  such that when for some $L_0 \geq c_{3,25}$, $a_0 \in (0,\alpha]$,
$u_0 \in [\xi^{L_0/d},1]$, it is true that
\begin{equation}
\label{induction first}
\phi_0:=c_{3,19} \widetilde{L}_1^{d-1}L_0 \E[\rho_0^{a_0}] \leq \xi^{\alpha u_0 L_0},
\end{equation}
then for all $k \geq 0$,
\begin{equation}
\label{induction term}
\phi_k:=c_{3,19} \widetilde{L}_{k+1}^{d-1}L_k \E[\rho_k^{a_k}] \leq (k+1) \xi^{\alpha u_k L_k},
\end{equation}
with $a_k:=a_0 2^{-k}$, $u_k:=u_0 8^{-k}$.

\end{corollary}

\begin{proof} We will use induction in $k$ to prove (\ref{induction term}). By hypothesis we only need to show $(\ref{induction term})$ for $n=k+1$ assuming that (\ref{induction term}) holds for $n=k$. To do this, with the help of Proposition \ref{prop 1} we have that for any $k \geq 0$

\begin{eqnarray*}
&\E[\rho_{k+1}^{a_{k+1}}]\leq c_{3,18} \left\{ \xi^{-c_2'L_{k+1}}\left(c_{3,19} \widetilde{L}_{k+1}^{(d-2)}\frac{L_{k+1}^3}{L_{k}^2}\widetilde{L}_k\E[q_k] \right)^{\frac{\widetilde{L}_{k+1}}{12N_k \widetilde{L}_k}}\right.\\
&+ \left.\sum_{0 \leq m \leq N_k+1} \left(c_{3,19}
  \widetilde{L}_{k+1}^{(d-1)}\E[\rho_k^{a_k}]\right)^{\frac{[N_k]+m-1}{2}}
  +e^{-c_{1} L_{k+1}
 \log \frac{1}{\xi^\alpha\eta_\alpha^2}}\right\},
\end{eqnarray*}

\noindent so that, for $k \geq 0$ and with the help of
 (\ref{specifications 1})
\begin{eqnarray*}
&\phi_{k+1} \leq c_{3,18}c_{3,19} \widetilde{L}_{k+2}^{(d-1)}L_{k+1}
\left\{\xi^{-c_2'L_{k+1}}\phi_k^{N_k^2/12}+\sum_{0 \leq m \leq N_k+1}
\phi_k^{\frac{[N_k]+m-1}{2}}\right\} \nonumber \\
&+ c_{3,18}c_{3,19} \widetilde{L}_{k+2}^{(d-1)}L_{k+1}e^{-c_{1}L_{k+1} \log \frac{1}{\xi^\alpha\eta_\alpha^2}}.
\end{eqnarray*}

\noindent Since  $\xi < c_{3,23} $, we can assert that $c_{3,18}c_{3,19} \widetilde{L}_{k+2}^{(d-1)}L_{k+1} e^{-c_{1} L_{k+1}  \log \frac{1}{\xi^\alpha\eta_\alpha^2}} \leq \xi^{ \alpha u_{k+1}L_{k+1}}$. Hence, we only need to prove that

\begin{equation}
\label{1st estimate}
c_{3,18}c_{3,19} \widetilde{L}_{k+2}^{(d-1)}L_{k+1}\left\{\xi^{-c_2'L_{k+1}}\phi_k^{N_k^2/12}+\sum_{0 \leq m \leq N_k+1}\phi_k^{\frac{[N_k]+m-1}{2}}\right\} \leq (k+1) \xi^{ \alpha u_{k+1}L_{k+1}}
\end{equation}

\noindent Firstly, note that for $L_0$ large enough
 by the induction hypothesis, (\ref{specifications 1}) and the fact
that $\xi < e^{-\frac{1}{24}}$

\begin{eqnarray*}
&\xi^{-c_2'L_{k+1}}\phi_k^{N_k^2/24} \leq  \xi^{-c_2'L_{k+1}} (k+1)^{\frac{N_k^2}{24}} \xi^{\frac{\alpha u_k N_k^2L_k}{24}}\\
& \leq  e^{c_2' \left(\log \frac{1}{\xi}+\frac{1}{24}-\frac{c_4 \log \frac{1}{\xi}}{24}\right) L_{k+1}} \leq 1.
\end{eqnarray*}
Substituting this estimate back into (\ref{1st estimate}) and using the hypothesis induction again, we obtain that

\begin{eqnarray*}
& c_{3,18} c_{3,19} \widetilde{L}_{k+2}^{(d-1)}L_{k+1}\left\{\xi^{-c_2'L_{k+1}}\phi_k^{N_k^2/12}+\sum_{0 \leq m \leq N_k+1}\phi_k^{\frac{[N_k]+m-1}{2}}\right\}\\
& \leq c_{3,18}c_{3,19}\widetilde{L}_{k+2}^{(d-1)}L_{k+1}\left\{\phi_k^{N_k^2/24}+(N_k+2)\phi_k^{N_k/4}\right\}\\
&  \leq  c_{3,24}\widetilde{L}_{k+2}^{(d-1)}L_{k+1}N_{k+1}\phi_{k}^{N_k/4}\\
&  \leq  c_{3,24}  \widetilde{L}_{k+2}^{(d-1)}L_{k+2}\phi_{k}^{N_k/8}(k+1)^{N_k/8} \xi^{\alpha u_k L_k N_k/8}\\
&  =  c_{3,24}  \widetilde{L}_{k+2}^{(d-1)}L_{k+2}\phi_{k}^{N_k/8}(k+1)^{N_k/8-1} (k+1)\xi^{\alpha u_{k+1} L_{k+1}},\\
\end{eqnarray*}
where $c_{3,24}:=2c_{3,18}c_{3,19}$.
Thus, in order to show that $\phi_{k+1} \leq (k+2)\xi^{\alpha u_{k+1}L_{k+1}}$ it is enough to prove that
\begin{equation}
\label{2nd estimate}
c_{3,24}\widetilde{L}_{k+2}^{(d-1)}L_{k+2}(k+1)^{N_k/8-1}\phi_k^{N_k/8} \leq 1.
\end{equation}

\noindent First, note that by the induction hypothesis,
\begin{equation}
\label{firste}
c_{3,24}\widetilde{L}_{k+2}^{(d-1)}L_{k+2}(k+1)^{N_k/8-1}\phi_k^{N_k/8} \leq
c_{3,24}\widetilde{L}_{k+2}^{(d-1)}L_{k+2}(k+1)^{N_k/4-1}\xi^{6 c_2' L_k}.
\end{equation}
 From (\ref{l0l1}), (\ref{specifications 1}) and (\ref{specifications 2}), we can say that

\begin{eqnarray}
&  c_{3,24} \widetilde{L}_{k+2}^{(d-1)}L_{k+2}(k+1)^{N_k/4-1}\xi^{6 c_2' L_k} \nonumber \\
& =c_{3,24} \left(\frac{L_{k+2}}{L_0}\right)^{3(d-1)}\widetilde{L}_0^{(d-1)}L_{k+2}(k+1)^{N_k/4-1}\xi^{6 c_2' L_k} \nonumber \\
& \leq c_{3,24} L_{k+2}^{3(d-1)+1}(k+1)^{N_k/4-1}\xi^{6 c_2' L_k}\nonumber \\
&  = c_{3,24} (N_{k+1}N_k)^{3d-2}L_k^{3d-2}(k+1)^{N_k/4-1}\xi^{ 6c_2' L_k} \nonumber \\
\label{3rd estimate}
&  \leq c_{3,24} 8^{3d-2}L_k^{3d-2}(k+1)^{N_k/4-1}\xi^{c_1 L_k}
N_k^{6d}\xi^{(6 c_1'-1)c_1L_k}.
\end{eqnarray}
But, note that

$$c_{3,24}8^{3d-2}L_k^{3d-2}(k+1)^{N_k/4-1}\xi^{c_1 L_k}
\leq 1.
$$
for $L_0$ large enough. Hence, substituting this estimate back into (\ref{3rd
  estimate}) and (\ref{firste}) we deduce that

\begin{equation}
\nonumber
c_{3,24}\widetilde{L}_{k+2}^{(d-1)}L_{k+2}(k+1)^{N_k/8-1}\phi_k^{N_k/8} \leq N_k^{6d}\xi^{(6 c_1'-1)c_1L_k} \leq N_k^{6d}\xi^{77c_1L_k},
\end{equation}
by our choice of $c_1'$. Finally, choosing $L_0$ large enough, the expression $N_k^{6d}\xi^{77c_1L_k} \leq 1$ for all $k
\geq 1$. In the case of $k=0$, we have that

$$
\left(\frac{c_4}{u_0}\right)^{6d}\xi^{77c_1L_0} \leq u_0^{-6d} \xi^{ 6L_0} \leq 1
$$ by our assumption on $u_0$. Then (\ref{2nd estimate}) follows and thus we get $(\ref{induction term})$ by induction and
choosing $L_0 \geq c_{3,25}$ for some constant $c_{3,25}>0$.

\end{proof}

\medskip

The following corollary
implies Proposition \ref{pec}. Since such a derivation
follows exactly the argument presented by Sznitman in \cite{Sz02}, we omit it.

\begin{corollary}

\label{cec} Let $l\in\mathbb S^{d-1}$, $d\ge 2$ and
%Assume that $\xi<\min\{c_{11},\eta_\alpha^{-2/\alpha}, %e^{-1/24}\}$.
 $\Upsilon=\max \left\{\frac{\alpha}{24}, \left(\frac{2c_1}{c_1-1}\right)\log \eta_{\alpha}^2 \right\}$.
Then, there exist constants $c_{3,26}=c_{3,26}(d)>0$ and $c_{3,27}=c_{3,27}(d)>0$ such that
if the following inequality is satisfied
\begin{equation}
\label{ec1}
c_{3,26}(d)
\inf_{L_0\ge c_{3,27}, 3\sqrt{d}\le\tilde L_0<L^3_0}\inf_{0<a\le \alpha}\left\{
\Upsilon^{3(d-1)}
 \tilde L_0^{d-1}L_0^{3(d-1)+1}\E[\rho^a_B]
\right\}<1,
\end{equation}
where $B=B(R,L_0-1,L_0+1,\tilde L_0)$,
then there exists a constant $c_{3,28}>0$ such that
\begin{equation}
\nonumber
\varlimsup_{L \to \infty}L^{-1}e^{c_{28}(\log L)^{1/2}}\log P_0 \left[\tilde{T}^l_{-\tilde{b}L}<T^l_{bL}\right]<0, \,\,\, \mbox{for}\,\, \mbox{all}\,\,\, b, \tilde{b}>0.
\end{equation}
\end{corollary}

\medskip

\noindent \begin{proof} If (\ref{ec1}) holds then there is a $\xi>0$ such that

\begin{equation}
\label{ec1a}
c_{3,26}(d)
\inf_{L_0\ge c_{3,27}, 3\sqrt{d}\le\tilde L_0<L^3_0}\inf_{0<a\le \alpha}\left\{
\left(\alpha \log \frac{1}{\xi}\right)^{3(d-1)}
 \tilde L_0^{d-1}L_0^{3(d-1)+1}\E[\rho^a_B]
\right\}<1,
\end{equation}
with $\xi <\{c_{3,23}, e^{-1/24}\}$. Then,
by (\ref{specifications 0}) and (\ref{l0l1}),
$$
\widetilde{L}_1^{d-1}L_0=\left(\frac{c_4}{u_0}\right)^{3(d-1)}\widetilde{L}_0^{d-1}L_0.
$$
Now, the maximum of  $u_0^{3(d-1)}\xi^{\alpha u_0L_0}$, as a function of $u_0$ for
$u_0 \in [\xi^{\frac{L_0}{d}},1]$,
 is given by $c_{3,29}(d)\left(\alpha L_0 \log \frac{1}{\xi}\right)^{-3(d-1)}$
for $u_0=\frac{3(d-1)}{\alpha L_0 \log \frac{1}{\xi}}$,
when $L_0$ is large enough,
 where $c_{3,29}(d):=\left(\frac{3(d-1)}{e}\right)^{3(d-1)}$. Thus if (\ref{ec1a}) holds, (\ref{induction first}) holds as well. Hence, applying Corollary \ref{inductive cor} we can say that (\ref{induction term}) is true for all $k \geq 0$.
The same reasoning used by Sznitman in \cite{Sz02} to derive
Proposition 2.3 of that paper gives the estimate

\begin{eqnarray*}
& P_0 \left(\tilde{T}^l_{-\tilde{b}L}<T^l_{bL}\right) \leq \left(|C|+\frac{bL}{L_k}+1\right) \xi^{\frac{\alpha}{2}u_kL_k} \nonumber \\
& \leq e^{-\widetilde{b}L e^{-c_{3,28}(\log\, \widetilde{b}L)^{\frac{1}{2}}}},
\end{eqnarray*}
for some constant $c_{3,28}>0$ and $L$ large enough, and where we have chosen $u_0=\frac{3(d-1)}{\alpha L_0 \log \frac{1}{\xi}}$,

\end{proof}

\section{An atypical quenched exit estimate}
\label{section-four}
\setcounter{equation}{0}
Here we will prove a crucial atypical quenched exit estimate
for tilted boxes, which will subsequently enable us in section
\ref{tail-regen} to show that the regeneration times of the random
walk are integrable. Let us first introduce some basic notation.

Without loss of generality, we will assume that $e_1$ is contained
in the open half-space defined by the asymptotic direction so that

\begin{equation}
\label{secondwlg}
\hat v\cdot e_1>0.
\end{equation}
Recall the definition of the hyperplane perpendicular to direction $e_1$
in (\ref{hiperplane}) so that
$$
H:=\{x\in\mathbb R^d: x\cdot e_1=0\}.
$$
Let $P:=P_{\hat v}$ (see (\ref{projl}))
be the projection on the asymptotic
direction along the hyperplane $H$ defined for $z\in\mathbb Z^d$

$$
Pz:=\left(\frac{z\cdot e_1}{\hat{v}\cdot e_1}\right)
\,\hat{v},
$$
and $Q:=Q_l$ (see (\ref{projql})) be
the projection of $z$ on $H$ along $\hat v$ so that

$$
Qz:=z-Pz.
$$
Now, for $x\in\mathbb Z^d$, $\beta>0$, $\varrho>0$ and $L>0$, define the {\it tilted boxes with
respect to the asymptotic direction $\hat v$} as

\begin{equation}
\label{tilted-box}
B_{\beta,L}(x):=\left\{y \in\mathbb Z^d: - L^{\beta} < (y-x) \cdot e_1 < L;
\,\, \|Q(y-x)\|_{\infty} <\varrho L^{\beta}\right\}.
\end{equation}
 and their {\it front boundary} by

\begin{equation}
\nonumber
\partial_{+} B_{\beta,L}(x):=\{y \in \partial B_{\beta,L}(x): y\cdot e_1-x
\cdot e_1\ge L\}.
\end{equation}

\begin{figure}[H]
\centering
\includegraphics[width=8cm]{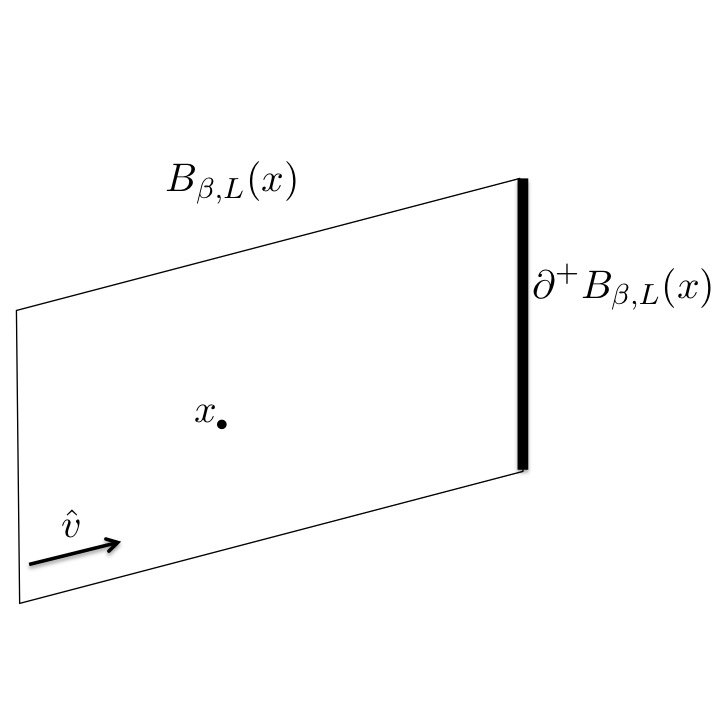}
\caption{The box $B_{\beta, L}(x)$.}
\label{box4.1}
\end{figure}

\noindent See Figure \ref{box4.1} for a picture of the box $B_{\beta,L}$ and
its front boundary.

\medskip

\begin{proposition}
\label{crucial-lemma}
Let $\alpha>0$ and assume that
$\eta_\alpha<\infty$ as defined in (\ref{eta-alpha}).
 Let $l\in\mathbb S^{d-1}$ and $M\ge 15d+5$. Assume that
$(P)_M|l$ is satisfied. Let
$\beta_0 \in (1/2,1)$, $\beta\in\left(\frac{\beta_0+1}{2},1\right)$ and
$\zeta \in (0,\beta_0)$.
 Then, for each
$\kappa>0$ we have that

$$
\limsup_{L\to\infty}
L^{-g(\beta_0,\beta,\zeta)}
\log \mathbb P\left(P_{0,\omega}\left(X_{T_{B_{\beta,L}(0)}} \in
 \partial_{+}B_{\beta,L}(0)\right)\leq e^{-\kappa L^\beta} \right)<0,
$$

\noindent where

\begin{equation}
\label{def-g}
g(\beta_0,\beta,\zeta):=\min\{\beta+\zeta,
3\beta-2+(d-1)(\beta-\beta_0)\}.
\end{equation}

\end{proposition}

\medskip

\noindent We will now prove Proposition \ref{crucial-lemma} following
  ideas similar in spirit from those presented by Sznitman in \cite{Sz02}.

\subsection{Preliminaries}

Firstly we need to define an appropriate mesoscopic scale
to perform a renormalization analysis.
Let $\beta_0\in (0.5,1)$,
$\beta\in (\beta_0,1)$
and $\chi:=\beta_0+1-\beta\in(\beta_0,1]$.
Define

$$
L_0:=\frac{L-\varrho L^{\beta_0}}{[ L^{1-\chi}]}.
$$
Now, for each $x \in \R^d$  we consider the
{\it mesoscopic box}

$$
\tilde{B}(x):=\left\{y \in\mathbb Z^d: - L^{\beta_0} < (y-x) \cdot e_1 <\varrho
L_0; \,\, \|y-x-P(y-x)\|_{\infty} < (1+\varrho)L^{\beta_0}\right\},
$$
and its {\it central part}

$$
\tilde{C}(x):=\left\{y \in\mathbb Z^d: 0 \leq (y-x) \cdot e_1 < \varrho L_0; \,\,
\|y-x-P(y-x)\|_{\infty} < L^{\beta_0}\right\}.
$$
Define also
$$
\partial_{+} \tilde B(x):=\{y \in \partial\tilde
 B(x): y\cdot e_1-x \cdot e_1 \ge\varrho L_0\}
$$
and

$$
\partial_{+} \tilde C(x):=\{y \in \partial\tilde
 C(x): y\cdot e_1-x \cdot e_1\ge \varrho L_0\}.
$$
\noindent We now say that a box $\tilde{B}(x)$ is  {\it good}
  if

$$
\sup_{x\in\tilde C(x)}P_{x,\omega}\left(X_{T_{\tilde B(x)}}\not \in\partial_+\tilde
B(x)\right)<\frac{1}{2},
$$
Otherwise the box is called {\it bad}. At this point,
by Theorem \ref{theorem1} proved in section \ref{pol-t}, we have
 the following version of Theorem \ref{sznitman01} (Theorem A.2
 of Sznitman  \cite{Sz02}).

\begin{theorem}
\label{sznitman02}
Let $l\in\mathbb S^{d-1}$ and $M\ge 15d+5$.
Consider  RWRE satisfying condition $(P)_M|l$
and the ellipticity condition $(E)_0$.
Then, for any $c>0$
and $\rho\in(0.5,1)$,

$$
\limsup_{u\to\infty}u^{-(2\rho-1)}\log P_0
\left(\sup_{0\le n\le T_u^{e_1}}|X_n-P(X_n)|\ge cu^\rho\right)<0,
$$
where $T_u^{e_1}$ is defined in (\ref{tele}).
\end{theorem}

\medskip

\noindent The following lemma  is an important corollary of
Theorem \ref{sznitman02}.

\medskip
\begin{lemma}
\label{bad box lemma}
Let $l\in\mathbb S^{d-1}$ and $M\ge 15d+5$.
Consider  RWRE satisfying the ellipticity condition $(E)_0$ and
condition $(P)_M|l$.
 Then

\begin{equation}
\label{bad box estimate}
\limsup_{L\to\infty}L^{-(\beta+\beta_0-1)}
\log\mathbb P(\tilde{B}(0)\ {\rm is}\ {\rm bad})<0.
\end{equation}
\end{lemma}

\begin{proof}
By Chebyshev's inequality we have that

\begin{eqnarray}
\nonumber
&\mathbb P(\tilde{B}(0)\ {\rm is}\ {\rm bad})\\
\nonumber
&\le
2^{d-1}L_0L^{\beta_0(d-1)}\left(P_0\left(\sup_{0\le n\le T_{\varrho L_0}^{\hat v}}
|X_n-PX_n|\ge (1+\varrho)L^{\beta_0}\right)
+P_0\left(\tilde T^{\hat v}_{- L^{\beta_0}}<\infty\right)\right).
\end{eqnarray}
By Theorem \ref{sznitman02}, the first summand can be estimated
as
$$
\limsup_{L\to\infty} L^{-(\beta+\beta_0-1)}
\log P_0\left( \sup_{0\le n\le T^{\hat v}_{\varrho L_0}}|X_n-PX_n|\ge (1+\varrho) L^{\beta_0}
\right)<0.
$$
To estimate the second summand, since $(P)_M|l$ is
satisfied, by Theorem \ref{theorem1} and the equivalence given
by Theorem \ref{sznitman-equivalence}, we
can chose $\gamma$ close enough to $1$ so that $\gamma\beta_0
\ge \beta_0+\beta-1$ and such that

$$
\limsup_{L\to\infty}L^{-\gamma\beta_0}\log P_0
\left(\tilde T^{\hat v}_{- L^{\beta_0}}<\infty\right)<0.
$$ 
\end{proof}

% For $n_1, n_2 \in \N$ we define the sublattices

%\begin{equation}
%\label{sublattice}
%{\mathcal{L}}_{n_1,n_2}:= \left \{x \in \Z^d: x \cdot e_1
%\in n_1 \Z, x \cdot e_j \in n_2\Z, \, \forall \, j \in \{2, \ldots, d\}
%\right \}
%\end{equation}
%In this paper, we will take $n_1:=[L_0]$ and $n_2:=[2K L^{\beta_0}]$.

\noindent Let $k_1,\ldots, k_d\in\mathbb Z$.
 From now on, we will use the notation $x=(k_1,\ldots,k_d)\in\mathbb R^d$ to denote the
point

$$
x=k_1\frac{\varrho}{\hat v\cdot e_1}L_0 \hat v
+\sum_{j=2}^d  2k_j (1+\varrho)L^{\beta_0}e_j.
$$
Define the following set of points which will correspond
to the centers of mesoscopic boxes.

$$
\mathcal L:=\left\{x\in\R^d: x=(k_1,\ldots,k_d)\ {\rm for}\ {\rm some}\
k_1,\ldots,k_d\in\Z\right\}.
$$
We will use subsequently the following property of the lattice
$\mathcal L$: there exist $2^d$ disjoint sub-lattices
$\mathcal L_1,\ldots, \mathcal L_{2^d}$ such that
$\mathcal L=\cup_{i=1}^{2^d} \mathcal L_i$ and for each $1\le i\le 2^d$,
the sub-lattice $\mathcal L_i$ corresponds to the centers of
mesoscopic boxes which are pairwise disjoint.
 Let $\mathcal{L}_0$ be the set defined by
$$
\mathcal{L}_0:=\{x=(k_1, \ldots, k_d) \in \mathcal{L}: k_1=0\}.
$$
For each $x \in \mathcal{L}_0$ we define the {\it column} of mesoscopic boxes
 as

\begin{equation}
\nonumber
C_x:=\bigcup_{k_1=-1}^{[L^{1-\chi}]}\tilde{B}\left(x+k_1 \frac{\varrho}{\hat{v}\cdot e_1} L_0  \hat{v}\right)
\end{equation}
See Figure \ref{column4.1} for a picture of the column $C_x$, for some $x \in \mathcal L_0$.

\begin{figure}[H]
\centering
\includegraphics[width=10cm]{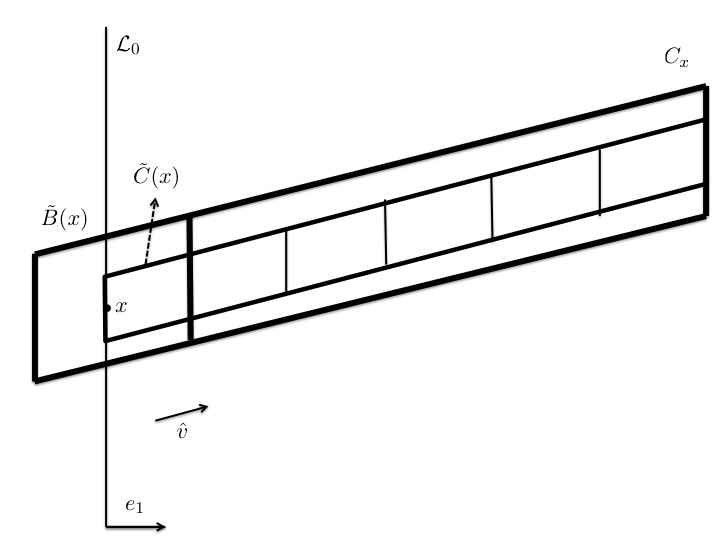}
\caption{A mesoscopic box $\tilde B$ with its corresponding middle part
  $\tilde C$, which belongs to the column $C_x$. }
\label{column4.1}
\end{figure}

\noindent The collection of these columns will be denoted by $\mathcal{C}$. Define now for each $C_x \in \mathcal{C}$
and $-1\le k\le [L^{1-\chi}]$ define

\begin{equation}
\nonumber
\partial_{k,1} C_{x}:=\partial_{+}\tilde{C}\left(x+k \frac{\varrho}{\hat{v}\cdot
  e_1} L_0  \hat{v}\right) \quad {\rm and}\quad \partial_{k,2} C_{x}
:=\partial_{+}\tilde{B}\left(x+k \frac{\varrho}{\hat{v}\cdot e_1} L_0
\hat{v}\right) \,\backslash
\partial_{k,1} {C}_{x}.
\end{equation}

%{\red Quiz\'as es mejor especificar que el primero camino est\'a dentro de $\tilde{B}$. Lo importante es que no se salga de la columna $C_x$}.
\noindent For each point $y \in \partial_{k,1}{C}_{x}$ we assign a
 path $\pi^{(k)}
=\{\pi_1^{(k)},\ldots,\pi_{n_1}^{(k)}\}$ in $\tilde B(x)$ with
$n_1:=\left[2c_1 \frac{\varrho}{\hat{v}\cdot e_1} L_0\right]$
steps from $y$ to $\partial_{k+1,1}{C}_{x}$,
so that $\pi_1^{(k)}=y$ and $\pi_{n_1}^{(k)}\in\partial_{k+1,1}{C}_{x}$.
For each point $z \in \partial_{k,2}{C}_{x}$ we assign a
 path $\bar \pi^{(k)}
=\{\bar\pi_1^{(k)},\ldots,\bar\pi_{n_2}^{(k)}\}$ with
$n_2:=\left[2c_1 \varrho L^{\beta_0}\right]$
steps from $z$ to $\partial_{k,1}{C}_{x}$,
so that $\bar\pi_1^{(k)}=z$ and $\bar\pi_{n_2}^{(k)}\in\partial_{k,1}{C}_{x}$.
We will also use the notation $\{m_1, \ldots, m_N\}$ to denote
some subset  of $\{-1, \ldots, [L^{1-\chi}] \}$ with $N$ elements.

Let $x\in\mathcal L_0$ and $\xi>0$. A column of boxes $\mathcal{C}_x \in \mathcal{C}$ will
be called {\it elliptically good} if it satisfies
the following two conditions

\begin{equation}
\label{egd1}
\sup_{N\le \left[\frac{L^{\beta}}{L_0}\right]}
\sup_{\{m_1, \ldots, m_N\}}
\sum_{j=1}^N
\sup_{y_{m_j}\in\partial_{m_j,1}C_x}
    \sum_{i=1}^{n_1}
\log
 \frac{1}{\omega(\pi^{(m_j)}_i,
 \Delta \pi^{(m_j)}_i)}\leq 2c_1\frac{\varrho}{\hat{v} \cdot e_1} \log \left(\frac{1}{\xi}\right)L^{\beta}
\end{equation}
and

\begin{equation}
\label{egd2}
\sum_{k=-1}^{[L^{1-\chi}]}
\sup_{z_k\in\partial_{k,2}C_x}
    \sum_{i=1}^{n_2}
\log
 \frac{1}{\omega(\bar\pi_i^{(k)}
 \Delta \bar\pi_i^{(k)})}\leq 2c_1\varrho
\log \left(\frac{1}{\xi}\right)L^{\beta}
.
\end{equation}
If either (\ref{egd1}) or (\ref{egd2}) is not  satisfied,
we will say that the column $C_x$ is {\it elliptically bad}.

\begin{lemma}
\label{elliptically bad column estimate}
For any $x \in \mathcal{L}_0$, $\beta \in \left[\frac{\beta_0+1}{2}, 1\right)$ and $\xi>0$ such that $\log \frac{1}{\xi^{2\alpha}\eta_{\alpha}^3}>0$ we have that

\begin{equation}
\label{bbce}
\limsup_{L \to \infty} L^{-\beta}\log \P \left(C_x \,\, \mbox{is} \,\, \mbox{elliptically} \, \, \mbox{bad}
\right)<0
\end{equation}

\end{lemma}
\begin{proof}
 Let us first
note that $\frac{L^{\beta}}{L_0} \geq1$ by our condition on $\beta$.
Now, it is clear that

\begin{equation}
\label{egd3}
\P \left(C_x \,\, \mbox{is} \,\, \mbox{elliptically} \, \, \mbox{bad}
\right) \leq \P\left( (\ref{egd1}) \,\mbox{is} \, \mbox{not} \, \mbox{satisfied}\right)+\P\left( (\ref{egd2}) \,\mbox{is} \, \mbox{not} \, \mbox{satisfied}\right)
\end{equation}
Regarding the first term on the right of (\ref{egd3})  and since $2 \beta-\beta_0-1<\beta-\beta_0<\beta$ we have that

\begin{eqnarray}
&\P\left( (\ref{egd1}) \,\mbox{is} \, \mbox{not} \, \mbox{satisfied}\right) \leq \sum_{N=1}^{[\frac{L^{\beta}}{L_0}]} \P \left(\exists \, \{m_1, \ldots, m_N\} \,\, \mbox{and}\,\, y_{m_j} \in \partial_{m_j,1}C_x \, \mbox{such} \,\,\mbox{that} \right. \nonumber \\
&\left.\sum_{j=1}^
    {N}
    \sum_{i=1}^{n_1}\log
 \frac{1}{\omega \left(\pi_i^{(m_j)},
 \Delta
      \pi_i^{(m_j)}\right)}  > 2c_1\frac{\varrho}{\hat{v}\cdot e_1} \log \left(\frac{1}{\xi}\right) L^{\beta} \right) \nonumber \\
      \label{bbce 1}
      &\leq \frac{L^{\beta}}{L_0}
       L^{(\beta-\beta_0) \frac{L^{\beta}}{L_0}}
 e^{(\log L)\beta_0(d-1)L^{\beta-\beta_0}} e^{2(\log \eta_{\alpha})c_1 \frac{\varrho}{\hat{v}\cdot e_1}L^{\beta}} e^{-2c_1 \frac{\varrho}{\hat{v}\cdot e_1} \left(\alpha \log \frac{1}{\xi}\right) L^{\beta}}
\leq e^{-c_{4,1}L^{\beta}}
\end{eqnarray}

\noindent for some constant $c_{4,1}>0$ if $L$ is large enough and $\log \frac{1}{\xi^{2\alpha}\eta_{\alpha}^3}>0$.

\medskip

\noindent Similarly for the rightmost term of (\ref{egd3}) we have that,

\begin{eqnarray}
 \nonumber
&\P\left( (\ref{egd2}) \,\mbox{is} \, \mbox{not} \, \mbox{satisfied}\right)\\
 \nonumber
&\leq \P \left(\exists \,  z_{k} \in \partial_{k,2}C_x \, \mbox{such} \,\,\mbox{that}
\sum_{k=-1}^{[L^{1-\chi}]}
    \sum_{i=1}^{n_2}
\log
 \frac{1}{\omega(\bar\pi_i^{(k)}
 \Delta \bar\pi_i^{(k)})}> 2c_1\varrho
\log \left(\frac{1}{\xi}\right)L^{\beta} \right) \nonumber \\
\label{bbce 2}
      & \leq e^{\log L(\varrho \beta_0(d-1)L^{1-\chi})} e^{2(\log \eta_{\alpha})c_1 \varrho L^{\beta}} e^{-2c_1 \varrho \left(\alpha \log \frac{1}{\xi}\right)L^{\beta}}
\leq e^{-c_{4,2}L^{\beta}}
\end{eqnarray}

\noindent for some constant $c_{4,2}>0$ if $L$ is large enough and $\log \frac{1}{\xi^{2\alpha}\eta_{\alpha}^3}>0$. Substituting (\ref{bbce 1}) and (\ref{bbce 2}) back into (\ref{egd3}), (\ref{bbce}) follows.
\end{proof}

The proof Proposition \ref{crucial-lemma}  will be
reduced to the control of the probability of the three events:
the first one, corresponding to subsection \ref{bad-boxes},
  gives a control on the number of bad boxes; the second
one, corresponding to subsection \ref{bad-column}, gives a control on the
number
 of elliptically good columns; the third
one, corresponding to subsection \ref{sec-confin},
 gives a control on the probability that the random walk can find
an appropriate path which leads to an elliptically good column.

\subsection{Control on the number of bad boxes}
\label{bad-boxes}

\noindent  We will need to consider only the mesoscopic boxes
which intersect the box $B_{\beta,L}(0)$ and whose $k_1$ index
is larger than or equal to $-1$. We hence define the collection
of mesoscopic boxes

\begin{equation}
\nonumber
\mathcal{B}:=\left\{\tilde B(x): \tilde B(x) \cap B_{\beta,L}(0) \neq
\emptyset, x=
(k_1,\ldots,k_d), \,\,
k_1,\ldots,k_d\in\Z, \,\, k_1 \geq -1 \right\}
\end{equation}

\noindent  In addition, we call the number of bad
mesoscopic boxes in $\mathcal B$,

$$
N(L):=\left|\left\{\tilde B \in \mathcal{B}: \tilde B \,\, \mbox{is} \,\,\mbox{bad}
\right\}\right|,
$$
and  for each $1\le i\le 2^d$, call the number of bad
mesoscopic boxes in $\mathcal B$ with centers in the sub-lattice
$\mathcal L_i$ as

$$
N_i(L):=\left|\left\{\tilde B(x) \in \mathcal{B}: \tilde B(x) \,\, \mbox{is}
\,\,\mbox{bad}\ \mbox{and}\ x\in\mathcal L_i
\right\}\right|.
$$
 Define

\begin{equation}
\nonumber
G_1:=\left\{\omega \in \Omega: N(L) \leq \frac{\varrho^{d-1} L^{(d-1)(\beta-\beta_0)}L^{\beta}}{2(1+\varrho)^{d-1}L_0}\right\}.
\end{equation}

\medskip

\begin{lemma}
\label{lemma43}
Assume that $\beta>\frac{\beta_0+1}{2}$. Then, there is a constant $c_{4,3}>0$ such that
for every $L>1$ we have that

$$
\P(G_1^c)\le
e^{-c_{4,3}L^{3\beta-2+(d-1)(\beta-\beta_0)}}.
$$

\end{lemma}
\begin{proof}
Note that the number of columns intersecting the
box $B_{\beta,L}(0)$ is equal to

$$
\left[\frac{\varrho^{d-1}
    L^{(d-1)(\beta-\beta_0)}}{(1+\varrho)^{d-1}}\right].
$$
Hence,
whenever $\omega \in G_1$, necessarily there exist at least
  $\left\lceil\frac{\varrho^{d-1}
    L^{(d-1)(\beta-\beta_0)}}{2(1+\varrho)^{d-1}}\right\rceil$ columns each one
with at most
 $\left\lceil\frac{L^{\beta}}{L_0}\right\rceil$  bad boxes.
Let us take
 $m_1:=\left[\frac{\varrho^{d-1} L^{(d-1)(\beta-\beta_0)}L^{\beta}}{2(1+\varrho)^{d-1}L_0}
 \right]$ and
$m_2:=|\mathcal B|=\left[\frac{\varrho^{d-2}L^{d(\beta-\beta_0)}}{(1+\varrho)^{d-1}}+
\frac{\varrho^{d-1}L^{(d-1)(\beta-\beta_0)}}
 {(1+\varrho)^{d-1}} \right]$.
Now, using the fact that the mesoscopic boxes in each sub-lattice
$\mathcal L_i$, $1\le i\le 2^d$, are disjoint, and
the estimate  (\ref{bad box estimate})
of Lemma \ref{bad box lemma},  we have by independence
that there exists a constant $c_{4,3}>0$ such that for every $L\ge 1$,

\begin{eqnarray*}
&\P(G_1^c)=\P\left(N(L) \geq m_1  \right)
\le\sum_{i=1}^{2^d} \P\left(N_i(L) \geq \frac{m_1}{2^d}  \right)
\nonumber\\
&\leq \sum_{i=1}^{2^d}\sum_{n=m_1/2^d}^{m_2}
{ m_2  \choose n}\P \left( \tilde{B}(0) \,\, \mbox{is} \,\, \mbox{bad}\right)^n \leq 2^d\sum_{n=m_1/2^d}^{m_2} m_2^{n}
e^{-nL^{\beta+\beta_0-1}} \nonumber\\
&\leq  e^{-c_{4,3}L^{\beta+\beta_0-1+(d-1)(\beta-\beta_0)+2\beta-\beta_0-1}} \leq e^{-c_{4,3}L^{3\beta-2+(d-1)(\beta-\beta_0)}}.
\end{eqnarray*}
 Note that in the second to last inequality
we have used the fact that $2 \beta+\beta_0-2>0$
which is equivalent to the condition $\beta>\frac{2-\beta_0}{2}$.
Now, this last condition is implied by the requirement
$\beta>\frac{\beta_0+1}{2}$.

\end{proof}

\subsection{Control on the number of elliptically bad columns}
\label{bad-column}

 Let
$m_3:=\left[\frac{\varrho^{d-1}
    L^{(d-1)(\beta-\beta_0)}}{2(1+\varrho)^{d-1}}\right]$
and define the event that any sub-collection of the set of columns
of cardinality larger than or equal to $m_3$ has at least one
elliptically good column

\begin{equation}
\label{event G2}
G_2:=\left\{\omega \in \Omega: \forall \, \mathcal{D} \subset \mathcal{C}, \,
|\mathcal{D}|\ge m_3,\, \exists\, C_x \in \mathcal{D} \,\, \mbox{such} \,\, \mbox{that}\,\, C_x \,\,\mbox{is} \,\, \mbox{elliptically}\,\, \mbox{good} \right\}.
\end{equation}
Here we will prove the following lemma.

\medskip
\begin{lemma}
\label{lemma44} There is a constant $c_{4,4}>0$ such that for every $L\ge 1$,

\begin{equation}
\nonumber
\P(G_2^c) \leq e^{-c_{4,4}L^{\beta+(d-1)(\beta-\beta_0)}}.
\end{equation}

\end{lemma}

\begin{proof} Note that the total number of columns intersecting the
box $B_{\beta,L}$ is equal to

$$
m_4:=\left[\frac{\varrho^{d-1}
    L^{(d-1)(\beta-\beta_0)}}{(1+\varrho)^{d-1}}\right].
$$
Using the fact that the events $\{C_x \,\, \mbox{is}\,\,
\mbox{elliptically}\,\, \mbox{bad}\}$, $\{C_y \,\, \mbox{is}\,\,
\mbox{elliptically}\,\, \mbox{bad}\}$ are independent if $x \neq y$, since
these columns are disjoint,
we conclude that there is a constant $c_{4,4}>0$ such that
for all $L\ge 1$,

\begin{eqnarray*}
&\P(G_2^c)=\P\left(\exists\, \mathcal{D} \subset \mathcal{C}, \, |\mathcal
D|\ge m_3 \,\, \mbox{such} \,\,\mbox{that}\,\, \forall C_x \in \mathcal{D}, \,\, C_x \,\, \mbox{is} \,\, \mbox{elliptically} \,\, \mbox{bad} \right)\\
&\leq \sum_{n=m_3}^{m_4}m_4^n \P\left(C_x \,\, \mbox{is}\,\,
\mbox{elliptically} \,\, \mbox{bad} \right)^n
\le e^{-c_{4,4}L^{\beta+(d-1)(\beta-\beta_0)}},
\end{eqnarray*}
 where in the last inequality
we have used the estimate (\ref{bbce})
of Lemma \ref{elliptically bad column estimate}
which provides a bound for the probability of a column to be elliptically bad.

\end{proof}

\subsection{The constrainment event}
\label{sec-confin}
Here we will obtain an adequate estimate for the probability that
the random walk hits an elliptically good column.
We will need to introduce some notation, corresponding to the
 the box where the random walk will move before
hitting the elliptically good column and a certain class of hyperplanes of this region.
Let first $\zeta \in (0, \beta_0)$, a parameter which gives the order
of width of the box $\bar B_{\zeta, \beta, L}$  where the random walk will be able to find
a reasonable path to the elliptically good column, so that

$$\bar B_{\zeta, \beta,L}:=\{x \in \Z^d: -L^{\zeta} \leq x \cdot e_1 \leq
L^{\zeta}, \|x-Px\|_{\infty} <L^{\beta}\}.$$
Note that this box is contained in $B_{\beta,L}(0)$ and that it also contains the
starting point $0$ of
the random walk.
Define now for each $0\le z\le L^\zeta$, the truncated hyperplane

$$
H_z:=\left\{x \in \bar B_{\zeta, \beta, L}: x \cdot e_1=z \right\},
$$
and consider the two collections of  truncated hyperplanes
defined as

$$
\mathcal H^{+}= \left\{H_z:z \in \Z, 0\leq z \leq L^{\zeta} \right\}
 \quad {\rm and}\quad \mathcal H^{-}= \left\{H_z:z \in \Z,
-L^{\zeta}\leq z <0\right \}.
$$

\noindent Whenever there is no risk of confusion, we will
drop the subscript from $H_z$ writing $H$ instead.
Let $r:=[2 \varrho L^{\beta}]$.
 Now, for each $H \in \mathcal H^{+} \cup \mathcal H^{-}$ and
each $j$ such that $e_j\ne\pm e_1$, we will consider the
set of paths $\Pi_j$ with $r$ steps defined by $\pi=
\{\pi_{1}, \ldots, \pi_{r}\}\in\Pi_j$ if and only if

\begin{equation}
\nonumber
\pi \subset H \quad {\rm and}\quad
\pi_{ i+1}-\pi_{i}=e_j.
\end{equation}

\medskip

\noindent In other words, $\pi$ is contained in the  truncated hyperplane $H$ and
it has steps which move only in the direction $e_j$.
 We now say that an hyperplane $H\in\mathcal H^+\cap\mathcal H^-$
is {\it elliptically good} if for all paths $\pi\in\cup_{j\ne 1, d+1}\Pi_j$
 one has that

\begin{equation}
\label{hyp1a}
\sum_{i=1}^{r}\log \frac{1}{\omega\left(\pi_{i}, \Delta \pi_{i}\right)} \leq 2\varrho \log \left(\frac{1}{\xi}\right)L^{\beta}.
\end{equation}
Otherwise $H$ will be called {\it elliptically bad} (See Figure \ref{coev4.1a}).

\begin{figure}[H]
\centering
\includegraphics[width=10cm]{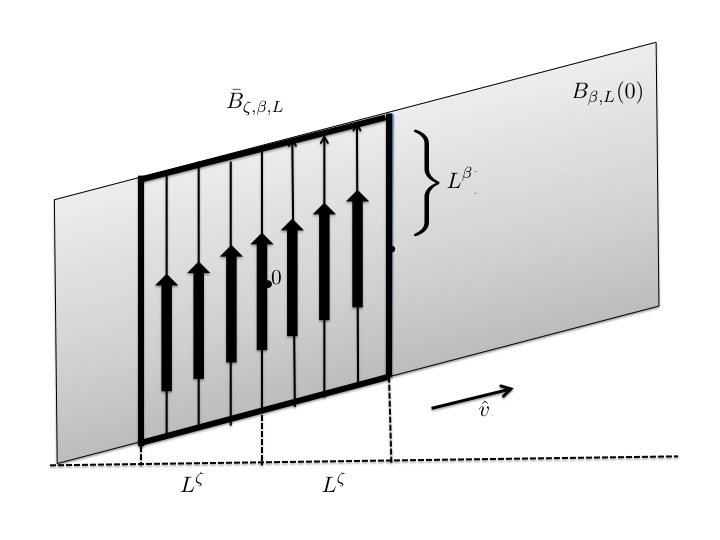}
\caption{The box $\bar B_{\zeta, \beta,L}$. The arrows indicate the
 ellipticity condition given by (\ref{hyp1a}), which implies that each hyperplane is elliptically good.}
\label{coev4.1a}
\end{figure}

\noindent From a routine counting argument and applying  Chebyshev inequality, note
that for each $H \in \mathcal H^{+} \cup \mathcal H^{-}$ and $\xi >0$ such that $\log \frac{1}{\xi^{\alpha}\eta_{\alpha}^2}>0$ there is a constant $c_{4,5}>0$ such that

\begin{equation}
\label{hyp1b}
\P\left( H \,\, \mbox{is} \,\, \mbox{elliptically} \, \mbox{bad}\right) \leq e^{-c_{4,5} L^{\beta}}.
\end{equation}

\noindent Now choose a rotation $\hat{R}$  such that
$\hat{R}(e_1)=\hat{v}$. Let $\hat{v}_j:=\hat{R}(e_j)$ for $j \geq 2$.
We now want to make a construction analogous to the
one which led to the concept of elliptically good
hyperplane. But now, we would need to define hyperplanes
perpendicular to the directions $\{\hat v_j\}$ which are not
necessarily equal to a canonical vector. Therefore, we will work here
with strips, instead of hyperplanes. For each $z\in\mathbb Z$ even and
$k \in \{2, \ldots, d\}$ consider the strip $I_{k,z}:=\{x \in \bar
B_{\zeta, \beta, L}(0): z-1<x \cdot \hat{v}_j<z+1\}$.
Consider also the two sets of strips, $\mathcal I_k^{+}$ and $\mathcal I_k^{-}$ defined by

\begin{equation}
\nonumber
\mathcal I_k^{+}:=\left\{I_{k,z}: z\ {\rm even}, 0 \leq z \leq \varrho L^{\beta}
 \right \} \quad{\rm and}\quad \mathcal I_k^{-}:=\left \{I_{k,z}:
z\ {\rm even}, -\varrho L^{\beta} \leq z<0 \right \}.
\end{equation}
Whenever there is no risk of confusion, we will drop the
subscripts from a strip $I_{k,z}$ writing $I$ instead.
We will need to work with the set of canonical directions
which are contained in the closed positive half-space
defined by the asymptotic direction, so that

$$
U^+:=\{e\in U: e\cdot\hat v\ge 0\}.
$$
 Let $s:=\left[2c_1\frac{L^{\zeta}}{\hat{v}\cdot e_1}\right]$.
For each $I \in \mathcal I_k^{+} \cup \mathcal I_k^{-}$ and each $y \in I$ we
associate a path $\hat{\pi}=\{\hat{\pi}_{1}, \ldots, \hat{\pi}_{ n}\}$,
with $s\le n\le s+1$, which satisfies

\begin{equation}
\nonumber
\hat{\pi} \subset I_{j,z}
\end{equation}
and

\begin{equation}
\nonumber
 \quad \hat{\pi}_{i+1}-\hat{\pi}_{i} \in U^+\quad{\rm for}\ 1\le i\le n-1, \quad \hat{\pi}_{n} \in H_{[L^{\zeta}]}.
\end{equation}
Note that by the fact that the strip $I$ has a Euclidean width
$1$, it is indeed possible to find a path satisfying these
 conditions and also that such a path is not necessarily unique. We will call $\hat \Pi_k$
 such a set of paths associated to all the points of the strip $I$.
 Now, a strip $I\in\mathcal I_k^+\cup\mathcal I_k^-$ will be called
 {\it elliptically good} if for all paths $\hat{\pi}\in\hat\Pi_k$ one has that

\begin{equation}
\label{hyp2a}
\sum_{i=1}^{n}\log \frac{1}{\omega\left(\hat{\pi}_{i}, \Delta \hat{\pi}_{i}\right)} \leq \log \left(\frac{1}{\xi}\right)n
\end{equation}
Otherwise $I$ will be called {\it elliptically bad} (See Figure \ref{coev4.1b}).

\begin{figure}[H]
\centering
\includegraphics[width=10cm]{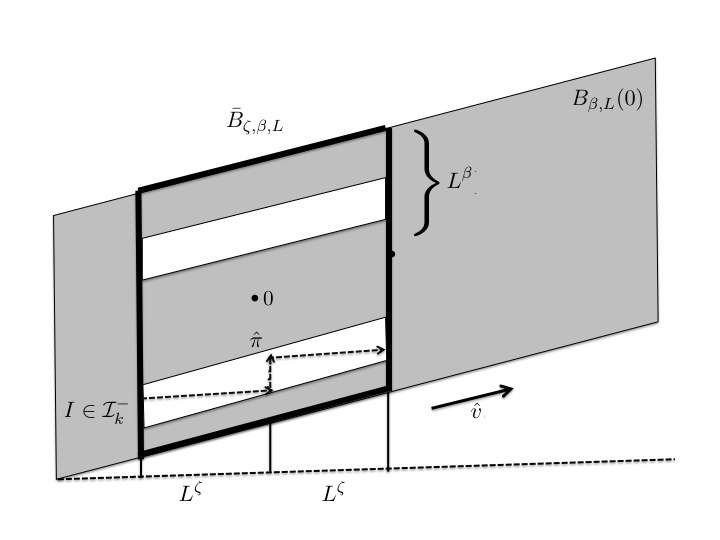}
\caption{In each strip $\mathcal I$, every path $\pi$ chosen previously satisfies the   ellipticity condition given by (\ref{hyp2a}). Then $\mathcal I$ is elliptically good.} 

\label{coev4.1b}
\end{figure}
%\noindent {\red Creo que el dibujo de $I$ est\'a correcto}

\noindent As before, from a routine counting argument and by
 Chebyshev inequality, note that for each $k\in\{2,\ldots,d\}$, $I \in \mathcal I_{k}^{+} \cup  \mathcal I_{k}^{-}$  and $\xi >0$ which satisfies $\log \frac{1}{\xi^{\alpha}\eta_{\alpha}^2}>0$, there exists a constant $c_{4,6}>0$ such that

\begin{equation}
\label{hyp2b}
\P\left( I \, \mbox{is} \, \mbox{elliptically} \, \mbox{bad}\right) \leq e^{-c_{4,6} L^{\zeta}}.
\end{equation}
We now define the  {\it constrainment event} as

\begin{eqnarray*}
\nonumber
&G_3:=\{\omega \in \Omega: \exists\, H_+\in{\mathcal H_+},
H_-\in{\mathcal H_-}, I_{+,2}\in{\mathcal I}^+_2,\ldots,
I_{+,d}\in{\mathcal I}^+_d,
I_{-,2}\in{\mathcal I}^-_2,\ldots,I_{-,d}\in{\mathcal I}^-_d\\
& {\rm such}\ {\rm that}\
H_+,H_-,I_{+,2},\ldots, I_{+,d}, I_{-,2},\ldots,I_{-,d}\ {\rm are}\ {\rm elliptically}\
  \mbox{good}\}.
\end{eqnarray*}
We can now state the following lemma which will eventually give
a control on the probability that the random walk hits an
elliptically good column.

\begin{lemma}
\label{lemma45}
There is a constant $c_{4,7}>0$ such that for every $L \geq 1$,

\begin{equation}
\label{cfe1}
\P \left(G_3^c \right) \leq e^{-c_{4,7}L^{\beta+\zeta}}.
\end{equation}

\end{lemma}

\begin{proof}
Note that
\begin{eqnarray}
\nonumber
&\P \left(G_3^c \right) \leq \P \left( \bigcap_{H \in \mathcal H^{+} }\{H \, \mbox{is} \, \mbox{elliptically} \, \mbox{bad} \} \right)  +\P \left( \bigcap_{H \in \mathcal H^{-} }\{H \, \mbox{is} \, \mbox{elliptically} \, \mbox{bad} \} \right)\\
\nonumber
&+ \sum_{k=2}^d\P\left(
 \bigcap_{I \in \mathcal I_k^{+}}\{I \, \mbox{is} \, \mbox{elliptically}\,
 \mbox{bad}\}\right)+
\sum_{k=2}^d \P\left(\bigcap_{I \in \mathcal I_k^{-}}\{I \, \mbox{is} \, \mbox{elliptically}\, \mbox{bad}\}\right),
\end{eqnarray}
Now, inequality (\ref{cfe1}) follows using the estimate  (\ref{hyp1b})
for the probability that a hyperplane is elliptically bad, the estimate
 (\ref{hyp2b}) for the probability that a strip is elliptically bad,
applying independence and  translation invariance.

\end{proof}

\subsection{Proof of Proposition \ref{crucial-lemma}}

Firstly, note that for any $\kappa>0$,
\begin{eqnarray}
\nonumber
&\P \left(P_{0,\omega} \left(X_{T_{B_{\beta,L}(0)}} \in \partial_{+}B_{\beta,L}(0) \right)\leq e^{-\kappa L^{\beta}}\right) \leq \P(G_1^c)+\P(G_2^c)+\P(G_3^c)\\
\label{plem41}
&+\P \left(P_{0,\omega} \left(X_{T_{B_{\beta,L}(0)}} \in \partial_{+}B_{\beta,L}(0) \right)\leq e^{-\kappa L^{\beta}}, G_1, G_2, G_3\right).
\end{eqnarray}
 Let us begin bounding the first three terms of
the right-hand side of
(\ref{plem41}).
 Let $\beta \in \left(\frac{\beta_0+1}{2},1 \right)$ and $\zeta\in (0,\beta_0)$. %$\beta \in \left(\frac{\beta_0+1}{2},1 \right)$.% 
By Lemma \ref{lemma43} of subsection \ref{bad-boxes}, Lemma \ref{lemma44}
of subsection \ref{bad-column} and Lemma \ref{lemma45}
of subsection \ref{sec-confin} we have that
there is a constant $c_{4,8}>0$ such that

\begin{equation}
\label{lem41aa}
\P(G_1^c)+\P(G_2^c)+\P(G_3^c) \leq
\frac{1}{c_{4,8}}e^{-c_{4,8}L^{3\beta-2+(d-1)(\beta-\beta_0)}}
+
\frac{1}{c_{4,8}}e^{-c_{4,8}L^{\beta+(d-1)(\beta-\beta_0)}}
+
\frac{1}{c_{4,8}}e^{-c_{4,8}L^{\beta+\zeta}}.
\end{equation}
Since $\beta<1$ is equivalent to
$\beta+(d-1)(\beta-\beta_0)>3\beta-2+(d-1)(\beta-\beta_0)$,
the sum in (\ref{lem41aa}) can be bounded as

\begin{equation}
\label{lem41a}
\P(G_1^c)+\P(G_2^c)+\P(G_3^c) \leq
\frac{1}{c_{4,9}}e^{-c_{4,9}L^{g(\beta,\beta_0,\zeta)}},
\end{equation}
for some constant $c_{4,9}>0$ and where
 $g(\beta, \beta_0, \zeta):=\min \{\beta+\zeta,
 3\beta-2+(d-1)(\beta-\beta_0)\}$.

We will now prove that the fourth term of
the right-hand side of inequality
(\ref{plem41}) satisfies for $L$ large enough

\begin{equation}
\label{e0}
\P \left(P_{0,\omega} \left(X_{T_{B_{\beta,L}(0)}} \in
\partial_{+}B_{\beta,L}(0) \right)\leq e^{-\kappa L^{\beta}}, G_1, G_2, G_3\right)=0.
\end{equation}
 In fact, we will show that for $L$ large enough
on the event $G_1\cap G_2\cap G_3$
one has that

\begin{equation}
\label{g123}
P_{0,\omega} \left(X_{T_{B_{\beta,L}(0)}} \in \partial_{+}B_{\beta,L}(0)
 \right)> e^{-\kappa L^{\beta}}.
\end{equation}
We will prove (\ref{g123})
showing that the walk can exit $B_{\beta, L}(0)$ through
 $\partial_{+}B_{\beta, L}(0)$ choosing a strategy which
corresponds to paths which go through an  elliptically good column. This implies, in particular, that the walk exits successively boxes $\tilde{B}(x)$ through $\partial_{+}\tilde{B}(x)$.
 The event $G_1$ implies that there exist at least
 $m_3=\left[\frac{\varrho^{d-1}
     L^{(d-1)(\beta-\beta_0)}}{2(1+\varrho)^{d-1}}\right]$ columns each one
 with at most $\left[\frac{L^{\beta}}{L_0}\right]$ of bad boxes. Meanwhile,
 the event $G_2$ asserts that in any collection of columns with cardinality
 $m_3$ or more, there is at least one elliptically good column.
Therefore, on the event $G_1\cap G_2$ there exists at least one
elliptically good column $D$ with at most $L^\beta/L_0$ bad boxes.
Thus, on $G_1\cap G_2$ we have that
for any point $y\in D$ and $\xi>0$,

\begin{equation}
\label{lem41b}
P_{y, \omega}\left(X_{T_{B_{\beta,L}(0)}} \in \partial_{+}B_{\beta,L}(0)
\right) \geq \left(\frac{1}{2}\right)^{L^{\beta-\beta_0}+1} \xi^
      {2c_1\frac{\varrho}{\hat v\cdot e_1}L^{\beta}}
\xi^{2c_1\varrho L^\beta},
\end{equation}
where the first factor is a bound for the probability that the random
walk exits all the good boxes of the column through their front side, 
%{\red quiz\'as basta que el camino se quede en $\tilde B$. Ya lo hab\'{\i}a comentado antes}
while the second factor is a bound for the probability that the
walk traverses each bad box (whose number is at most
$L^\beta/L_0$) exiting through its front side and
following a path with at most $\frac{2c_1\rho L_0}{\hat v\cdot e_1}$ steps
and is given by the condition (\ref{egd1}) for elliptically good columns,
while the third factor is a bound for the probability
that once the walk exits a box (whose number is at most
$L^{\beta-\beta_0}+1$) it moves through its front boundary to
the central point of this front boundary following a
path with at most $[2c_1\rho L^\beta_0]$ steps and is given by the condition
(\ref{egd2}) for elliptically good columns.

Now, the constrainment event $G_3$ ensures that
with a high enough probability the random walk will reach
the elliptically good column $D$ which has at most $L^\beta/L_0$ bad boxes.
 More precisely, a.s. on $G_3$,  the random walk reaches either
an elliptically good hyperplane $H\in\mathcal H_+\cup\mathcal H_-$,
an elliptically good strip $I\in \mathcal I^+_2\cup\cdots\cup\mathcal I^+_d$
or an elliptically good
strip $I\in \mathcal I^-_2\cup\cdots\cup\mathcal I^-_d$
(recall the definitions of elliptically good hyperplanes and strips
given in (\ref{hyp1a}) and (\ref{hyp2a}) of subsection \ref{bad-column}).
Now, once the walk reaches either an elliptically
good hyperplane or strip, we know  by (\ref{hyp1a}) or (\ref{hyp2a}),
choosing an appropriate path that the probability that it hits
the column $D$ is at least $\xi^{c_{4,10}\varrho L^\beta}$
for some constant $c_{4,10}>0$.
Thus, we know that there is a constant $c_{4,10}>0$ such that
\begin{equation}
\label{lem41c}
P_{0,\omega}\left(\mbox{the} \, \mbox{walk} \, \mbox {reaches} \, D \cap \bar B_{\zeta, \beta, L} (0) \right) \geq \xi^{c_{4,10}\varrho L^{\beta}}.
\end{equation}
Therefore, combining (\ref{lem41b}) and (\ref{lem41c}),
 we conclude that there is a constant $c_{4,11}>0$ such that
for all $\varrho\in (0,1)$ on the event $G_1\cap G_2\cap G_3$ the following estimate is satisfied,

$$
P_{0,\omega} \left(X_{T_{B_{\beta,L}(0)}} \in \partial_{+}B_{\beta,L}(0)
 \right)> e^{-c_{4,11}\rho L^{\beta}}.
$$
Hence, choosing $\varrho$ sufficiently small, we have that
on $G_1\cap G_2\cap G_3$,

\begin{equation}
\label{plem41a}
P_{0, \omega}\left(X_{T_{B_{\beta,L}(0)}} \in \partial_{+}B_{\beta,L}(0)
\right) >e^{-\kappa L^{\beta}}
\end{equation}
for $L$ larger than a deterministic constant depending only on $\varrho$.
This proves (\ref{e0}).

Finally, with the help of (\ref{plem41}), (\ref{lem41a}) and (\ref{plem41a}) the Proposition \ref{crucial-lemma} is proved.

\section{Moments of the regeneration time}
\label{tail-regen}
\setcounter{equation}{0}

Here we will prove Theorems \ref{theorem2} and
\ref{theorem3}.  Our method is inspired
by some ideas used by Sznitman  to prove
Proposition 3.1 of \cite{Sz01}, which give tail estimates
on the distribution of the regeneration times.
 Theorem \ref{theorem2}
will follow from part $(a)$ of  Theorem \ref{lwlln}, while
part $(a)$ of Theorem \ref{theorem3} from
part $(b)$ of Theorem \ref{lwlln}. Part $(b)$
of Theorem \ref{theorem3} will follow from Theorem \ref{lwlln-quenched}.
The non-degeneracy of the
covariance matrix in the annealed part $(a)$ of Theorem \ref{theorem3}
can be proven exactly as in section IV of \cite{Sz00}
using (\ref{defa}).
\medskip

\begin{proposition}
\label{tail prop}
Let $l\in\mathbb S^{d-1}$, $\alpha>0$ and $M\ge 15d+5$.
Assume that $(P)_M|l$ is satisfied and that $(E')_\alpha$ holds
towards the asymptotic direction (cf. (\ref{kappa}), (\ref{finiteness}),
(\ref{direction})).
Then

\begin{equation}
\nonumber
\limsup_{u\to \infty} (\log u)^{-1}\log P_0[\tau_1^{\hat v}>u]\leq -\alpha.
\end{equation}

\end{proposition}

\medskip

\noindent The proof of the above proposition
is based on the atypical quenched exit estimate
corresponding to Proposition \ref{crucial-lemma}  of section \ref{section-four}.
Some slight modifications in the proof of
Proposition \ref{crucial-lemma}, would lead to a version of
it, which could be used to show that
Proposition \ref{tail prop}
remains valid if the regeneration time $\tau_1^{\hat v}$
is replaced by $\tau_1^l$ for any direction $l$ such that
$l\cdot\hat v>0$. Note also that Proposition \ref{tail prop} implies that
whenever  $(E')_1$ is satisfied towards the asymptotic direction,  then
the first regeneration time is integrable. Through
part $(a)$ of
Theorem \ref{lwlln}, this implies
 Theorem \ref{theorem2}.
Similarly we can conclude part $(a)$ of Theorem \ref{theorem3}.
Part $(b)$ of Theorem \ref{theorem3} can be derived analogously
through Theorem \ref{lwlln-quenched}.

Let us now proceed with the proof of Proposition \ref{tail prop}.
As in section \ref{section-four}, we will assume that $\hat v\cdot e_1>0$ (c.f.
(\ref{secondwlg})).
Let us take a rotation $\hat R$ in $\R^d$ such that $\hat R(e_1)=\hat v$
and fix
$\beta\in \left(\frac{5}{6},1\right)$ and $M>0$.
For each $u>0$  define the scale

$$
L=L(u):=\left(\frac{1}{4M\sqrt{d}
}\right)^{\frac{1}{\beta}}(\log u)^{\frac{1}{\beta}},
$$
and the box

$$
C_L:=\left\{x\in\mathbb Z^d:\frac{-L}{2(\hat v \cdot e_1)}
\le x\cdot \hat R(e_i)\le \frac{L}{2(\hat v \cdot
  e_1)}\ {\rm for}\ 1\le i\le 2d\right\}.
$$

%\begin{equation}
%\label{fundamental condition}
%\nu_{\alpha}:=2\alpha \left(\frac{M-1}{M}\right)-\frac{4\log \eta_{\alpha}}{M (\hat v \cdot e_1)\log \left(\frac{1}{\xi}\right)}>1
%\end{equation}

\noindent Throughout the rest of this proof we will
continue writing $\tau_1$ instead of $\tau_1^{\hat v}$.
Now note that
% {\red no entiendo por qu\'e deber\'{\i}a cambiar a continuaci\'\ on el $0$ por $1$?}

\begin{equation}
\label{tau-t}
P_0(\tau_1 >u)\leq P_0\left(\tau_1>u,T_{C_{L(u)}}\leq
\tau_1\right)+P_0\left(T_{C_{L(u)}}>u\right),
\end{equation}
where $T_{C_{L(u)}}$ is the first exit time from
the set $C_{L(u)}$ defined in (\ref{exit-time}).
 For the first term of
the right-hand side of inequality (\ref{tau-t}), we
can use Corollary \ref{regeneration time lemma},
 to conclude that for every $\gamma \in (\beta,1)$ there exists a constant
$c_{5,1}$ such that
\begin{equation}
\label{tail estimate 1}
P_0 \left(\tau_1>u,T_{C_{L(u)}}\leq \tau_1 \right) \leq
\frac{1}{c_{5,1}}e^{-c_{5,1}L^{\gamma}(u)}.
\end{equation}
For the second term of the right-hand side of
inequality (\ref{tau-t}), following Sznitman \cite{Sz01}
 we introduce the event

\begin{equation}
\nonumber
F_1:=\left\{\omega \in \Omega: t_{\omega}\left(C_{L(u)}\right)> \frac{u}{(\log u)^{\frac{1}{\beta}}}\right\},
\end{equation}
where for each $A\subset\mathbb Z^d$ we define

$$
t_{\omega}(A):=\inf \left\{n \geq 0: \sup_{x}P_{x,
  \omega}\left(T_A>n\right)\leq \frac{1}{2}\right\}.
$$
Trivially,

\begin{equation}
\label{trivial}
P_0 \left(
T_{C_{L(u)}}>u \right) \leq \E \left[F_1^c, P_{0,\omega} \left(T_{C_{L(u)}}>u
  \right)\right]+\P(F_1).
\end{equation}

\noindent To bound the first term of the right-hand side
of (\ref{trivial}), on the event $F_1^c$ we apply the strong Markov property
$[(\log u)^{\frac{1}{\beta}}]$ times to conclude that

\begin{equation}
\label{tail estimate 2}
\E \left[F_1^c, P_{0,\omega} \left(T_{C_{L(u)}}>u \right)\right] \leq \left(\frac{1}{2}\right)^{[(\log u)^{\frac{1}{\beta}}]}.
\end{equation}
\noindent To bound the second term of the right-hand side
of (\ref{trivial}), we will use the fact that
 for each $\omega \in \Omega$ there exists $x_0 \in C_{L(u)}$ such that

\begin{equation}
\label{hitting time estimation}
P_{x_0,\omega}(\tilde{H}_{x_0}>T_{C_{L(u)}}) \leq \frac{2|C_{L(u)}|}{t_{\omega}(C_{L(u)})}
\end{equation}
where for $y\in\mathbb Z^d$,

$$
\tilde{H}_y=\inf\{n \geq 1: X_n=y\}.
$$
(\ref{hitting time estimation})  can be derived
using the fact that for every subset $A\subset \mathbb Z^d$
and $x\in A$,

$$
E_{x,\omega}(T_{A})=\sum_{y\in A}\frac{P_{x,\omega}(H_y<T_A)}
{P_{y,\omega}(\tilde H_y>T_A)}
$$
(see for example Lemma 1.3 of Sznitman \cite{Sz01}).
Now note that (\ref{hitting time estimation}) implies

\begin{equation}
\label{tail-s}
\P(F_1) \leq \P\left(\omega \in \Omega: \exists \,\, x_0 \in C_{L(u)} \,\,{\rm such} \,\, {\rm that}\,\, P_{x_0,\omega}(\tilde H_{x_0}>T_{C_{L(u)}})\leq \frac{2(\log u)^{\frac{1}{\beta}}}{u}|C_{L(u)}|\right).
\end{equation}
 Choose for each $x \in C_{L(u)}$ a point $y_x$ as any
  point in $\Z^d$ which is closest to
the point
$x+\left(\frac{\log u}{2M\sqrt{d}(\hat{v} \cdot e_1)}\right)\,\, \hat{v}=x+2\left(\frac{L^\beta(u)}{\hat{v} \cdot e_1}\right)\,\, \hat{v}
$.
 It is straightforward to see that
%for $u$ large enough

\begin{equation}
\nonumber
 N-1 \leq |y_x-x|_1 \leq
 N+1,
\end{equation}
% {\red creo que esto se cumple a\'un con la norma $|\ldots|_1$}
where

$$
N:= \frac{ |\hat v|_1\log u}{2M \sqrt{d}(\hat v \cdot e_1)}.
$$
Let us call

$$
V(y_x):=\{y\in\mathbb Z^d:|y_x-y|_1\le 4d\},
$$
the closed $l_1$ ball centered at $y_x$. We furthermore define

\begin{align}
\label{defvprime}
V'(y_x):=
\left\{
 \begin{array}{r}
V(y_x)\cap C_{L(u)}\quad {\rm\ if\ } y_x\in C_{L(u)}\\
V(y_x)\cap (C_{L(u)})^c\quad  {\rm\ if\ } y_x\notin C_{L(u)}.
\end{array} \right.
\end{align}
Now, there are constants $K_1$ and $K_2$ such that for each $1\le i\le 2d$, we can  find $4d-2$ different paths
 $\{\pi^{(i,j)}:1\le j\le 4d-2\}$,
each one with $n_j$ steps,
with $\pi^{(i,j)}:=\{\pi_1^{(i,j)}, \ldots, \pi_{n_j}^{(i,j)}\}$ for each $1\le j\le
 4d-2$
such that the following conditions are satisfied (recall that $U$
are the canonical vectors in $\mathbb Z^d$):

\begin{enumerate}

\item {\bf (the paths connect $x$ or $x+e_i$ to $V'(y_x)$)} For each $1\le j\le 2d-1$, the path $\pi^{(i,j)}$ goes from $x$ to
 $V'(y_x)$, so that
$\pi^{(i,j)}_1=x$ and $\pi^{(i,j)}_{n_i}\in V'(y_x)$.
For each $2d\le j\le 4d-2$, the path $\pi^{(i,j)}$ goes from
$x+e_i$ to $V'(y_x)$, so that
$\pi^{(i,j)}_1=x+e_i$ and $\pi^{(i,j)}_{n_j}\in V'(y_x)$.

\item {\bf (the paths start using $4d-2$ directions as shown
in Figure \ref{paths5.1})} For each $1\le j\le 2d-1$,
$\pi^{(i,j)}_2-\pi_1^{(i,j)}=f_j$ where $f_j\in U-\{e_i\}$
and $f_j\ne f_{j'}$ for $1\le j\ne j'\le 2d-1$.
For each $2d\le j\le 4d-2$,
$\pi^{(i,j)}_2-\pi_1^{(i,j)}=g_j$ where $g_j\in U-\{-e_i\}$
and $g_j\ne g_{j'}$ for $2d\le j\ne j'\le 4d-2$.

\item {\bf (the paths intersect at most $K_1$ times)} For each $1\le j\le 4d-2$, the
paths $\pi^{(i,j)}$ and $\pi^{(i,j')}$ intersect at at most
$K_1$ vertexes. Furthermore, after each point of intersection, both
paths perform jumps in different directions in $H_{\hat v}$.

\item  {\bf (the number of steps of all paths is close to $N$)} The number of steps $n_j$ of each path satisfies $N-1\le n_j\le N+K_2$.

\item  {\bf (increments in $H_{\hat v}$)} For each $1\le j\le 2d-1$,  we have that
$\Delta\pi^{(i,j)}_k\in H_{\hat v}\cup \{f_j\}$ for $2\le k\le n_j-1$,
while for each $2d\le j\le 4d-2$,  we have that
$\Delta\pi^{(i,j)}_k\in H_{\hat v}\cup \{g_j\}$ for $2\le k\le n_j-1$.
\end{enumerate}

\begin{figure}[H]
\centering
\includegraphics[width=10cm]{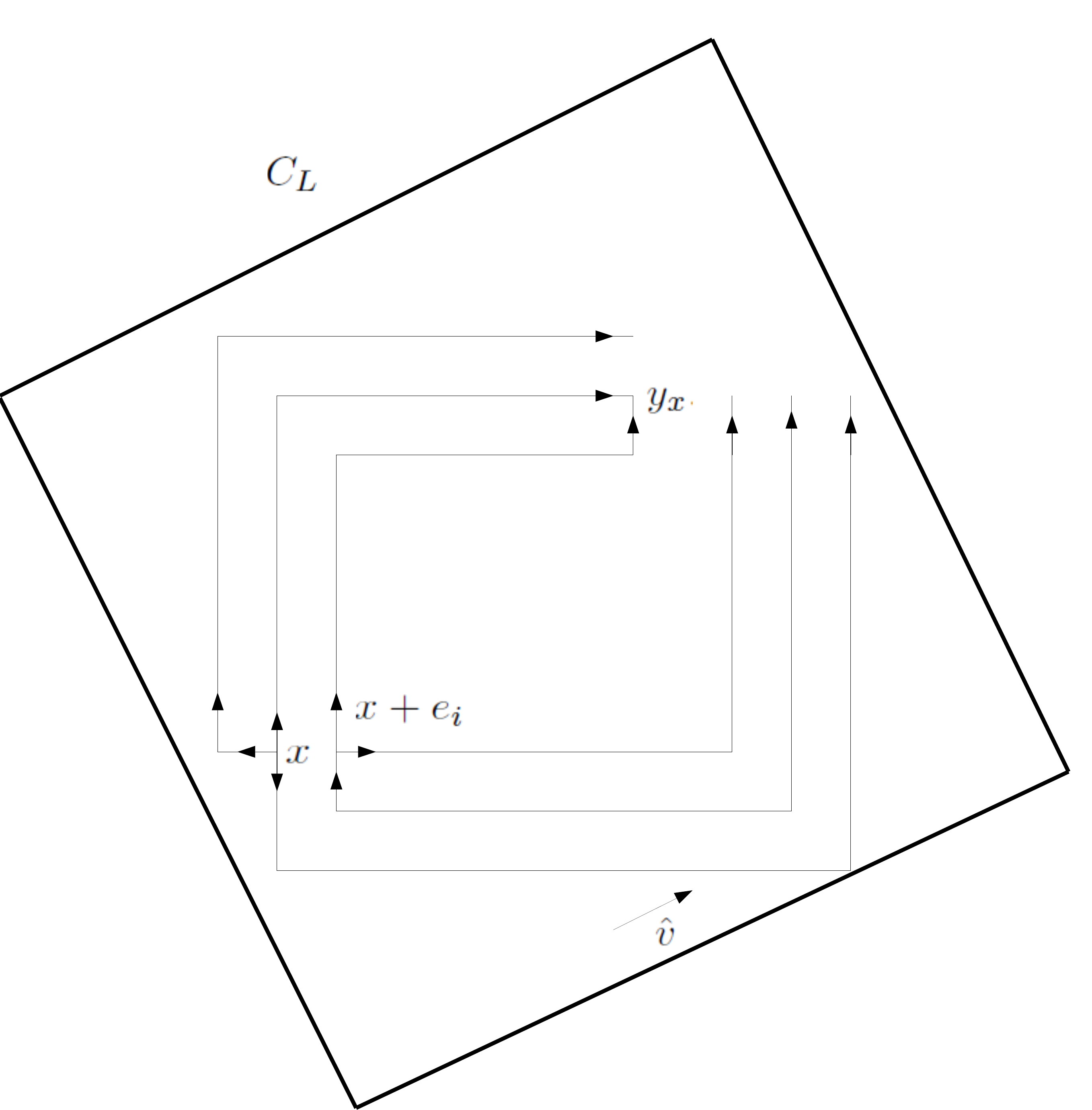}
\caption{The $4\times 2-2=6$ paths from $x$ to $V'(y_x)$ and from $x+e_i$ to $V'(y_x)$ are represented by  the arrowed lines. The set $V'(y_x)$ is given by the $5$
endpoints of the paths.}
\label{paths5.1}
\end{figure}

\noindent In Figure \ref{paths5.1} it is seen how can one construct such a set
of paths for dimensions $d=2$ (a similar construction works for
dimensions $d\ge 3$). From Figure \ref{paths5.1}, note that
the maximal number of steps of each path is given by
$|y_x-x|_1+7$, where the $7$ corresponds to the extra steps which have to be
performed when a path exits the point $x$ (or $x+e_i$)
or enters some point in $V'(y_x)$.
 Let us now introduce for each $1\le i\le 2d$ the  event

\begin{eqnarray*}
&F_{2,i}:= \left\{\omega \in \Omega: \mbox{for}\,\,
 \mbox{each} \,\, x \in C_{L(u)},\,\, \exists \, j \in \{1, \ldots, 2d\}\, \, \mbox{such} \, \, \mbox{that}
 \right.\\
& \left. \sum_{k=1}^{n_j}\log \frac{1}{\omega(\pi^{(i,j)}_k,\Delta \pi^{(i,j)}_k)} \leq \frac{2(M-1)(\hat v \cdot e_1)\sqrt{d}}{|\hat v|_1} n_j \right \}
\end{eqnarray*}
and also

$$
F_{2}:=\cap_{i=1}^{2d}F_{2,i}.
$$
Then, with the help of (\ref{tail-s})
 we have that

\begin{equation}
\label{tail estimate 3}
\P(F_1) \leq \P\left(\exists\,\, x_0 \in C_{L(u)} \,\,{\rm such}\,\,{\rm that}\ P_{x_0,\omega}(\tilde{H}_{x_0}>T_{C_{L(u)}})\leq \frac{2(\log u)^{\frac{1}{\beta}}}{u}|C_{L(u)}|,F_2 \right)+\P(F_2^c).
\end{equation}

\noindent Let us define

$$
F_3:=\left\{\omega \in \Omega:
\exists\,\, x_0 \in C_{L(u)} \,\,{\rm such}\,\,{\rm that}\
P_{x_0,\omega}(\tilde{H}_{x_0}>T_{C_{L(u)}})\leq \frac{2(\log
  u)^{\frac{1}{\beta}}}{u}|C_{L(u)}|,F_2 \right\}.
$$
Note that on the event $F_3$, which appears in the
probability of the right-hand side of
(\ref{tail estimate 3}), we can use the definition of the event $F_2$ to join
for each $1\le i\le 2d$,
$x_0$ with $V'(y_{x_0})$  using one of the paths $\pi^{(i,j)}$ to conclude that

\begin{equation}
\label{ibound}
\omega_i e^{-\frac{8(M-1)(\hat v \cdot e_1)\sqrt{d}}{|\hat v|_1}} u^{-\left(1-\frac{1}{M}\right)}\inf_{z\in V'(y_{x_0})}
P_{z,\omega} \left(T_{C_{L(u)}}<H_{x_0}\right) \leq  P_{x_0,\omega}\left(T_{C_{L(u)}}<\tilde{H}_{x_0}\right)
 \leq \frac{2(\log u)^{\frac{1}{\beta}}}{u}|C_{L(u)}|,
\end{equation}
where

$$
\omega_i:=\omega(x_0,e_i).
$$
The factor $\omega_i$ above corresponds to the probability of jumping
from $x_0$ to $x_0+e_i$ (in the case that the path $\pi^{(i,j)}$ starts from
$x_0+e_i$).
Summing up over $i$  in (\ref{ibound}) and using the equality

\begin{equation}
\label{keysum}
\sum_{i=1}^{2d}\omega_i=1,
\end{equation}
 we conclude that on the event $F_3$
one has that

$$
e^{-\frac{8(M-1)(\hat v \cdot e_1)\sqrt{d}}{|\hat v|_1}} u^{-\left(1-\frac{1}{M}\right)}
\inf_{z\in V'(y_{x_0})} P_{z,\omega} \left(T_{C_{L(u)}}<H_{x}\right)
 \leq 4d\frac{(\log u)^{\frac{1}{\beta}}}{u}|C_{L(u)}|.
$$
 In particular, on $F_3$ we can see that for $u$ large enough
 $ V'(y_{x_0}) \subset C_{L(u)}$ (see (\ref{defvprime})).  As a result, on $F_3$ we have that
for $u$ large enough

\begin{equation}
\nonumber
\inf_{z\in V'(y_{x_0})} P_{z,\omega}\left(X_{T_{z+U_{\beta,L}}}\cdot e_1>z\cdot e_1\right)  \leq
\inf_{z\in V'(y_{x_0})} P_{z,\omega} \left(T_{C_{L(u)}}<H_{x}\right)
 \leq  \frac{1}{u^{\frac{1}{2M}}}=e^{-2\sqrt{d} L(u)^{\beta}},
\end{equation}
where

$$
U_{\beta,L}:=\{x \in \Z^d: -L^{\beta} <x \cdot e_1<L\}.
$$

\noindent From this and using the translation invariance of the measure $\P$, we conclude that

\begin{eqnarray}
\nonumber
&\P\left(\exists\,\, x_0 \in C_{L(u)} \,\,{\rm such}\,\, {\rm that} \,\, P_{x_0,\omega}[\tilde{H}_{x_0}>T_{C_{L(u)}}]\leq 4d\frac{(\log u)^{\frac{1}{\beta}}}{u}|C_{L(u)}|,F_2 \right)\\
\nonumber
&\le\P\left(\exists\,\, x_0 \in C_{L(u)} \,\,{\rm such}\,\, {\rm that} \,\,
\inf_{z\in V'(y_{x_0})} P_{z,\omega}(X_{T_{z+U_{\beta,L(u)}}}\cdot  e_1>z\cdot e_1)\leq e^{-2\sqrt{d} L(u)^{\beta}}\right)\\
\nonumber
&\leq  |V'(0)|
|C_{L(u)}|\P \left(P_{0,\omega}\left(X_{T_{U_{\beta,L(u)}}}\cdot  e_1
>0 \right)\leq e^{-2\sqrt{d} L(u)^{\beta}}\right)\\
\nonumber
&\leq |V'(0)| |C_{L(u)}|\P \left(P_{0,\omega}\left(X_{T_{B_{\beta,L(u)}}}\cdot  e_1 >0 \right)\leq e^{-2\sqrt{d} L(u)^{\beta}}\right),
\end{eqnarray}
where the titled box $B_{\beta,L}$ was defined in (\ref{tilted-box})
of section \ref{section-four}. Therefore, we can estimate the first term of the right-hand side
of (\ref{tail estimate 3}) using Proposition {\ref{crucial-lemma}}
to conclude that there is a constant $c_{5,2} >0$
such that
for each $\beta_0\in \left(\frac{1}{2},1\right)$ one has that

\begin{equation}
\label{tail estimate 4}
\P\left(\exists\,\, x_0 \in C_{L(u)} \,\,{\rm such}\,\, {\rm that} \,\, P_{x_0,\omega}[\tilde{H}_{x_0}>T_U]\leq \frac{2(\log u)^{\frac{1}{\beta}}}{u}|C_{L(u)}|,F_2 \right) \leq \frac{1}{c_{5,2}} e^{-c_{5,2}L(u)^{g(\beta_0,\beta,  \zeta)}},
\end{equation}
where $g(\beta_0,\beta,\zeta)$ is defined in (\ref{def-g})
of Proposition \ref{crucial-lemma}.
 We next have to control the probability $\P(F_2^c)$.
To simplify the notation in the calculations that
follow, we drop the super-index $i$ in the paths writing
$\pi^{(j)}:=\pi^{(i,j)}$. Furthermore we will define

$$
M_1:=2(M-1)\frac{(\hat v\cdot e_1)\sqrt{d}}{|\hat v|_1}.
$$
 Then,
 by the independence property $(e)$ of the paths $\{\pi^{(i,j)}\}$, we can say that

\begin{eqnarray}
\nonumber
&
\P(F_2^c)\le\\
\nonumber
& \sum_{i=1}^{2d}\P\left(\exists x\in C_{L(u)}\ {\rm such}\ {\rm that}\
\forall\ j,\
\log\frac{1}{\omega(\pi_1^{(j)},\Delta \pi_1^{(j)})}
+\sum_{k=2}^{n_j}\log\frac{1}{\omega(\pi_k^{(j)},\Delta \pi_k^{(j)})}
> M_1n_j
\right)\\
\nonumber
&\le \sum_{i=1}^{2d}\sum_{x\in C_{L(u)}}\P\left(
\forall\ 1\le j\le 2d-1,\
\log\frac{1}{\omega(\pi_1^{(j)},\Delta \pi_1^{(j)})}
+\sum_{k=2}^{n_j}\log\frac{1}{\omega(\pi_k^{(j)},\Delta \pi_k^{(j)})}
>M_1n_j
\right)\\
\label{impr1}
&\times \P\left(
\forall\ 2d\le j\le 4d-2,\
\log\frac{1}{\omega(\pi_1^{(j)},\Delta \pi_1^{(j)})}
+\sum_{k=2}^{n_j}\log\frac{1}{\omega(\pi_k^{(j)},\Delta \pi_k^{(j)})}
>M_1n_j
\right).
\end{eqnarray}
 Now, using
Chebyshev's inequality, the first step property $(b)$ of
the paths, the intersection property $(c)$,  the
bound on the number of steps $(d)$, $N-1\le n_j$,
and the property of the increments $(e)$, we
can bound  the rightmost-hand side of (\ref{impr1}) by the following expression,
where for $1\le j\le 2d-1$, $\alpha_j:=\alpha_1$
 for $f_j\in H_{\hat v}$ while
$\alpha_j:=\alpha(f_j)\le\alpha_1$ for $f_j\notin H_{\hat v}$
(cf. (\ref{direction}), and similarly
for $2d\le j\le 4d-2$.

\begin{eqnarray*}
&|C_{L(u)}|\sum_{i=1}^{2d}\E\left[e^{\sum_{j=1}^{2d-1} \alpha_j\log\frac{1}{\omega(0,f_j)}}\right]^{2K_1}
e^{-M_1\sum_{j=1}^{2d-1}\alpha_jn_j}
\Pi_{j=1}^{2d-1} \E\left[e^{\alpha_j\sum_{k=2}^{n_j}\log\frac{1}{\omega\left(\pi_k^{(j)},
\Delta\pi_k^{(j)}\right)}}\right]
\\
&\times
\E\left[e^{\sum_{j=2d}^{4d-2}\alpha_j \log\frac{1}{\omega(0,g_j)}}\right]^{2K_1}
e^{-M_1\sum_{j=2d}^{4d-2}\alpha_jn_j}
\Pi_{j=2d}^{4d-2} \E\left[e^{\alpha_j\sum_{k=2}^{n_j}\log\frac{1}{\omega\left(\pi_k^{(j)},
\Delta\pi_k^{(j)}\right)}}\right]\\
&\le
|C_{L(u)}|\sum_{i=1}^{2d}(\eta_i^+\eta_i^-)^{2K_1}e^{(\log\eta_\alpha)\sum_{j=1}^{4d-2}
  n_j
+
}
e^{-\kappa
\left(1-\frac{1}{M}\right)\log u
},
\end{eqnarray*}
where $\kappa$ is defined in (\ref{kappa}), $\eta_\alpha$ is defined
in (\ref{eta-alpha}) for $\alpha<\min_{e\in U}\alpha(e)$ (cf. (\ref{eprime-eo})) and where

$$
\eta^+_i:=\E\left[e^{\sum_{j=2d}^{4d-2} \alpha_j\log\frac{1}{\omega(0,g_j)}}\right]
$$
and

$$
\eta^-_i:=\E\left[e^{\sum_{j=1}^{2d-1}\alpha_j \log\frac{1}{\omega(0,g_j)}}\right].
$$

\noindent Then, using the inequality $n_j\le N+K_2$ (property $(g)$ of the paths), we have that

\begin{equation}
\nonumber
\P \left(F_2^c\right) \leq
|C_{L(u)}|\sum_{i,i'=1}^{2d}(\eta_i^+
\eta_i^-)^{2K_1}
 e^{(4d-2)(N+K_2)\log\eta_\alpha}
 e^{-\kappa \left(1-\frac{1}{M}\right)\log u}.
\end{equation}
Now, for $\alpha<\kappa$, 
using the definition of $N$, and of condition $(E')_\alpha$
(cf. (\ref{kappa}) and (\ref{finiteness})) we see from here
that  if we choose
$M$ large enough,
one has that for $u$ large enough

\begin{equation}
\label{tail estimate 5}
\P(F_2^c)\le c_{5,3} u^{-\alpha}.
\end{equation}
\noindent for some constant $c_{5,3}>0$.
Now note that for each $\beta\in \left(\frac{5}{6},1\right)$ there exists a
$\beta_0 \in \left(\frac{1}{2},\beta\right)$ such that
for every $\zeta\in \left(0,\frac{1}{2}\right)$ one has that

\begin{equation}
\label{gbb}
g(\beta,\beta_0,\zeta)>\beta.
\end{equation}
 Therefore, substituting (\ref{tail estimate 4}) and (\ref{tail estimate 5})
 back into (\ref{tail estimate 3}) and using
(\ref{gbb}) we can see that there is a constant $c_{5,4}>0$ such that for $u$
 large enough

\begin{equation}
\label{tail est 6}
\P(F_1) \leq c_{5,4} u^{-\alpha}.
\end{equation}

\noindent Now with the help of (\ref{trivial}), (\ref{tail estimate 2}) and (\ref{tail est 6}) there exists a constant $c_{5,5}>0$ such that for $u$ large

\begin{equation}
\nonumber
P_0 \left(T_{C_{L(u)}} >u \right) \leq c_{5,5} u^{-\alpha}.
\end{equation}
Finally, since $\gamma\in (\beta,1)$ in (\ref{tail estimate 1}), using
(\ref{tau-t}) we conclude the proof,  since we see that for $u$ large enough

\begin{equation}
\nonumber
P_0(\tau_1 >u) \leq c_{5,6} u^{-\alpha},
\end{equation}
for a certain constant $c_{5,6} >0$. This proves part  Proposition
\ref{tail prop}.

\medskip

{\bf Acknowledgments.} The authors would like
to thank Christophe Sabot for pointing out that
equality (\ref{keysum}) gives the right integrability
condition for RWDRE, to Alexander Drewitz, Alexander Fribergh and Laurent Tournier
for useful discussions, and to Elodie Bouchet and the referee
for several corrections.

\end{document}